\documentclass[11pt,reqno]{amsart}
\usepackage{amssymb,amscd,amsbsy,mathrsfs}
\newcommand{\lb}{\linebreak}

\renewcommand{\a}{\alpha}
\renewcommand{\b}{\beta}
\newcommand{\g}{\gamma}

\newcommand{\e}{\varepsilon}

\newcommand{\z}{\zeta}

\renewcommand{\l}{\lambda}

\newcommand{\s}{\sigma}
\renewcommand{\t}{\tau}

\newcommand{\f}{\varphi}
\renewcommand{\o}{\omega}

\newcommand{\D}{\Delta}
\renewcommand{\L}{\Lambda}

\renewcommand{\O}{\Omega}

\newcommand{\B}{{\mathscr B}}

\newcommand{\E}{{\mathscr E}}

\newcommand{\F}{{\mathscr F}}

\newcommand{\h}{{\mathscr H}}

\newcommand{\cL}{{\mathscr L}}
\newcommand{\M}{{\mathscr M}}

\newcommand{\cR}{{\mathscr R}}

\newcommand{\X}{{\mathscr X}}
\newcommand{\Y}{{\mathscr Y}}

\newcommand{\W}{{\mathscr W}}

\newcommand{\C}{{\Bbb C}}
\newcommand{\T}{{\Bbb T}}
\newcommand{\pp}{{\Bbb P}}
\newcommand{\dd}{{\Bbb D}}
\newcommand{\R}{{\Bbb R}}
\newcommand{\Z}{{\Bbb Z}}

\newcommand{\0}{{\boldsymbol{0}}}

\newcommand{\bs}{\boldsymbol}

\newcommand{\m}{{\boldsymbol m}}
\newcommand{\bS}{{\boldsymbol S}}

\newcommand{\rf}[1]{(\ref{#1})}

\newcommand{\df}{\stackrel{\mathrm{def}}{=}}

\newcommand{\re}{\operatorname{Re}}

\newcommand{\supp}{\operatorname{supp}}

\newcommand{\trace}{\operatorname{trace}}
\newcommand{\rank}{\operatorname{rank}}
\newcommand{\const}{\operatorname{const}}
\newcommand{\tr}{\operatorname{trace}}
\newcommand{\eeq}{\end{equation}}
\newcommand{\beq}{\begin{equation}}
\newcommand{\bay}{\begin{eqnarray}}
\newcommand{\ba}{\begin{align*}}
\newcommand{\ea}{\end{align*}}
\newcommand{\ey}{\end{eqnarray}}
\newcommand{\bey}{\begin{eqnarray*}}
\newcommand{\eey}{\end{eqnarray*}}

\newcommand{\eq}{\Leftrightarrow}
\newcommand{\imp}{\Rightarrow}
\newcommand{\be}{\infty}

\newcommand{\bl}{\blacksquare}
\newcommand{\ess}{\operatorname{ess}}

\newcommand{\Range}{\operatorname{Range}}

\newcommand{\Pf}{{\bf Proof. }}
\newcommand{\im}{\operatorname{Im}}
\renewcommand{\re}{\operatorname{Re}}
\newcommand{\ov}{\overline}

\newtheorem{thm}{\hspace{\parindent}Theorem}[section]

\newtheorem{cor}[thm]{\hspace{\parindent}Corollary}
\newtheorem{lem}[thm]{\hspace{\parindent}Lemma}

\usepackage{amsmath}

\makeatletter
\def\upintkern@{\mkern-7mu\mathchoice{\mkern-3.5mu}{}{}{}}
\def\upintdots@{\mathchoice{\mkern-4mu\@cdots\mkern-4mu}%
 {{\cdotp}\mkern1.5mu{\cdotp}\mkern1.5mu{\cdotp}}%
 {{\cdotp}\mkern1mu{\cdotp}\mkern1mu{\cdotp}}%
 {{\cdotp}\mkern1mu{\cdotp}\mkern1mu{\cdotp}}}

\newcommand{\UpMultiIntegral}[1]{%
  \edef\ints@c{\noexpand\upintop
    \ifnum#1=\z@\noexpand\upintdots@\else\noexpand\upintkern@\fi
    \ifnum#1>\tw@\noexpand\upintop\noexpand\upintkern@\fi
    \ifnum#1>\thr@@\noexpand\upintop\noexpand\upintkern@\fi
    \noexpand\upintop
    \noexpand\ilimits@
  }%
  \futurelet\@let@token\ints@a
}
\makeatother

\DeclareFontFamily{OMX}{mdbch}{}
\DeclareFontShape{OMX}{mdbch}{m}{n}{ <->s * [0.8]  mdbchr7v }{}
\DeclareFontShape{OMX}{mdbch}{b}{n}{ <->s * [0.8]  mdbchb7v }{}
\DeclareFontShape{OMX}{mdbch}{bx}{n}{<->ssub * mdbch/b/n}{}

\DeclareSymbolFont{uplargesymbols}{OMX}{mdbch}{m}{n}
\SetSymbolFont{uplargesymbols}{bold}{OMX}{mdbch}{b}{n}
\DeclareMathSymbol{\upintop}{\mathop}{uplargesymbols}{82}
\DeclareMathSymbol{\upointop}{\mathop}{uplargesymbols}{"48}

\DeclareFontEncoding{MDB}{}{}
\DeclareFontFamily{MDB}{mdbch}{}
\DeclareFontShape{MDB}{mdbch}{m}{n}{ <->s * [0.8]  mdbchrmb }{}
\DeclareFontShape{MDB}{mdbch}{b}{n}{ <->s * [0.8]  mdbchbmb }{}
\DeclareFontShape{MDB}{mdbch}{bx}{n}{<->ssub * mdbch/b/n}{}
\DeclareFontSubstitution{MDB}{cmr}{m}{n}
\DeclareSymbolFont{mathdesignB}{MDB}{mdbch}{m}{n}%
\SetSymbolFont{mathdesignB}{bold}{MDB}{mdbch}{b}{n}%
\DeclareMathSymbol{\upintclockwise}{\mathop}{mathdesignB}{128}
\DeclareMathSymbol{\upointclockwise}{\mathop}{mathdesignB}{130}
\DeclareMathSymbol{\upointctrclockwise}{\mathop}{mathdesignB}{132}
\DeclareMathSymbol{\upoiint}{\mathop}{mathdesignB}{134}
\DeclareMathSymbol{\upoiiint}{\mathop}{mathdesignB}{136}

\makeatletter
\newcommand{\upint}{\DOTSI\upintop\ilimits@}
\newcommand{\upoint}{\DOTSI\upointop\ilimits@}
\makeatother

\theoremstyle{remark}

\newtheorem*{rem*}{Remark}

\newcommand\Li{{\rm Lip}}
\newcommand\fM{\frak M}

\newcommand\dg{\frak D}

\newcommand{\fI}{{\frak I}}

\newcommand\mC{\mathcal{C}}

\newcommand\dt{\frak t}

\newcommand\sG{\mathscr G}

\newcommand{\qm}{\quad\mbox{and}\quad}

\newcommand\mB{\mathcal{B}}
\newcommand\fB{\frak B}

\newcommand{\fF}{{\frak F}}

\newcommand{\Lip}{\operatorname{Lip}}

\begin{document}

\newcommand{\vse}{\vspace{.2in}}
\numberwithin{equation}{section}

\title{Multiple operator integrals in perturbation theory}
\author{V.V. Peller}
\thanks{The author is partially supported by NSF grant DMS 1300924}

\begin{abstract}
The purpose of this survey article is to give an introduction to double operator integrals and multiple operator integrals and to discuss various applications of such operator integrals in perturbation theory. We start with the Birman--Solomyak approach to define double operator integrals and consider applications in estimating operator differences $f(A)-f(B)$ for self-adjoint operators $A$ and $B$. Next, we present the Birman--Solomyak approach to the Lifshits--Krein trace formula that is based on double operator integrals. We study the class of operator Lipschitz functions, operator differentiable functions, operator H\"older functions, obtain Schatten--von Neumann estimates for operator differences. Finally, we consider in Chapter 1  estimates of functions of normal operators and functions of $d$-tuples of commuting self-adjoint operators under perturbations. 

In Chapter 2 we define multiple operator integrals in the case when the integrands belong to the integral projective tensor product of $L^\be$ spaces. We consider applications of such multiple operator integrals to the problem of the existence of higher operator derivatives and to the problem of estimating higher operator differences. We also consider connections with trace formulae for functions of operators under perturbations of class $\bS_m$, $m\ge2$.

In the last chapter we define Haagerup-like tensor products of the first kind and of the second kind and we use them to study functions of noncommuting self-adjoint operators under perturbation. We show that for functions $f$ in the Besov class $B_{\be,1}^1(\R^2)$
and for $p\in[1,2]$ we have a Lipschitz type estimate in the Schatten--von Neumann norm $\bS_p$ for functions of pairs of noncommuting self-adjoint operators, but there is no such a Lipschitz type estimate in the norm of $\bS_p$ with $p>2$ as well as in the operator norm. We also use triple operator integrals to estimate the trace norms of commutators of functions of almost commuting self-adjoint operators and extend the Helton--Howe trace formula for arbitrary functions in the Besov space $B_{\be,1}^1(\R^2)$.
\end{abstract}

\clearpage

\maketitle

\tableofcontents

%
%

\

\setcounter{section}{0}
\section{Introduction}
\label{In}

\medskip

In this survey article we study the role of double operator integrals and multiple operator integrals in perturbation theory. Double operator integrals appeared in the paper
\cite{DK} by Yu.L. Daletskii and S.G. Krein. In that paper they considered the problem of differentiating the operator-valued function $t\mapsto f(A+tK)$, where $A$ and $K$ are self-adjoint operators on Hilbert space. They discovered the following formula that expresses the derivative in terms of double operator integrals:
$$
\frac{d}{dt}\big(f(A+tK)\big)\Big|_{t=0}=
\iint\limits_{\R\times\R}\frac{f(x)-f(y)}{x-y}\,dE_A(x)K\,dE_A(y)
$$
for sufficiently nice functions $f$. Here $E_A$ stands for the spectral measure of $A$.

That time there was no rigorous theory of double operator integrals. Such a theory was developed later by Birman and Solomyak in \cite{BS1}, \cite{BS2} and \cite{BS3}.

In general double operator integrals are expressions of the form
$$
\iint\Phi(x,y)\,dE_1(x)T\,dE_2(y),
$$
where $\Phi$ is a measurable function, $T$ is a linear operator, and $E_1$ and $E_2$ are spectral measures on Hilbert space.

The Birman--Solomyak approach allows one to define such integrals in the case when $T$ is a Hilbert Schmidt operator and $\Phi$ is an arbitrary bounded measurable function. This, in turn, permits us to define double operator integrals for arbitrary bounded linear operators $T$ and for functions $\Phi$ satisfying certain assumptions (such functions are called {\it Schur multipliers}).

It turned out that double operator integrals play a very important role in perturbation theory. They appear naturally when estimating various norms of operator differences $f(A)-f(B)$, where $A$ is an unperturbed operator and $B$ is a perturbed operator. In particular, double operator integrals are very helpful when studying the class of {\it operator Lipschitz functions}, i.e., functions $f$ on $\R$, for which
\bay
\label{oLnervo}
\|f(A)-f(B)\|\le\const\|A-B\|.
\ey
It turns out that if inequality \rf{oLnervo} holds for all bounded self-adjoint operators $A$ and $B$, then the same inequality holds for unbounded $A$ and $B$ once $A-B$ is bounded. 

Roughly speaking, a functions $f$ on $\R$ is operator Lipschitz if and only if the divided difference $(x,y)\mapsto\frac{f(x)-f(y)}{x-y}$ is a Schur multiplier.

It is obvious that operator Lipschitz functions $f$ must be {\it Lipschitz}, i.e.,
the inequality $|f(x)-f(y)|\le\const|x-y|$ must hold for $x,\,y\in\R$. The question whether the converse is true was resolved in negative by Farforovskaya in \cite{F}.
Later McIntosh \cite{Mc} and Kato \cite{Ka} proved that the function $x\mapsto|x|$ is not
operator Lischitz. Then in \cite{JW} it was shown that operator Lipschitz functions must be differentiable everywhere on $\R$ (but {\it not necessarily continuously differentiable}, see \cite{KSh0}. Later in \cite{Pe2} necessary conditions for operator Lipschitzness were found in terms of Besov spaces and Carleson measures
(see also \cite{Pe}).

In Chapter 1 we give an introduction to the theory of double operator integrals and define and characterize the class of Schur multipliers. Then we consider various applications of double operator integrals in perturbation theory. Namely, we study operator Lipschitz functions, operator H\"older functions, operator differentiable functions. We obtain sharp estimates for Schatten--von Neumann norms of operator differences $f(A)-f(B)$ for functions $f$ in the H\"older class $\L_\a(\R)$. We present the Birman--Solomyak approach to the Lifshits--Krein trace formula that is based on double operator integrals. We also consider similar problems for functions of normal operators and for functions of $m$-tuples of commuting self-adjoint operators.

In Chapter 2 we proceed to multiple operator integrals, i.e., expressions of the form
$$
\underbrace{\int\cdots\int}_m\Phi(x_1,x_2,\cdots,x_m)
\,dE_1(x_1)T_1\,dE_2(x_2)T_2\cdots T_{m-1}\,dE_m(x_m).
$$
We follow the approach to multiple operator integrals given in \cite{Pe4} and define such multiple operator integrals in the case when the {\it integrand} $\Psi$ belongs to the (integral) projective tensor products of the spaces $L^\be(E_j)$, $1\le j\le m$.
We use this approach to study the problem of the existence of higher operator derivatives of the function $t\mapsto f(A+tK)$ and express higher operator derivatives in terms of multiple operator integrals. We also use multiple operator integrals to obtain sharp estimates of higher operator differences
$$
\sum_{j=0}^m(-1)^{m-j}\left(\begin{matrix}m\\j\end{matrix}\right)f\big(A+jK\big).
$$
Finally, in the last section of Chapter 2 we apply multiple operator integrals to trace formulae for functions of self-adjoint operators of class $\bS_m$ with $m\in\Z$, $m\ge2$.

An alternative approach to multiple operator integrals is given in \cite{JTT}. That approach is based on the Haagerup tensor product of $L^\be$ spaces. We define in Chapter 3 triple operator integrals whose integrands belong to the Haaherup tensor product of three $L^\be$ spaces. We study Schatten--von Neumann properties of such triple operator integrals and we see that their Schatten--von Neumann properties are not as nice as in the case of triple operator integrals with integrands in the integral projective tensor product. 

We are going to use triple operator integrals to estimate functions of pairs of noncommuting self-adjoint operators under perturbation. It turns out that for our purposes none of the approaches based on the integral projective tensor product and on the Haagerup tensor product of $L^\be$ spaces works. We define new tensor products and call them Haagerup-like tensor products of the first kind and of the second kind. Then we define triple operator integrals with integrands in such Haagerup-like tensor products 
and use them to estimate the norms
$\|f(A_1,B_1)-f(A_2,B_2)\|$, where $(A_2,B_2)$ is a perturbation of $(A_1,B_1)$ and $f$ is a function in the Besov space $B_{\be,1}^1(\R^2)$. 

Note that functions $f(A,B)$ for not necessarily commuting self-adjoint operators are defined as double operator integrals
$$
f(A,B)=\iint f(x,y)\,dE_A(x)\,dE_B(y).
$$

We show that for $p\in[1,2]$, we have a Lipschitz type estimate in the Schatten--von Neumann norm $\bS_p$, but such Lipschitz type estimates do not hold in $\bS_p$ with $p>2$ as well as in the operator norm. We conclude the chapter with estimating commutators of almost commuting self-adjoint operators ($A$ and $B$ are called {\it almost commuting} is $AB-BA\in\bS_1$). Such estimates allow us to extend the Helton--Howe trace formula for arbitrary functions in the Besov class 
$B_{\be,1}^1(\R^2)$. The results of the last chapter were obtained recently in \cite{ANP1}, \cite{ANP2}, \cite{ANP3} and \cite{AP10}.

I am grateful to A.B. Aleksandrov for helpful remarks.

\

\section{Preliminaries}
\label{prel}

\

In this section we collect necessary information on function spaces and operator ideals.

\medskip

{\bf 2.1. Besov classes of functions on Euclidean spaces and Littlewood--Paley type expansions.} The technique of Littlewood--Paley type expansions of functions or distributions on Euclidean spaces 
is a very important tool in Harmonic Analysis. 

Let $w$ be an infinitely differentiable function on $\R$ such
that
\bay
\label{w}
w\ge0,\quad\supp w\subset\left[\frac12,2\right],\quad\mbox{and} \quad w(s)=1-w\left(\frac s2\right)\quad\mbox{for}\quad s\in[1,2].
\ey

We define the functions $W_n$, $n\in\Z$, on $\R^d$ by 
$$
\big(\F W_n\big)(x)=w\left(\frac{\|x\|}{2^n}\right),\quad l\in\Z, \quad x=(x_1,\cdots,x_d),
\quad\|x\|\df\left(\sum_{j=1}^dx_j^2\right)^{1/2},
$$
where $\F$ is the {\it Fourier transform} defined on $L^1\big(\R^d\big)$ by
$$
\big(\F f\big)(t)=\!\int\limits_{\R^d} f(x)e^{-{\rm i}(x,t)}\,dx,\!\quad 
x=(x_1,\cdots,x_d),
\quad t=(t_1,\cdots,t_d), \!\quad(x,t)\df \sum_{j=1}^dx_jt_j.
$$
Clearly,
$$
\sum_{n\in\Z}(\F W_n)(t)=1,\quad t\in\R^d\setminus\{0\}.
$$

With each tempered distribution $f\in{\mathscr S}^\prime\big(\R^d\big)$, we
associate the sequence $\{f_n\}_{n\in\Z}$,
\bay
\label{fn}
f_n\df f*W_n.
\ey
The formal series
$
\sum_{n\in\Z}f_n
$
is a Littlewood--Paley type expansion of $f$. This series does not necessarily converge to $f$. Note that in this paper a significant role is played by the Besov spaces $B_{\be,1}^1(\R^d)$ (see the definition below).
For functions $f\in B_{\be,1}^1(\R^d)$, we have
$$
f(x)-f(y)=\sum_{n\in\Z}\big(f_n(x)-f_n(y)\big),\quad x,~y\in\R^d,
$$
and the series on the right converges uniformly.

Initially we define the (homogeneous) Besov class $\dot B^s_{p,q}\big(\R^d\big)$,
$s>0$, $1\le p,\,q\le\be$, as the space of all
$f\in{\mathscr S}^\prime(\R^n)$
such that
\bay
\label{Wn}
\{2^{ns}\|f_n\|_{L^p}\}_{n\in\Z}\in\ell^q(\Z)
\ey
and put
$$
\|f\|_{B^s_{p,q}}\df\big\|\{2^{ns}\|f_n\|_{L^p}\}_{n\in\Z}\big\|_{\ell^q(\Z)}.
$$
According to this definition, the space $\dot B^s_{p,q}(\R^n)$ contains all polynomials
and all polynomials $f$ satisfy the equality $\|f\|_{B^s_{p,q}}=0$. Moreover, the distribution $f$ is determined by the sequence $\{f_n\}_{n\in\Z}$
uniquely up to a polynomial. It is easy to see that the series 
$\sum_{n\ge0}f_n$ converges in ${\mathscr S}^\prime(\R^d)$.
However, the series $\sum_{n<0}f_n$ can diverge in general. It can easily be proved that the series
\bay
\label{ryad}
\sum_{n<0}\frac{\partial^r f_n}{\partial x_1^{r_1}\cdots\partial x_d^{r_d}},\qquad \mbox{where}\quad r_j\ge0,\quad\mbox{for}\quad
1\le j\le d,\quad\sum_{j=1}^dr_j=r,
\ey
converges uniformly on $\R^d$ for every nonnegative integer
$r>s-d/p$. Note that in the case $q=1$ the series \rf{ryad}
converges uniformly, whenever $r\ge s-d/p$.

Now we can define the modified (homogeneous) Besov class $B^s_{p,q}\big(\R^d\big)$. We say that a distribution $f$
belongs to $B^s_{p,q}(\R^d)$ if \rf{Wn} holds and
$$
\frac{\partial^r f}{\partial x_1^{r_1}\cdots\partial x_d^{r_d}}
=\sum_{n\in\Z}\frac{\partial^r f_n}{\partial x_1^{r_1}\cdots\partial x_d^{r_d}},\quad
\mbox{whenever}\quad 
r_j\ge0,\quad\mbox{for}\quad
1\le j\le d,\quad\sum_{j=1}^dr_j=r.
$$
in the space ${\mathscr S}^\prime\big(\R^d\big)$, where $r$ is
the minimal nonnegative integer such that $r>s-d/p$ ($r\ge s-d/p$ if $q=1$). Now the function $f$ is determined uniquely by the sequence $\{f_n\}_{n\in\Z}$ up
to a polynomial of degree less than $r$, and a polynomial $g$ belongs to 
$B^s_{p,q}\big(\R^d\big)$
if and only if $\deg g<r$. 

In the case when $p=q$ we use the notation $B_p^s(\R^d)$ for $B_{p,p}^s(\R^d)$.

Consider now the scale $\L_\a(\R^d)$, $\a>0$, of H\"older--Zygmund classes. They can be defined by $\L_\a(\R^d)\df B_\be^\a(\R^d)$. We need a description of $\L_\a$ in terms of convolutions with de la Vall\'ee Poussin type kernel $V_n$.

To define a de la Vall\'ee Poussin type kernel $V_n$, we define the $C^\be$ function $v$ on $\R$ by
\bay
\label{dlVP}
v(x)=1\quad\mbox{for}\quad x\in[-1,1]\quad\mbox{and}\quad v(x)=w(\|x\|)\quad\mbox{if}\quad \|x\|\ge1,
\ey
where $w$ is the function defined by \rf{w}. We define $V_n$, $n\in\Z$, by
$$
\F V_n(x)\df v\left(\frac{\|x\|}{2^n}\right),\quad n\in\Z,\quad x\in\R^d.
$$
In the definition of the classes $\L_\a(\R^d)$, $\a>0$, we can replace the condition  
$\|f_n\|_{L^\be}\le\const2^{-n\a}$, $n\in\Z$, with the condition
\bay
\label{dlVPa}
\|f-f*V_n\|_{L^\be}\le\const2^{-n\a},\quad n\in\Z.
\ey

In the case of Besov classes $B_{\be,q}^s(\R^d)$ the functions $f_n$, defined by \rf{fn} have the following properties: $f_n\in L^\be(\R^d)$ and 
$\supp\F f\subset\{\xi\in\R^d:~\|\xi\|\le2^{n+1}\}$. Such functions can be characterized by the following Paley--Wiener--Schwartz type theorem  (see \cite{R}, Theorem 7.23 and exercise 15 of Chapter 7):

{\it Let $f$ be a continuous function
on $\R^d$ and let $M,\,\s>0$. The following statements are equivalent:

{\em(i)} $|f|\le M$ and $\supp\F f\subset\{\xi\in\R^d:\|\xi\|\le\s\}$;

{\em(ii)} $f$ is a restriction to $\R^d$ of an entire function on $\C^d$ such that 
$$
|f(z)|\le Me^{\s\|\im z\|}
$$
for all $z\in\C^d$.}

Besov classes admit many other descriptions.
We give here the definition in terms of finite differences.
For $h\in\R^d$, we define the difference operator $\D_h$,
$$
(\D_hf)(x)=f(x+h)-f(x),\quad x\in\R^d.
$$
It is easy to see that $B_{p,q}^s\big(\R^d\big)\subset L^1_{\rm loc}\big(\R^d\big)$ for every $s>0$
and $B_{p,q}^s\big(\R^d\big)\subset C\big(\R^d\big)$ for every $s>d/p$. Let $s>0$ and let $m$ be the integer such that $m-1\le s<m$.
The Besov space $B_{p,q}^s\big(\R^d\big)$ can be defined as the set of
functions $f\in L^1_{\rm loc}\big(\R^d\big)$ such that
$$
\int_{\R^d}|h|^{-d-sq}\|\D^m_h f\|_{L^p}^q\,dh<\be\quad\mbox{for}\quad q<\be
$$
and
$$
\sup_{h\not=0}\frac{\|\D^m_h f\|_{L^p}}{|h|^s}<\be\quad\mbox{for}\quad q=\be.
$$
However, with this definition the Besov space can contain polynomials of higher degree than in the case of the first definition given above.

We refer the reader to \cite{Pee} and \cite{Tr} for more detailed information on Besov spaces.

\medskip

{\bf 2.2. Besov classes of periodic functions.} Studying periodic functions on $\R^d$ is equivalent to studying functions on the $d$-dimensional torus $\T^d$. To define Besov spaces on $\T^d$, we consider a function $w$ satisfying \rf{w} and define the trigonometric polynomials $W_n$, $n\ge0$, by
$$
W_n(\z)\df\sum_{j\in\Z^d}w\left(\frac{|j|}{2^n}\right)\z^j,\quad n\ge1,
\quad W_0(\z)\df\sum_{\{j:|j|\le1\}}\z^j,
$$
where 
$$
\z=(\z_1,\cdots,\z_d)\in\T^d,\quad j=(j_1,\cdots,j_d),\quad\mbox{and}\quad
|j|=\big(|j_1|^2+\cdots+|j_d|^2\big)^{1/2}.
$$
For a distribution $f$ on $\T^d$ we put
$$
f_n=f*W_n,\quad n\ge0,
$$
and we say that $f$ belongs the Besov class $B_{p,q}^s(\T^d)$, $s>0$, 
$1\le p,\,q\le\be$, if
\bay
\label{Bperf}
\big\{2^ns\|f_n\|_{L^p}\big\}_{n\ge0}\in\ell^q.
\ey

Note that locally the Besov space $B_{p,q}^s(\R^d)$ coincides with the Besov space
$B_{p,q}^s$ of periodic functions on $\R^d$.

\medskip

{\bf 2.3. Operator ideals.}
For a bounded linear operator $T$ on Hilbert space, we consider its singular
values $s_j(T)$, $j\ge0$, 
$$
s_j(T)\df\inf\big\{\|T-R\|:~\rank R\le j\big\}.
$$

Let $\bS_p$, $0<p<\be$, be the Schatten--von Neumann class of operators $T$ on Hilbert space
such that
$$
\|T\|_{\bS_p}\df\left(\sum_{j\ge0}\big(s_j(T)\big)^p\right)^{1/p}<\be.
$$
This is a normed ideal for $p\ge1$. 
The class $\bS_1$ is called {\it trace class}. For a linear operators $T$ on a Hilbert space $\h$ its trace is defined by
$$
\trace T\df\sum_{j\ge0}(Te_j,e_j),
$$
where $\{e_j\}_{j\ge0}$ is an orthonormal basis in $\h$. The right-hand side does not depend on the choice of a basis.

The class $\bS_2$ is called the {\it Hilbert--Schmidt class}. It is a Hilbert space with inner product
$$
(T,R)_{\bS_2}\df\trace(TR^*).
$$

For $p\in(1,\be)$, the dual space $(\bS_p)^*$ can be isometrically identified with
$\bS_{p'}$ with respect to the pairing
$$
\langle T,R\rangle\df\trace(TR).
$$

The dual space to $\bS_1$ can be identified with the space of bounded linear operators, while the dual space to the space of compact operators can be identified with $\bS_1$ with respect to the same pairing.

We refer the reader to \cite{GK} for detailed information on singular values and operator ideals.

\

\begin{center}
\bf\huge Chapter 1
\end{center}

\

\begin{center}
\bf\Large Applications of double operator integrals \\
in perturbation theory
\end{center}

\


\addtocontents{toc}{\vspace*{.3cm}\textbf{{\sc Chapter 1}. Applications of double operator integrals in perturbation \\ \hspace*{2.75cm}theory}\hfill\pageref{Adoipt}}
\renewcommand{\thesection}{1.\arabic{section}} 
\setcounter{section}{0}
\label{Adoipt}

\

In the first chapter we give an introduction to the theory of double operator integrals that was developed by Birman and Solomyak. We discuss the problem of a representation for operator differences $f(A)-f(B)$ in terms of double operator integrals. This allows us to obtain necessary conditions and sufficient conditions for a function on the real line to be operator Lipschitz. In particular, we show that if $f$ belongs to the Besov class $B_{\be,1}^1(\R)$, then $f$ is operator Lipschitz. It turns out that the same condition
$f\in B_{\be,1}^1(\R)$ is also sufficient for operator differentiability.
Next, we present the Birman--Solomyak approach to the Lifshits--Krein trace formula. Their approach is based on double operator integrals.
We also discuss H\"older type estimates and Schatten--von Neumann estimates for operator differences. 

Finally, we consider perturbations of functions of normal operators and perturbations of functions of $m$-tuples of commuting self-adjoint operators.

\

\section{An introduction to double operator integrals}
\label{dois}

\

Double operator integrals appeared in the paper \cite{DK} by Daletskii and S.G. Krein. It was Birman and Solomyak who developed later the beautiful theory of double operator integrals in \cite{BS1}, \cite{BS2}, and \cite{BS3}. 

Let $(\X,E_1)$ and $(\Y,E_2)$ be spaces with spectral measures $E_1$ and $E_2$
on a Hilbert space $\h$. The idea of Birman and Solomyak is to define first
double operator integrals
\bay
\label{doi}
\int\limits_\X\int\limits_\Y\Phi(x,y)\,d E_1(x)T\,dE_2(y)
\ey
for bounded measurable functions $\Phi$ and operators $T$
of Hilbert Schmidt class $\bS_2$. Consider the spectral measure $\E$ whose values are orthogonal
projections on the Hilbert space $\bS_2$, which is defined by
$$
\E(\L\times\D)T=E_1(\L)TE_2(\D),\quad T\in\bS_2,
$$
$\L$ and $\D$ being measurable subsets of $\X$ and $\Y$. 
Obviously, left multiplication by $E_1(\L)$ commutes with right multiplication by 
$E_2(\D)$.
It was shown in \cite{BS} that $\E$ extends to a spectral measure on
$\X\times\Y$ and if $\Phi$ is a bounded measurable function on $\X\times\Y$, by definition,
$$
\int\limits_\X\int\limits_\Y\Phi(x,y)\,d E_1(x)T\,dE_2(y)\df
\left(\,\,\int\limits_{\X\times\Y}\Phi\,d\E\right)T.
$$
Clearly,
$$
\left\|\int\limits_\X\int\limits_\Y\Phi(x,y)\,dE_1(x)T\,dE_2(y)\right\|_{\bS_2}
\le\|\Phi\|_{L^\be}\|T\|_{\bS_2}.
$$
If
$$
\int\limits_\X\int\limits_\Y\Phi(x,y)\,d E_1(x)T\,dE_2(y)\in\bS_1
$$
for every $T\in\bS_1$, we say that $\Phi$ is a {\it Schur multiplier of $\bS_1$ associated with
the spectral measures $E_1$ and $E_2$}.

To define double operator integrals of the form \rf{doi} for bounded linear operators $T$,
we consider the transformer
$$
Q\mapsto\int\limits_{\Y}\int\limits_{\X}\Phi(y,x)\,d E_2(y)\,Q\,dE_1(x),\quad Q\in\bS_1,
$$
and assume that the function $(y,x)\mapsto\Phi(y,x)$ is a Schur multiplier of $\bS_1$ associated with $E_2$ and $E_1$.

In this case the transformer
\bay
\label{tra}
T\mapsto\int\limits_\X\int\limits_\Y\Phi(x,y)\,d E_1(x)T\,dE_2(y),\quad T\in \bS_2,
\ey
extends by duality to a bounded linear transformer on the space of bounded linear operators on $\h$
and we say that the function $\Phi$ is {\it a Schur multiplier (with respect to $E_1$ and $E_2$) of the space of bounded linear operators}.
We denote the space of such Schur multipliers by $\fM(E_1,E_2)$.
The norm of $\Phi$ in $\fM(E_1,E_2)$ is, by definition, the norm of the
transformer \rf{tra} on the space of bounded linear operators.

The function $\Phi$ in \rf{tra} is called the {\it integrand} of the double operator integral.

Note that the term Schur multiplier in the context of double operator integrals was introduced in \cite{Pe2}. This is a generalization of the notion of a matrix Schur multiplier. Indeed, consider the very special case when the Hilbert space is the sequence space $\ell^2$ and both spectral measures $E_1$ and $E_2$ are defined on the $\s$-algebra of all subsets of $\Z_+$ as follows: $E_1(\D)=E_2(\D)$ is the orthogonal projection onto the closed linear span of the vectors $e_n$, $n\in\D$, where $\{e_n\}_{n\ge0}$ is the standard orthonormal basis of $\ell^2$. In this case a function $\Phi$ on $\Z_+\times\Z_+$ is a Schur multiplier if and only if the matrix $\{\Phi(m,n)\}_{m,n\ge0}$ (for which we keep the notation 
$\Phi$) is a {\it matrix Schur multiplier}, i.e., 
$$
T=\{t_{jk}\}_{j,k\ge0}\in\mB\quad\Longrightarrow\quad \Phi\star T\in\mB,
$$
where $\mB$ is the space of matrices that induce bounded linear operators on $\ell^2$ and $\Phi\star T$ is the {\it Hadamard--Schur product} of the matrices 
$\Phi$ and $T$. Recall that the Hadamard--Schur product $A\star B$ of matrices
 $A=\{a_{jk}\}_{j,k\ge0}$ and $B=\{b_{jk}\}_{j,k\ge0}$ is defined by
$$
(A\star B)_{jk}=a_{jk}b_{jk},\quad j,~k\in\Z_+.
$$

It is easy to see that if a function $\Phi$ on $\X\times\Y$ belongs to the {\it projective tensor
product}
$L^\be(E_1)\hat\otimes L^\be(E_2)$ of $L^\be(E_1)$ and $L^\be(E_2)$ (i.e., $\Phi$ admits a representation
$$
\Phi(x,y)=\sum_{n\ge0}\f_n(x)\psi_n(y),
$$
where $\f_n\in L^\be(E_1)$, $\psi_n\in L^\be(E_2)$, and
$$
\sum_{n\ge0}\|\f_n\|_{L^\be}\|\psi_n\|_{L^\be}<\be),
$$
then $\Phi\in\fM(E_1,E_2)$.
For such function $\Phi$ we have
$$
\int\limits_\X\int\limits_\Y\Phi(x,y)\,dE_1(x)T\,dE_2(y)=
\sum_{n\ge0}\left(\,\int\limits_\X\f_n\,dE_1\right)T\left(\,\int\limits_\Y\psi_n\,dE_2\right).
$$

More generally, $\Phi\in\fM(E_1,E_2)$ if $\Phi$
belongs to the {\it integral projective tensor product} $L^\be(E_1)\hat\otimes_{\rm i}
L^\be(E_2)$ of $L^\be(E_1)$ and $L^\be(E_2)$, i.e., $\Phi$ admits a representation
\bay
\label{ipt}
\Phi(x,y)=\int_\O \f(x,w)\psi(y,w)\,d\l(w),
\ey
where $(\O,\l)$ is a $\s$-finite measure space, $\f$ is a measurable function on $\X\times \O$,
$\psi$ is a measurable function on $\Y\times \O$, and
\bay
\label{ir}
\int_\O\|\f(\cdot,w)\|_{L^\be(E_1)}\|\psi(\cdot,w)\|_{L^\be(E_2)}\,d\l(w)<\be.
\ey
If $\Phi\in L^\be(E_1)\hat\otimes_{\rm i}L^\be(E_2)$, then
$$
\int\limits_\X\int\limits_\Y\Phi(x,y)\,dE_1(x)T\,dE_2(y)=
\int\limits_\O\left(\,\int\limits_\X\f(x,w)\,dE_1(x)\right)T
\left(\,\int\limits_\Y\psi(y,w)\,dE_2(y)\right)\,d\l(w).
$$
Clearly, the function
$$
w\mapsto \left(\,\int_\X\f(x,w)\,dE_1(x)\right)T
\left(\,\int_\Y\psi(y,w)\,dE_2(y)\right)
$$
is weakly measurable and
$$
\int\limits_\O\left\|\left(\,\int\limits_\X\f(x,w)\,dE_1(x)\right)T
\left(\,\int\limits_\Y\psi(y,w)\,dE_2(w)\right)\right\|\,d\l(w)<\be.
$$
It is easy to see that 
$$
\|\Phi\|_{\fM(E_1,E_2)}\le\|\Phi\|_{L^\be(E_1)\hat\otimes_{\rm i}L^\be(E_2)},
$$
where $\|\Phi\|_{L^\be\hat\otimes_{\rm i}L^\be}$ is, by definition, the infimum of
the left-hand side of \rf{ir} over all representations of $\Phi$ of the form 

It turns out that all Schur multipliers can be obtained in this way (see Theorem \ref{tomSc} below).

Another sufficient condition for a function to be a Schur multiplier can be stated in terms of the Haagerup tensor products of $L^\be$ spaces. The {\it Haagerup tensor product} $L^\be(E_1)\!\otimes_{\rm h}\!L^\be(E_2)$ can be defined as the space of functions $\Phi$ of the form 
\bay
\label{FiH}
\Phi(x,y)=\sum_{n\ge0}\f_n(x)\psi_n(y),
\ey
where $\f_n\in L^\be(E_1)$, $\psi_n\in L^\be(E_2)$ and
$$
\{\f_n\}_{n\ge0}\in L_{E_1}^\be(\ell^2)\quad\mbox{and}\quad
\{\psi_n\}_{n\ge0}\in L_{E_2}^\be(\ell^2).
$$
The {\it norm of $\Phi$ in $L^\be(E_1)\!\otimes_{\rm h}\!L^\be(E_2)$ is defined as} the infimum of
$$
\big\|\{\f_n\}_{n\ge0}\big\|_{L_{E_1}^\be(\ell^2)}
\big\|\{\psi_n\}_{n\ge0}\big\|_{L_{E_2}^\be(\ell^2)}
$$
over all representations of $\Phi$ of the form \rf{FiH}. Here
$$
\big\|\{\f_n\}_{n\ge0}\big\|_{L_{E_1}^\be(\ell^2)}\df
\Big\|\sum_{n\ge0}|\f_n|^2\Big\|_{L^\be(E_1)}^{1/2}\quad\!\!\mbox{and}\quad\!\!
\big\|\{\psi_n\}_{n\ge0}\big\|_{L_{E_1}^\be(\ell^2)}\df
\Big\|\sum_{n\ge0}|\psi_n|^2\Big\|_{L^\be(E_2)}^{1/2}.
$$
It can easily be verified that if $\Phi\in L^\be(E_1)\!\otimes_{\rm h}\!L^\be(E_2)$, then $\Phi\in\fM(E_1,E_2)$ and
$$
\iint\Phi(x,y)\,dE_1(x)T\,dE_2(y)=
\sum_{n\ge0}\Big(\int\f_n\,dE_1\Big)T\Big(\int\psi_n\,dE_2\Big).
$$
It is also easy to see that the series on the right converges in the weak operator topology and
$$
\|\Phi\|_{\fM(E_1,E_2)}\le\|\Phi\|_{L^\be(E_1)\otimes_{\rm h}L^\be(E_2)}.
$$

As the following theorem says, the condition 
$\Phi\in L^\be(E_1)\!\otimes_{\rm h}\!L^\be(E_2)$ is not only sufficient, but also necessary. 

\begin{thm}
\label{tomSc} 
Let $\Phi$ be a measurable function on
$\X\times\Y$ and let $\mu$ and $\nu$ be positive measures on $\X$ and $\Y$ that are mutually absolutely continuous with respect to $E_1$ and $E_2$. The following are equivalent:

{\rm (i)} $\Phi\in\fM(E_1,E_2)$;

{\rm (ii)} $\Phi\in L^\be(E_1)\hat\otimes_{\rm i}L^\be(E_2)$;

{\rm (iii)} $\Phi\in L^\be(E_1)\!\otimes_{\rm h}\!L^\be(E_2)$;

{\rm (iv)} there exist measurable functions $\f$ on $\X\times\O$ and $\psi$ on $\Y\times\O$ such that
{\em\rf{ipt}} holds and
\bay
\label{bs}
\left\|\left(\int_\O|\f(\cdot,w)|^2\,d\l(w)\right)^{1/2}\right\|_{L^\be(E_1)}
\left\|\left(\int_\O|\psi(\cdot,w)|^2\,d\l(w)\right)^{1/2}\right\|_{L^\be(E_2)}<\be;
\ey

{\rm (v)} if the integral operator
$f\mapsto\int k(x,y)f(y)\,d\nu(y)$ from $L^2(\nu)$ to $L^2(\mu)$ belongs to $\bS_1$, then the same is true for the integral operator 
$f\mapsto\int\Psi(x,y)k(x,y)f(y)\,d\nu(y)$.
\end{thm}

The implications (iv)$\imp$(i)$\eq$(v) were established in \cite{BS3}.
In the case of matrix Schur multipliers the fact that (i) implies (ii) was proved in \cite{Be}. We refer the reader to \cite{Pe2} for the proof of the equivalence of (i), (ii), and (iv) and to \cite{Pi} for the proof of the fact that (i) is equivalent to (iii).

Suppose that $\frak F_1$ and $\frak F_2$ are closed subsets of $\R$. We denote by
$\fM_{\fF_1,\fF_2}$ the space of functions that belong to $\fM(E_1,E_2)$ for arbitrary spectral measures $E_1$ and $E_2$ such that $\supp E_1\subset\fF_2$ and 
$\supp E_2\subset\fF_2$. 

It is well known (see \cite{KSh1} and \cite{KSh2}) that if $\Phi$ is a continuous function 
on $\frak F_1\times\frak F_2$ and
$E_1$ and $E_2$ are Borel spectral measures such that $\supp E_1=\fF_2$ and 
$\supp E_2=\fF_2$, then $\Phi\in\fM_{\fF_1,\fF_2}$ if and only if $\Phi\in\fM(E_1,E_2)$.
The same conclusion under the weaker assumption that $\Phi$ is continuous in each variable
was established in \cite{AP5}.

It is easy to see that conditions (i) - (iv) are also equivalent to the fact that $\Phi$ is a Schur multiplier of $\bS_1$. It follows that if $\fI$ is an operator ideal that is an interpolation ideal between the space of bounded linear operators and trace class $\bS_1$ and $\Phi$ satisfies one of the conditions (i) - (iv), then
$\Phi$ is a Schur multiplier of $\fI$, i.e.,
$$
T\in\fI\quad\Longrightarrow\iint\Phi(x,y)\,dE_1(x)T\,dE_2(y)\in\fI.
$$
In particular, this is true when $\fI$ is the Schatten--von Neumann class $\bS_p$,
$1<p<\be$.

If $\fI$ is a separable (symmetrically normed) operator ideal (see \cite{GK}) , we say that a function $\Phi$ is a {\it Schur multiplier of} $\fI$ if the transformer $T\mapsto\iint\Phi\,dE_1T\,dE_2$ defined on $\bS_1$ admits an extension to a bounded linear operator on $\fI$. In the case when $\fI$ is an operator ideal dual to separable,
we can define Schur multiplier of $\fI$ by duality.
We denote the space of Schur multipliers of $\fI$ with respect to $E_1$ and $E_2$ by 
$\fM_\fI(E_1,E_2)$.

Consider now the case when $E_1=E_2=E$ and $T\in\bS_1$. It follows easily from Theorem
\ref{tomSc} that functions in the space $\fM(E,E)$ of Schur multipliers have traces on the diagonal
$$
\Phi\mapsto\Phi\Big|\{(x,x):~x\in\X\}
$$
and the traces of functions in $\M(E,E)$ belong to $L^\be(E)$.

The following useful fact was established in \cite{BS3}.

\begin{thm}
\label{sled}
Let $E$ be a spectral measure and $\Phi\in\fM(E,E)$. Suppose that $T\in\bS_1$. Then
\bay
\label{sleddoi}
\trace\left(\iint\Phi(x,y)\,dE(x)T\,dE(y)\right)=\int\Phi(x,x)\,d\mu(x),
\ey
where $\mu$ is the signed measure defined by
$$
\mu(\D)=\trace\big(TE(\D)\big).
$$
\end{thm}

\Pf It follows easily from Theorem \ref{tomSc} that it suffices to establish formula \rf{sleddoi} in the case $\Phi(x,y)=\f(x)\psi(y)$, $\f,\,\psi\in L^\be(E)$. We have
\begin{align*}
\trace\left(\iint\Phi(x,y)\,dE(x)T\,dE(y)\right)&=
\trace\left(\left(\int\f\,dE\right)T\left(\int\psi\,dE\right)\right)\\[.2cm]
&=\trace\left(\left(\int\psi\,dE\right)\left(\int\f\,dE\right)T\right)\\[.2cm]
&=\trace\left(\left(\int\f\psi\,dE\right)T\right)=\int\Phi(x,x)\,d\mu(x).
\quad\bl
\end{align*}

\

\section{A representation of operator differences in terms of double operator integrals} 
\label{razni}

\

In the paper \cite{DK} by Daletskii and S.G. Krein  under certain assumptions on a function $f$ on $\R$ the following formula was discovered:
\bay
\label{fDK}
f(A)-f(B)=
\int\limits_\R\int\limits_\R\frac{f(s)-f(t)}{s-t}\,dE_A(s)(A-B)\,dE_B(t)
\ey
for bounded self-adjoint operators $A$ and $B$. Here $E_A$ and $E_B$ are the spectral measures of $A$ and $B$. Later in Birman and Solomyak in \cite{BS3} proved formula \rf{fDK} in a much more general situation.

Consider first the case when $A-B$ belongs to the Hilbert--Schmidt class $\bS_2$. Suppose that $f$ is a {\it Lipschitz function} on $\R$, i.e.,
$$
|f(s)-f(t)|\le\const|s-t|,\quad s,~t\in\R,\quad
\|f\|_{\Lip}\df\sup_{s\ne t}\frac{|f(s)-f(t)|}{|s-t|}.
$$
Consider the divided difference $\dg f$ defined by
$$
(\dg f)(s,t)=\frac{f(s)-f(t)}{s-t},\quad s\ne t.
$$
It was established in \cite{BS3} that in the case $A-B\in\bS_2$, formula \rf{fDK} holds for arbitrary Lipschitz functions $f$.
To understand the right-hand side of \rf{fDK} for Lipschitz functions, we have to define the divided difference $\dg f$ on the diagonal $\D\df\{(x,x):\lb~x\in\R\}$. It turns out that
no matter how we can define $\dg f$ on the diagonal, formula \rf{fDK} holds. If $f$ is differentiable, it is natural to assume that $(\dg f)(s,s)=f'(s)$. We can also define $\dg f$ to be zero on the diagonal. 
Put
$$
(\dg_0f)(s,t)=\left\{\begin{array}{ll}
(\dg f)(s,t),& s\ne t,\\[.2cm]0,&s=t.
\end{array}\right.
$$

\begin{thm}
\label{LeS2}
Suppose that $A$ and $B$ are (not necessarily bounded) self-adjoint operators such that $A-B\in\bS_2$. Let $\bs{\frak f}$ be a bounded Borel function on $\R^2$ such that 
$\bs{\frak f}\big|_{\R^2\setminus\D}=\dg f\big|_{\R^2\setminus\D}$. 
Then
\bay
\label{frk}
f(A)-f(B)=
\int\limits_\R\int\limits_\R\bs{\frak f}(s,t)\,dE_A(s)(A-B)\,dE_B(t).
\ey
\end{thm}

\begin{cor}
If $A$ and $B$ are self-adjoint operators such that $A-B\in\bS_2$ and $f$ is a Lipschitz functions on $\R$, then $f(A)-f(B)\in\bS_2$ and
$$
\|f(A)-f(B)\|_{\bS_2}\le\|f\|_{\Lip}\|A-B\|_{\bS_2}.
$$
\end{cor}

\Pf It suffices to put $\bs{\frak f}=\dg_0f$. $\bl$

\medskip

We refer the reader to \cite{BS3} for the proof of Theorem \ref{LeS2}. Here we prove an analog of Theorem \ref{LeS2} in the case when $A-B$ is a bounded operator.

\begin{thm}
\label{dlyaogr}
Suppose that $A$ and $B$ are (not necessarily bounded) self-adjoint operators such that the operator $A-B$ is bounded. Let $\bs{\frak f}$ be a function on $\R^2$ such that 
$\bs{\frak f}\in\fM(E_A,E_B)$ and 
$\bs{\frak f}\big|_{\R^2\setminus\D}=\dg f\big|_{\R^2\setminus\D}$. 
Then formula {\em\rf{frk}} holds and
$$
\|f(A)-f(B)\|\le\|\bs{\frak f}\|_{\fM(E_A,E_B)}\|A-B\|.
$$
\end{thm}

\Pf Consider first the case when $A$ an $B$ are bounded operators. We have
\begin{align*}
\iint\bs{\frak f}(s,t)\,dE_A(s)(A-B)\,dE_B(t)&=
\iint\bs{\frak f}(s,t)\,dE_A(s)A\,dE_B(t)\\[.2cm]
&-
\iint\bs{\frak f}(s,t)\,dE_A(s)B\,dE_B(t).
\end{align*}
It is easy to see from the definition of double operator integrals in the Hilbert--Schmidt case that
$$
\iint\bs{\frak f}(s,t)\,dE_A(s)A\,dE_B(t)=
\iint s\bs{\frak f}(s,t)\,dE_A(s)\,dE_B(t)
$$
and
$$
\iint\bs{\frak f}(s,t)\,dE_A(s)B\,dE_B(t)=
\iint t\bs{\frak f}(s,t)\,dE_A(s)\,dE_B(t).
$$
Thus
$$
\iint\bs{\frak f}(s,t)\,dE_A(s)(A-B)\,dE_B(t)
=\iint(s-t)\bs{\frak f}(s,t)\,dE_A(s)\,dE_B(t).
$$
Clearly, $(s-t)\bs{\frak f}(s,t)=f(s)-f(t)$ for all $s,\,t\in\R$.
It follows that
\begin{align*}
\iint\bs{\frak f}(s,t)\,dE_A(s)(A-B)\,dE_B(t)&=
\iint f(s)\,dE_A(s)\,dE_B(t)\\[.2cm]
&-
\iint f(t)\,dE_A(s)\,dE_B(t)=f(A)-f(B).
\end{align*}

Suppose now that $A$ and $B$ are unbounded self-adjoint operators. 

Clearly, $f$ must be a Lipschitz function.
It follows easily that the domain of $f(A)$ contains the domain of $A$ and the same is true for the operator $B$. Hence, $f(A)-f(B)$ is a densely defined operator. Let us prove that it extends to a bounded operator and its extension (for which we keep the same notion $f(A)-f(B)$) satisfies \rf{frk}.

Consider the orthogonal projections 
$$
P_N\df E_A([-N,N])\quad\mbox{and}\quad Q_N\df E_B([-N,N])
$$
and define bounded self-adjoint operators $A_{[N]}$ and $B_{[N]}$ by
$$
A_{[N]}=P_NA\quad\mbox{and}\quad B_{[N]}=Q_NB.
$$
Obviously, 
$$
\lim_{N\to\be}P_N\left(\iint\bs{\frak f}(s,t)\,dE_A(s)(A-B)\,dE_B(t)\right)Q_N=
\iint\bs{\frak f}(s,t)\,dE_A(s)(A-B)\,dE_B(t)
$$
in the strong operator topology.

On the other hand, it is easy to see that
\begin{align*}
P_N&\left(\iint\bs{\frak f}(s,t)\,dE_A(s)(A-B)\,dE_B(t)\right)Q_N\\[.2cm]
&=
P_N\left(\iint\bs{\frak f}(s,t)\,dE_{A_{[N]}}(s)(A_N-B_N)\,dE_{B_{[N]}}(t)\right)Q_N
\\[.2cm]
&=P_N\big(f(A_{[N]})-f(B_{[N]})\big)Q_N
\end{align*}
because equality \rf{frk} holds for bounded self-adjoint operators. It remains to observe that
$$
P_N\big(f(A_{[N]})-f(B_{[N]})\big)Q_Nx=P_N\big(f(A)-f(B)\big)Q_Nx
~\to~(f(A)-f(B))x\quad\mbox{as}\quad N\to\be
$$
for all vectors $x$ in the dense subset $\bigcup_N\Range E_B([-N,N])$. $\bl$

\medskip

{\bf Remark 1.}
Suppose now that $\fI$ is a separable (or dual to separable) operator ideal of $\B(\h)$ equipped with a norm that makes it
a Banach space. Let $A$ and $B$ be self-adjoint operators such that $A-B\in\fI$ and let $f$ be a Lipschitz function on $\R$.
As above one can show that if the divided difference $\dg f$ can be extended to the diagonal $\D$ and the resulting function ${\frak f}$ on $\R^2$ belongs to the space
$\fM_\fI(E_A,E_B)$ of Schur multipliers of $\fI$, then formula \rf{frk} holds,
$$
f(A)-f(B)\in\fI\quad\mbox{and}\quad
\|f(A)-f(B)\|_\fI\le\|{\frak f}\|_{\fM_\fI(E_A,E_B)}\|A-B\|_\fI.
$$
We refer the reader to \cite{BS3} for more detail.

\medskip

{\bf Remark 2.}
Note also that the significance of formula \rf{fDK} is in the fact that it allows us to {\it linearize} the nonlinear problem of estimating $f(A)-f(B)$. Indeed, one can obtain desired estimates of $f(A)-f(B)$ by studying properties of the linear transformer
$$
T\mapsto\iint(\dg f)(s,t)\,dE_1(s)T\,dE_2(t).
$$

\medskip

{\bf Remark 3.} Similar results hold for functions of unitary operators, see \cite{BS3}. Analogs of the above results can also be obtained for analytic functions of contractions and of dissipative operators, see \cite{Pe0}, \cite{Pe7} and \cite{AP6}. However, in the case of contractions and in the case of dissipative operators one has to consider double operator integrals with respect to semi-spectral measures.

\

\section{Commutators and quasicommutators} 
\label{komutqua}

\

In the previous section we have seen that operator differences $f(A)-f(B)$ can be represented as double operator integrals with integrand equal to the divided difference
$\dg f$. Birman and Solomyak observed (see \cite{BS5}) that similar formula hold
for commutators $f(A)Q-Qf(A)$ and quasicommutators $f(A)Q-Qf(B)$.
The proof of the following result is practically the same as the proof of 
Theorem \ref{dlyaogr}.

\begin{thm}
\label{quasikomu}
Suppose that $A$ and $B$ are (not necessarily bounded) self-adjoint operators and $Q$ is a bounded linear operator such that the operator $AQ-QB$ is bounded. Let $\bs{\frak f}$ be a a function in $\fM_{\R,\R}$ such that 
$\bs{\frak f}\big|_{\R^2\setminus\D}=\dg f\big|_{\R^2\setminus\D}$. 
Then the quasicommutator $f(A)Q-Qf(B)$ is bounded,
\bay
\label{razquazi}
f(A)Q-Qf(B)=
\int\limits_\R\int\limits_\R\bs{\frak f}(s,t)\,dE_A(s)(AQ-QB)\,dE_B(t)
\ey
and 
$$
\|f(A)Q-Qf(B)\|\le\const\|{\frak f}\|_{\fM_{\R,\R}}\|AQ-QB\|.
$$
\end{thm}

Note that in the special case $A=B$ we obtain commutators $f(A)Q-Qf(A)$, while in the special cade $Q=I$ we obtain operator differences $f(A)-f(B)$.

Similar result holds in the case when $AQ-QB$ belongs to the Hilbert Schmidt class or other operator ideals.

In the rest of the paper we discuss in details estimates of $f(A)-f(B)$. Practically all the results are also valid for commutators and quasicommutators though we are not going to dwell on them.

\

\section{Operator Lipschitz functions} 
\label{oLoD}

\

In \S~\ref{razni} we have observed that if $f$ is a differentiable function on $\R$ such that the divided difference $\dg f$ belongs to the space of Schur multipliers 
$\fM_{\R,\R}$, then $f$ is operator Lipschitz. It turns out that the converse is also true. First of all, if $f$ is an operator Lipschitz function on $\R$, then $f$ is differentiable everywhere on $\R$ which was established in \cite{JW} (but not necessarily continuously differentiable: the function $x\mapsto x^2\sin(1/x)$ is operator Lipschitz, see \cite{KSh0}). On the other hand, it was shown in \cite{Pe2} (see also \cite{Pe3}) that if $f$ is a differentiable operator Lipschitz function, then $\dg f\in\fM_{\R,\R}$. Similar results hold for functions on the unit circle.

In this section we discuss some necessary conditions and sufficient conditions for a function to be operator Lipschitz. 

We start with necessary conditions for functions on the unit circle. The following result was established in \cite{Pe2}.

\begin{thm}
\label{B11}
Let $f$ be an operator Lipschitz function on $\T$. Then $f$ belongs to the Besov class
$B_1^1(\T)$.
\end{thm}

\Pf As we have discussed in \S~\ref{razni}, the divided difference
$$
(\z,\t)\mapsto\frac{f(\z)-f(\t)}{\z-\t}
$$
belongs to the space of Schur multipliers $\fM_{\T,\T}$. Trivially, this implies that the function
$$
(\z,\t)\mapsto\frac{f(\z)-f(\t)}{1-\ov{\t}\z}
$$
belongs to $\fM_{\T,\T}$. 

Consider the rank one operator $P$ on $L^2(\T)$ defined by 
$$
(P h)(\z)=\int_\T h(\t)\,d\m(\t).
$$
By Theorem \ref{tomSc}, the integral operator $\mC_f$ defined by
$$
(\mC_f h)(\z)=\int_\T\frac{f(\z)-f(\t)}{1-\ov{\t}\z} h(\t)\,d\m(\t)
$$
belongs to trace class $\bS_1$. An elementary calculation shows that
$$
\mC_f h=\pp_-fh_+-\pp_+fh_-,
$$
where $\pp_+$ is the orthogonal projection from $L^2$ onto the Hardy class $H^2$, $\pp_-$ is the 
orthogonal projection onto $H^2_-\df L^2\ominus H^2$, $h_+\df\pp_+ h$, and $h_-\df\pp_-h$
(see \cite{Pe5}, Ch. 1, \S~1). It is easy to see that both Hankel operators $H_f$ and 
$H_{\ov{f}}$ belong to $\bS_1$. Recall that the {\it Hankel operator} 
$H_f:H^2\to H^2_-$ is defined by $H_f\f=\pp_-f\f$.

By the trace class criterion for Hankel operators \cite{Pe1} (see also \cite{Pe5}, 
Ch. 6, \S~1), we find that $f\in B_1^1(\T)$. $\bl$

The following stronger necessary condition was also obtained in \cite{Pe2} by using 
the trace class criterion for Hankel operators.

\begin{thm}
\label{snu}
Let $f$ be an operator Lipschitz function on $\T$. Then the Hankel operators
$H_f$ and $H_{\ov{f}}$ map the Hardy class $H^1$ into the Besov space $B_1^1(\T)$.
\end{thm}

Note that Semmes observed (see the proof in \cite{Pe0} that the Hankel operators $H_f$ and $H_{\ov{f}}$ map $H^1$ into $B_1^1(\T)$ if and only if the measure $\mu$ defined by
$$
d\mu(\z)=\big(|(f_+)''(\z)|+|(\ov{f}_+)''(\z)|\big)\,dm_2(\z)
$$
is a Carleson measure on the unit disk $\dd$.

Theorem \ref{B11} implies easily that a continuously differentiable function on $\T$ does not have to be operator Lipschitz. Indeed, it follows from \rf{Bperf} that the lacunary Fourier coefficients of the derivative of a function $f$ in $B_1^1(\T)$ must satisfy the condition
$$
\sum_{k\ge0}|\widehat{f'}(2^k)|<\be,
$$
while it is well known that an arbitrary sequence in $\ell^2$ can be the sequence of lacunary coefficients of the derivative of a continuously differentiable function.

The above results can be extended to functions on $\R$. The analog of Theorem \ref{B11}
is that if $f$ is an operator Lipschitz function on $\R$, then $f$ belongs to the Besov space $B_1^1$ locally. An analog of Theorem \ref{snu} also holds as well as the characterization of the last necessary condition in terms of Carleson measures. We refer the reader to \cite{Pe3} and \cite{Pe0}.

We proceed now to sufficient conditions for operator Lipschitzness. 
The following result was obtained in \cite{Pe2}. 

\begin{thm}
\label{dosu}
Let $f$ be a function on $\T$ of Besov class $B_{\be,1}^1(\T)$. Then $f$ is operator Lipschitz.
\end{thm}

We give here an idea of the proof of Theorem \ref{dosu}. It is easy to see that it suffices to prove the following inequality: {\it suppose that $\f$ is an analytic polynomial (i.e., a polynomial of $z$) of degree $m$, then the norm of $\dg\f$ in the projective tensor product $C(\T)\hat\otimes C(\T)$ admits the following estimate:}
\bay
\label{nerB}
\|\dg\f\|_{C(\T)\hat\otimes C(\T)}\le\const m\|\f\|_{C(\T)}.
\ey
Note that the projective tensor product $C(\T)\hat\otimes C(\T)$ can be defined in the same way as the projective tensor product of $L^\be$ spaces.

It can easily be verified that
$$
(\dg\f)(z_1,z_2)=\sum_{j,k\ge0}\hat\f(j+k+1)z_1^jz_2^k.
$$
Clearly,
\begin{align*}
\sum_{j,k\ge0}\hat\f(j+k+1)z_1^iz_2^j&=
\sum_{k=0}^{m-1}\left(\sum_{j\ge0}\a_{jk}\hat\f(j+k+1)z_1^j\right)z_2^k\\[.2cm]
&+
\sum_{j=0}^{m-1}\left(\sum_{k\ge0}\b_{jk}\hat\f(j+k+1)z_2^k\right)z_1^j,
\end{align*}
where
$$
\a_{jk}=\left\{\begin{array}{ll}\frac12,&j=k=0,\\[.2cm]
\frac{j}{j+k},& j+k\ne0,
\end{array}\right.
\quad
\mbox{and}\quad
\b_{jk}=\left\{\begin{array}{ll}\frac12,&j=k=0,\\[.2cm]
\frac{k}{j+k},& j+k\ne0.
\end{array}\right.
$$
It can be shown that 
$$
\left\|\sum_{j\ge0}\a_{jk}\hat\f(j+k+1)z^j\right\|_{C(\T)}\le\const\|\f\|_{C(\T)}
$$
and
$$
\left\|\sum_{k\ge0}\b_{jk}\hat\f(j+k+1)z^j\right\|_{C(\T)}\le\const\|\f\|_{C(\T)}
$$
(see \cite{Pe2} and \cite{Pe4} for details).
This implies easily inequality \rf{nerB}.

A similar fact holds for functions on $\R$. The following result was obtained in \cite{Pe4}.

\begin{thm}
\label{neprana}
If $f$ belongs to the Besov class 
$B_{\be,1}^1(\R)$, then $f$ is an operator Lipschitz function on $\R$.
\end{thm}

It follows from the definition of $B_{\be,1}^1(\R)$ (see \rf{Wn}) that to prove Theorem \ref{neprana}, it suffices to establish the following fundamental inequality: 
\bay
\label{fundrr}
\|\dg f\|_{\fM_{\R,\R}}\le\const\s\|f\|_{L^\be}
\ey
for an arbitrary bounded function $f$ on $\R$ with Fourier transform supported in $[-\s,\s]$.
Inequality \rf{fundrr} together with formula \rf{fDK} implies that 
\bay
\label{fundner}
\|f(A)-f(B)\|\le\const\s\|f\|_{L^\be}\|A-B\|
\ey
for arbitrary self-adjoint operators $A$ and $B$ with bounded $A-B$ and for an arbitrary function $f$ in $L^\be(\R)$ whose Fourier transform is supported in $[-\s,\s]$. In \cite{AP4} it was shown that inequality \rf{fundner} holds with constant 1 on the right.

To prove inequality \rf{fundrr} we introduce the functions $r_u$, $u>0$, whose Fourier transforms $\F r_u$ are defined by
$$
(\F r_u)(s)=\left\{\begin{array}{ll}1,&|s|\le u,\\[.2cm]\frac{u}{|s|},&|s|>u.
\end{array}\right.
$$
It is easy to show that $r_u\in L^1(\R)$ and $\|r_u\|_{L^1}\le\const$. It follows that the function $1-r_u$ is the Fourier transform of finite signed measure. We denote this measure by $\mu_u$. We have
$$
(\F\mu_u)(s)=\left\{\begin{array}{ll}0,&|s|\le u,\\[.2cm]\frac{|s|-u}{|s|},&|s|>u.
\end{array}\right.
$$

To prove inequality \rf{fundrr}, we establish the following integral representation for the divided difference $\dg f$:

\begin{lem}
\label{ifrr}
Let $f$ be a bounded function on $\R$ whose Fourier transform has compact support in 
$[0,\be)$. Then the following representation holds:
\begin{align}
\label{tozhrr}
(\dg f)(s,t)={\rm i}\int_{\R_+}\big(f*\mu_u)(s)e^{-{\rm i}su}e^{{\rm i}tu}\,du
+{\rm i}\int_{\R_+}\big(f*\mu_u)(t)e^{-{\rm i}tu}e^{{\rm i}su}\,du.
\end{align}
\end{lem}

To prove identity \rf{tozhrr}, we can first consider the special case when $f$ is the Fourier transform of an $L^1$ function, in which case this is an elementary exercise, and then consider suitable approximation, see \cite{Pe3} and \cite{Pe4} for details.

\begin{cor}
\label{orr0s}
Let $f$ be a bounded function on $\R$ whose Fourier transform is supported in $[0,\s]$.
Then inequality {\em\rf{fundrr}} holds.
\end{cor}

\Pf Clearly, $f*\mu_u=\0$ for $u>\s$. Representation \rf{tozhrr} gives us the following estimates:
$$
\|\dg f\|_{L^\be\hat\otimes_{\rm i}L^\be}
\le2\int_0^\s\|f*\mu_u\|_{L^\be}\,du\le
2\int_0^\s(\|f\|_{L^\be}+\|f*r_u\|_{L^\be})\,du
\le\const\s\|f\|_{L^\be}.\quad\bl
$$

To prove inequality \rf{fundrr} in the general case we represent $f$ as the sum of \lb$f_n=f*W_n$ (see \rf{fn}). Consider now the function $(f_n)_+$ whose Fourier transform is equal to $\chi_{\R_+}\F f_n$. It remains to observe that 
$\|(f_n)_+\|_{L^\be}\le\const\|f_n\|_{L^\be}$, see \cite{Pe3} and \cite{Pe4} for details.

\medskip

{\bf Proof of Theorem \ref{neprana}.} By \rf{fundrr}, we have
$$
\|\dg f\|_{\fM_{\R,\R}}\le\sum_{n\in\Z}\|\dg f_n\|_{\fM_{\R,\R}}
\le\const\sum_{n\in\Z}2^n\|f_n\|_{L^\be},
$$
where, as usual, $f_n\df f*W_n$. $\bl$
\medskip

Note that inequality \rf{fundner} and its version for Schatten--von Neumann norms will play a very important role in H\"older type inequalities, in Schatten--von Neumann estimates of operator differences, see Sections
\ref{oHmc} and \ref{SchvNe}.

To conclude the section, I would like to mention that similar results also hold for functions of contractions and functions of dissipative operators, see \cite{Pe0},
\cite{KSh4} and 
\cite{AP6}.

\

\section{Operator differentiable functions} 
\label{oDf}

\

In the previous section we have shown that the condition $f\in B_{\be,1}^1(\R)$ is sufficient for $f$ to be operator Lipschitz on $\R$.
It turns out that the same condition $f\in B_{\be,1}^1(\R)$ is also sufficient for operator differentiability. 

\medskip

{\bf Definition.}
A function $f$ on $\R$ is called {\it operator differentiable} if the 
limit
$$
\lim_{\dt\to0}\dt^{-1}\big(f(A+\dt K)-f(A)\big)
$$
exists in the operator norm for an arbitrary self-adjoint operator $A$ and an arbitrary
bounded self-adjoint operator $K$.

\medskip

The following result can be found in \cite{Pe3} and \cite{Pe4}.

\begin{thm}
\label{diffe}
Let $f$ be a function in $B_{\be,1}^1(\R)$. Then $f$ is operator differentiable and
\bay
\label{DalKr}
\lim_{t\to0}t^{-1}\big(f(A+t K)-f(A)\big)=\iint\frac{f(s_1)-f(s_2)}{s_1-s_2}\,dE_A(s_1)K\,dE_A(s_2).
\ey
whenever $A$ is a self-adjoint operator and $K$ is a bounded self-adjoint operator.
\end{thm}

Formula \rf{DalKr} is called the {\it Daletskii--Krein formula}. It was established in \cite{DK} under considerably stronger assumptions. Later Birman and Solomyak proved in \cite{BS3}
formula \rf{DalKr} under less restrictive assumptions.

Let me give an idea of the proof of Theorem \ref{diffe}. Under the hypotheses of the theorem, we have
$$
f(A+t K)-f(A)=t\iint\frac{f(s_1)-f(s_2)}{s_1-s_2}\,dE_{A_t}(s_1)K\,dE_A(s_2),
$$
where $A_t\df A+t K$ (see \S~\ref{razni}). To establish formula \rf{DalKr}, we can represent the divided difference $\dg f$ as an element of the integral projective tensor product
$$
(\dg f)(s_1,s_2)=\int\f_x(s_1)\psi_x(s_2)\,d\s(x),
$$
where $\s$ is a $\s$-finite measure and
$$
\int\|\f_x\|_{L^\be}\|\psi_x\|_{L^\be}\,d\s(x)<\be
$$
Moreover, the functions $\f_x$ satisfy the following:
$$
\lim_{t\to0}\|\f_x(A+t K)-\f_x(A)\|=0
$$
for all $x$. This can be deduced easily from Lemma \ref{ifrr}, see \cite{Pe3} for details.

We have
$$
\iint\frac{f(s_1)-f(s_2)}{s_1-s_2}\,dE_{A_t}(s_1)K\,dE_A(s_2)=
\int\f_x(A_t)K\psi_x(A)\,d\s(x).
$$
Clearly, the above conditions easily imply that
\begin{align*}
\lim_{t\to0}\int\f_x(A_t)K\psi_x(A)\,d\s(x)&=\int\f_x(A)K\psi_x(A)\,d\s(x)\\[.2cm]
&=
\iint\frac{f(s_1)-f(s_2)}{s_1-s_2}\,dE_A(s_1)K\,dE_A(s_2)
\end{align*}
which implies \rf{diffe}.

\medskip

{\bf Remark 1.} The problem of differentiability of the function $t\mapsto f(A+tK)$ is the problem of the existence of the G\^ateaux derivative of the map
$K\mapsto f(A+K)-f(A)$ defined on the real space of bounded self-adjoint operators.
We have proved that this map is differentiable in the sense of G\^ateaux for functions $f$ in $B_{\be,1}^1(\R)$. The reasoning given above allows one to prove that under the same assumptions this map is differentiable in the sense of Fr\'echet and the differential of this map is the double operator integral
$$
\iint\frac{f(s_1)-f(s_2)}{s_1-s_2}\,dE_A(s_1)K\,dE_A(s_2).
$$

\medskip

{\bf Remark 2.} 
The above argument shows that under the hypotheses of Theorem \ref{diffe}, the function $t\mapsto f(A_t)$ is actually continuously differentiable in the operator norm.

\

Finally, I would like to mention that similar results hold for functions of unitary operators, functions of contractions and functions of dissipative operators, see \cite{Pe2}, \cite{Pe7}  and \cite{AP6}. In particular, in the case of functions of unitary operators we can consider the problem of differentiability of the function 
$t\mapsto f\big(e^{{\rm i}tA}U)$, $t\in\R$, where $f$ is a function on the unit circle $\T$, $U$ is a unitary operator and $A$ is a bounded self-adjoint operator. It was proved in \cite{Pe2} that under the assumption $f\in B_{\be,1}^1(\T)$, the function $t\mapsto f\big(e^{{\rm i}tA}U)$ is differentiable in the operator norm and its derivative is equal to the following double operator integral:
$$
{\rm i}\left(\iint\frac{\f(\z)-\f(\tau)}{\z-\tau}\,dE_U(\z)A\,dE_U(\tau)\right)U.
$$

\

\section{The Lifshits--Krein trace formula}
\label{fsLK}

\

The notion of the spectral shift function was introduced be I.M. Lifshits in \cite{L}.
He discovered in that paper a trace formula for $f(A)-f(B)$ where $A$ is the initial operator and $B$ is a perturbed operator that involves the spectral shift function.
Later M.G. Krein in \cite{Kr} generalized the trace formula to a considerably more general situation when $A$ is an arbitrary self-adjoint operator and $B$ is a trace class perturbation of $A$. 

Let $A$ be a self-adjoint operator on Hilbert space and let $B$ be a perturbed self-adjoint operator with $A-B\in\bS_1$. It was shown in \cite{Kr} that there exists a unique real function $\xi$ in $L^1(\R)$ such that
\bay
\label{LKss}
\trace\big(f(B)-f(A)\big)=\int_\R f'(s)\xi(s)\,ds,
\ey
whenever $f$ is a differentiable function on $\R$ whose derivative is the Fourier transform of an $L^1$ function. The function $\xi$ is called the {\it spectral shift function} associated with the pair $(A,B)$.

Moreover, it was shown in \cite{Kr} that under the same assumptions
$$
\|\xi\|_{L^1}=\|B-A\|_{\bS_1}\quad\mbox{and}\quad\int_\R\xi(s)\,ds=\trace(B-A).
$$

The right-hand side of formula \rf{LKss} is well defined for an arbitrary Lipschitz function $f$. Krein asked in \cite{Kr} whether formula \rf{LKss} holds for an arbitrary Lipschitz function $f$. It turns out, however, that the Lipschitzness of $f$ does not imply that the operator $f(A)-f(B)\in\bS_1$ whenever $A-B\in\bS_1$. 
This was first observed in \cite{F0}. 

On the other hand, it can be shown that a function $f$ preserves trace class perturbations, i.e.,
\bay
\label{soyav}
A-B\in\bS_1\quad\Longrightarrow\quad f(A)-f(B)\in\bS_1
\ey
if and only if $f$ is operator Lipschitz (the operators $A$ and $B$ do not have to be bounded). This implies that the necessary conditions for operator Lipschitzness mentioned in \S~\ref{oLoD} are also necessary for property \rf{soyav}.

On the other hand, it was proved in \cite{Pe3} that the condition $f\in B_{\be,1}^1(\R)$ sufficient for operator Lipschitzness (see \S~\ref{oLoD}) is also sufficient for 
trace formula \rf{LKss} to hold.

In this section we use the Birman--Solomyak approach \cite{BS6} that is based on double operator integrals. Actually, their approach allows to prove the existence of a finite real signed Borel measure $\nu$ such that
\bay
\label{slednu}
\trace\big(f(A)-f(B)\big)=\int_\R f'(s)\,d\nu(s),\quad\mbox{and}\quad
\|\nu\|\le\|T\|_{\bS_1}
\ey
for sufficiently nice functions $f$. It follows from the results of \cite{Kr} that $\nu$ is absolutely continuous with respect to Lebesgue measure and $d\nu=\xi\,d\m$, where $\xi$ is the spectral shift function. Moreover, we combine the Birman--Solomyak approach with the observation that for $f\in B_{\be,1}^1(\R)$, the function $t\mapsto f(A+t(B-A))$, $t\in\R$, is continuously differentiable in the trace norm (this can be proved in the same way as Theorem \ref{DalKr}) and the Daletskii--Krein formula holds for the derivative of this operator function. This allows us to prove the following extension of the Birman--Solomyak result:

\begin{thm}
\label{orBS}
Let $A$ and $B$ be self-adjoint operators on Hilbert space such that $B-A\in\bS_1$.
Then there exists a real signed Borel measure $\nu$ on $\R$ such that formula
{\em\rf{slednu}} holds for an arbitrary function $f$ in
$B_{\be,1}^1(\R)$. 
\end{thm}

Recall that under the hypotheses of the theorem $f(A)-f(B)\in\bS_1$.

\Pf
Put $A_t\df A+t(B-A)$. Then the function $t\mapsto f(A_t)$ is differentiable in the trace norm and
$$
\frac{d}{dt}f(A_t)\Big|_{t=u}=\iint(\dg f)(s_1,s_2)\,dE_{A_u}(s_1)(B-A)\,dE_{A_u}(s_2).
$$
This can be proved in exactly the same way as Theorem \ref{DalKr}.

Since $\dg f$ is a Schur multiplier (see \S~\ref{oLoD}), it follows that
$$
\iint(\dg f)(s_1,s_2)\,dE_{A_u}(s_1)(B-A)\,dE_{A_u}(s_2)\in \bS_1.
$$
By Theorem \ref{sled} , we have
\begin{align*}
\trace\left(\iint(\dg f)(s_1,s_2)\,dE_{A_u}(s_1)(B-A)\,dE_{A_u}(s_2)\right)
=\int f'(s)\,d\nu_u(s),
\end{align*}
where the signed measure $\nu_s$ is defined on the Borels sets by
$$
\nu_u(\D)=\trace\big(E_{A_u}(\D)T\big).
$$
Clearly, $\|\nu_u\|\le\|T\|_{\bS_1}$.

It is easy to verify that
\begin{align*}
\trace\big(f(A)-f(B)\big)=
\int\trace\left(\frac{d}{dt}f(A_t)\Big|_{t=u}\right)\,du=\int f'(s)\,d\nu(s),
\end{align*}
where the signed measure $\nu$ is defined by
$$
\nu=\int_0^1\nu_u\,du.
$$
Note that the function $u\mapsto\nu_u$ is continuous in the space of measures
equipped with the weak-$*$ topology, and so integration makes sense.
Clearly, $\|\nu\|\le\|T\|_{\bS_1}$. $\bl$

We have already mentioned that $d\nu=\xi\,d\m$, where $\xi$ is the spectral shift function. This implies the following extension of the Krein theorem.

\begin{thm}
\label{otKr}
Let $A$ and $B$ be self-adjoint operators such that $B-A\in\bS_1$. Suppose that
$f\in B_{\be,1}^1(\R)$. Then trace formula {\em\rf{LKss}} holds.
\end{thm}

The original proof of Theorem \ref{otKr} by a different method was obtained in \cite{Pe3}.

\

\section{Operator H\"older functions. Arbitrary moduli of continuity}
\label{oHmc}

\

In this section we obtain norm estimates for $f(A)-f(B)$, where $A$ and $B$ are self-adjoint operators and $f$ is a H\"older function of order $\a$, $0<\a<1$. Then we consider the more general problem of estimating $f(A)-f(B)$ in terms of the modulus of continuity of $f$.

By analogy with the notion of operator Lipschitz functions. We say that a function $f$ on $\R$ is {\it operator H\"older of order} $\a$, $0<\a<1$, if
$$
\|f(A)-f(B)\|\le\const\|A-B\|^\a.
$$
The problem of whether a H\"older function of order $\a$ (recall that the class of such functions is denoted by $\L_\a(\R)$, see \S~\ref{prel}) is necessarily
operator H\"older of order $\a$ remained open for 40 years and it was solved in \cite{AP1} (see also \cite{AP2} for a detailed presentation). The solution is given by the following theorem:

\begin{thm}
\label{Hoa}
Let $\a\in(0,1)$. Then 
$$
\|f(A)-f(B)\|\le\const(1-\a)^{-1}\|A-B\|^\a,
$$
whenever $A$ and $B$ are self-adjoint operators with bounded $A-B$.
\end{thm}

Thus the term ``an operator H\"older function of order $\a$" turns out to be short-lived.

We prove here Theorem \ref{Hoa} for bounded self-adjoint operators and refer the reader to
\cite{AP4} for details how to treat the case of unbounded operators.

\medskip

{\bf Proof of Theorem \ref{Hoa}.}
Let $N$ be an integer. Then $f(A)-f(B)$ admits a representation 
\bay
\label{f}
f(A)\!-\!f(B)=\!\sum_{n=-\be}^N\!\big(f_n(A)\!-\!f_n(B)\big)\!+\!\big((f-f*V_N)(A)\!-\!(f-f*V_N)(B)\big)
\ey
and the series converges absolutely in the operator norm. Here $f_n=f*W_n$ (see \rf{fn}) and $V_N$ is the de la Vall\'ee Poussin type kernel defined by \rf{dlVP}. Suppose that $M<N$. It is easy to see that
\begin{align*}
f(A)-f(B)&-
\left(\sum_{n=M+1}^N\big(f_n(A)-f_n(B)\big)+\big((f-f*V_N)(A)-(f-f*V_N)(B)\big)\right)\\[.2cm]
&=(f*V_M)(A)-(f*V_M)(B).
\end{align*}
Clearly, the Fourier transform of $f*V_M$ is is supported in 
$\big[-2^{M+1},2^{M+1}\big]$. Thus it follows from fundamental inequality
\rf{fundner} that for $M\le0$,
\begin{align*}
\big\|(f*V_M)(A)-(f*V_M)(B)\big\|&\le\const2^M\|f*V_M\|_{L^\be}\|A-B\|\\[.2cm]
&\le\const2^M\|f*V_0\|_{L^\be}\|A-B\|
\to0\quad\mbox{as}\quad M\to-\be. 
\end{align*}

Suppose  now that $N$ is the integer satisfying 
$$
2^{-N}<\|A-B\|\le2^{-N+1}.
$$
By \rf{dlVPa}, we have
\begin{align*}
\big\|(f-f*V_N)(A)&-(f-f*V_N)(B)\big\|\le2\|f-f*V_N\|_{L^\be}\\[.2cm]
&\le\const\|f\|_{\L_\a(\R)}2^{-N\a}\le\const\|f\|_{\L_\a(\R)}\|A-B\|^\a.
\end{align*}
On the other hand, it follows from fundamental inequality \rf{fundner} and from \rf{Wn} that
\begin{align*}
\sum_{n=-\be}^N\|f_n(A)-f_n(B)\|&\le\const\sum_{n=-\be}^N2^n\|f_n\|_{L^\be}\|A-B\|\\[.2cm]
&\le\const\sum_{n=-\be}^N2^n2^{-n\a}\|f\|_{\L_\a(\R)}\|A-B\|\\[.2cm]
&=\const2^{N(1-\a)}\sum_{k\ge0}2^{-k(1-\a)}\|f\|_{\L_\a(\R)}\|A-B\|\\[.2cm]
&=\const2^{N(1-\a)}(1-\a)\|f\|_{\L_\a(\R)}\|A-B\|\\[.2cm]
&\le\const(1-\a)\|f\|_{\L_\a(\R)}\|A-B\|^\a.\quad\bl
\end{align*}

Suppose now that $\o$ is an arbitrary {\it modulus of continuity}, i.e., $\o$ is a continuous nondecreasing function on $[0,\be)$ such that $\o(s+t)\le\o(s)+\o(t)$, $s,\,t\ge0$á and $\o(0)=0$.
We associate with $\o$ the function $\o_*$ defined by
$$
\o_*(x)=x\int_x^\be\frac{\o(t)}{t^2}\,dt=\int_1^\infty\frac{\omega(sx)}{s^{2}}ds,\quad x>0.
$$
It is easy to see that if $\o_*(x)<\be$ for some $x>0$, then $\o_*(x)<\be$ for all $x>0$ in which case $\o_*$ is also a modulus of continuity.

The following result was obtained in \cite{AP1} and \cite{AP2}.

\begin{thm}
\label{modne}
Let $\o$ be a modulus of continuity. Then
$$
\|f(A)-f(B)\|\le\const\o_*\big(\|A-B\|\big)
$$
for arbitrary self-adjoint operators $A$ and $B$ with bounded $A-B$.
\end{thm}

The proof of Theorem \ref{modne} is similar to the proof of Theorem \ref{Hoa}. 

Slightly weaker results were obtained independently in \cite{FN2}.

Theorem \ref{modne} implies the following result proved in \cite{AP2}:

\begin{cor}
\label{FaN}
Suppose that $A$ and $B$ be self-adjoint operators with spectra in an interval $[a,b]$. Then for a continuous function $f$ on $[a,b]$ the following inequality holds:
$$
\|f(A)-f(B)\|\le\const\,\log\left(e\,\frac{b-a}{\|A-B\|}\right)
\,\o_f\big(\|A-B\|\big).
$$
\end{cor}

Theorem \ref{FaN}
improves earlier estimates obtained in \cite{F3}.

We refer the reader to \cite{AP5} for more detailed information and more sophisticated estimates of $f(A)-f(B)$.

Note that similar results hold for functions of unitary operators, contractions and dissipative operators, see \cite{AP2} and \cite{AP6}.

\

\section{Schatten--von Neumann estimates of operator differences}
\label{SchvNe}

\

In this section we list several results on estimates of the norms of $f(A)-f(B)$ in operator ideals and, in particular, in Schatten--von Neumann classes.

Fundamental inequality \rf{fundrr} together with formula \rf{fDK} allows us to use Mityagin's interpolation theorem \cite{Mi} to generalize and generalize inequality 
\rf{fundner}
to arbitrary separable (or dual to separable) ideals $\fI$:
\bay
\label{obid}
\|f(A)-f(B)\|_\fI\le\const\s\|f\|_{L^\be}\|A-B\|_\fI
\ey
for arbitrary self-adjoint operators $A$ and $B$ with bounded $A-B$ and for an arbitrary bounded function $f$ on $\R$ whose Fourier transform is supported in $[-\s,\s]$.

In particular, inequality \rf{obid} holds in the case $\fI=\bS_p$, $p\ge1$.

This implies the following result (see \cite{Pe2} and \cite{Pe3}).

\begin{thm}
\label{idealy}
Let $\fI$ be a separable (symmetrically normed) operator ideal or an operator ideal dual to separable and let $f$ be a function in the Besov class $B_{\be,1}^1(\R)$. Suppose that 
$A$ and $B$ are self-adjoint operators such that $A-B\in\fI$. Then $f(A)-f(B)\in\fI$ and
$$
\|f(A)-f(B)\|_\fI\le\const\|f\|_{B_{\be,1}^1}\|A-B\|_\fI.
$$
\end{thm}

In the case when $\fI=\bS_p$, $1<p<\be$, Theorem \rf{idealy} was improved significantly in \cite{PoS}:

\begin{thm}
\label{PoSu}
Let $1<p<\be$ and let $f$ be a Lipschitz function on $\R$. Suppose that  $A$ and $B$ are self-adjoint operators such that $A-B\in\bS_p$. Then $f(A)-f(B)\in\bS_p$ and
$$
\|f(A)-f(B)\|_{\bS_p}\le c_p\|f\|_{\Li}\|A-B\|_{\bS_p},
$$
where $c_p$ is a positive number that depends only on $p$.
\end{thm}

We proceed now to estimating Schatten--von Neumann norms of $f(A)-f(B)$ for functions $f$ in the H\"older class $\L_\a(\R)$, $0<\a<1$. For a nonnegative integer $l$
and for $p\ge1$, we define the following norm on the space of bounded linear operators on Hilbert space:
$$
\|T\|_{\bS_p^l}\df\left(\sum_{j=0}^l (s_j(T))^p\right)^{1/p},
$$
where $s_j(T)$ is the $j$th singular value of $T$.

The following result obtained in \cite{AP3} is crucial.

\begin{thm}
\label{singSp}
Let $0<\a<1$. Then there exists a positive number $c>0$ such that for every
$l\ge0$,  $p\in[1,\be)$,  $f\in\L_\a(\R)$, and for arbitrary self-adjoint operators $A$ and $B$ on Hilbert space with bounded $A-B$, the following inequality holds:
$$
s_j\big(f(A)-f(B)\big)\le c\,\|f\|_{\L_\a(\R)}(1+j)^{-\a/p}\|A-B\|_{\bS_p^l}^\a
$$
for every $j\le l$.
\end{thm}

\Pf Put $f_n\df f*W_n$, $n\in\Z$, and fix an integer $N$. We have
by \rf{obid} and \rf{Wn},
\begin{align*}
\left\|\sum_{n=-\be}^N\big(f_n(A)-f_n(B)\big)\right\|_{\bS_p^l}
&\le\sum_{n=-\be}^N\big\|f_n(A)-f_n(B)\big\|_{\bS_p^l}\\[.2cm]
&\le\const\sum_{n=-\be}^N2^n\|f_n\|_{L^\be}\|A-B\|_{\bS_p^l}\\[.2cm]
&\le\const\|f\|_{\L_\a(\R)}\sum_{n=-\be}^N2^{n(1-\a)}\|A-B\|_{\bS_p^l}\\[.2cm]
&\le\const2^{N(1-\a)}\|f\|_{\L_\a(\R)}\|A-B\|_{\bS_p^l}.
\end{align*}
On the other hand,
\begin{align*}
\left\|\sum_{n>N}\big(f_n(A)-f_n(B)\big)\right\|
&\le2\sum_{n>N}\|f_n\|_{L^\be}\\[.2cm]
&\le\const\|f\|_{\L_\a(\R)}\sum_{n>N}2^{-n\a}
\le\const2^{-N\a}\|f\|_{\L_\a(\R)}.
\end{align*}
Put
$$
R_N\df\sum_{n=-\be}^N\big(f_n(A)-f_n(B)\big)
\qm Q_N\df\sum_{n>N}\big(f_n(A)-f_n(B)\big).
$$
Clearly, for $j\le l$,
\begin{align*}
s_j\big(f(A)-f(B)\big)&\le s_j(R_N)+\|Q_N\|
\le(1+j)^{-1/p}\|f(A)-f(B)\|_{\bS_p^l}+\|Q_N\|\\[.2cm]
&\le\const\left((1+j)^{-1/p}2^{N(1-\a)}\|f\|_{\L_\a(\R)}\|A-B\|_{\bS_p^l}
+2^{-N\a}\|f\|_{\L_\a(\R)}\right).
\end{align*}
To obtain the desired estimate, it suffices to choose the number $N$ so that
$$
2^{-N}<(1+j)^{-1/p}\|A-B\|_{\bS_p^l}\le2^{-N+1}.\quad\bl
$$

The following result can be deduced from Theorem \ref{singSp}. We refer the reader to \cite{AP3} for details.

\begin{thm}
\label{Scha}
Let $0<\a<1$ and $1<p<\be$ and let
$f\in\L_\a(\R)$. Supposed that $A$ and $B$ are self-adjoint operators such that $A-B\in\bS_p$. Then the operator $f(A)-f(B)$ belongs to $\bS_{p/\a}$ and
$$
\big\|f(A)-f(B)\big\|_{\bS_{p/\a}}\le c_{\a,p}\,\|f\|_{\L_\a(\R)}\|A-B\|^\a_{\bS_p},
$$
where $c_{\a,p}$ depends only on $\a$ and $p$.
\end{thm}

Note that for $p=1$, the conclusion of Theorem \ref{Scha} does not hold, see \cite{AP3}. The example given in \cite{AP3} is based on the $S_p$ criterion for Hankel operators, see \cite{Pe1} and \cite{Pe5}.

Nevertheless, the conclusion of Theorem \ref{Scha} can be obtained under stronger assumptions on $f$. The following result was obtained in \cite{AP3}.

\begin{thm}
\label{S1s}
Let $0<\a\le1$ and let $f$ be a function in the Besov class
$B_{\be1}^\a(\R)$.  Supposed that $A$ and $B$ are self-adjoint operators such that 
$A-B\in\bS_1$. Then the operator $f(A)-f(B)$ belongs to $\bS_{1/\a}$ and
$$
\big\|f(A)-f(B)\big\|_{\bS_{1/\a}}\le c_\a\,\|f\|_{B_{\be1}^\a(\R)}
\|A-B\|_{\bS_1}^\a.
$$
where $c_{\a}$ depends only on $\a$.
\end{thm}

Note that in \cite{AP3} the above results were generalized to the case of considerably more general operator ideals.

\medskip

{\bf Remark.}
As before, I would like to mention that similar results hold for functions of unitary operators, contractions and dissipative operators, see \cite{AP3} and \cite{AP6}.

\

\section{Functions of normal operators}
\label{fnorm}

\

We proceed now to the study of functions of normal operators under perturbation. 
The results of this section were obtained in \cite{APPS}. Earlier weaker results were obtained in \cite{FN1}.

The spectral theorem allows us to define functions of a normal operator as integrals with respect to its spectral measure:
$$
f(N)=\int_\C \z\,dE_N(\z).
$$
Here $N$ is a normal operator and $E_N$ is its spectral measure.

We are going to study estimates of $f(N_1)-f(N_2)$ in terms of $N_1-N_2$. 

As in the case of functions of self-adjoint operators we can consider the divided difference
$$
(\dg f)(\z_1,\z_2)=\frac{f(\z_1)-f(\z_2)}{\z_1-\z_2},\quad\z_1,~\z_2\in\C,
$$
and prove the formula
\bay
\label{nordd}
f(N_1)-f(N_2)=\iint(\dg f)(\z_1,\z_2)\,dE_1(\z_1)(N_1-N_2)\,dE_2(\z_2),
\ey
whenever $N_1$ and $N_2$ are normal operators with $N_1-N_2\in\bS_2$ and $f$ is a Lipschitz function. Again, it does not matter how we define $\dg f$ on the diagonal of 
$\C\times\C$.

If $f$ is a function on $\C$ such that the divided difference $\dg f$ can be extended to the diagonal and the extension belongs to the space of Schur multipliers $\fM_{\C,\C}$, then formula \rf{nordd} as soon as $N_1$ and $N_2$ are normal operators with bounded difference. Moreover, such functions are necessarily {\it operator Lipschitz}, i.e.,
$$
\|f(N_1)-f(N_2)\|\le\const\|N_1-N_2\|.
$$
The trouble is that such functions are necessarily linear which follows from the results of \cite{JW}.

In \cite{APPS} we used the following representation for $f(N_1)-f(N_2)$:
\begin{align}
\label{irnoro}
f(N_1)-f(N_2)&=\iint\limits_{\C^2}\big(\dg_yf\big)(z_1,z_2)\,
dE_1(z_1)(B_1-B_2)\,dE_2(z_2)\nonumber\\[.2cm]
&+\iint\limits_{\C^2}\big(\dg_xf\big)(z_1,z_2)\,
dE_1(z_1)(A_1-A_2)\,dE_2(z_2),
\end{align}
where 
$$
A_j=\re N_j,\quad B_j=\im N_j,\quad x_j=\re z_j,\quad y_j=\im z_j,\quad j=1,~2,
$$
and the divided differences $\dg_xf$ and $\dg_yf$
are defined by
$$
\big(\dg_xf\big)(z_1,z_2)\df\frac{f(x_1,y_2)-f(x_2,y_2)}{x_1-x_2},
\quad z_1,\,z_2\in\C,
$$
and
$$
\big(\dg_yf\big)(z_1,z_2)\df\frac{f(x_1,y_1)-f(x_1,y_2)}{y_1-y_2},
\quad z_1,\,z_2\in\C.
$$

It was established in \cite{APPS} that for a function $f$ in the Besov class
$B_{\be,1}^1(\R^2)$, both divided differences $\dg_xf$ and $\dg_yf$ belong to the space of Schur multipliers $\fM_{\C,\C}$. This follows from the following analog of fundamental inequality \rf{fundrr}.

\begin{thm}
\label{fnner}
Let $f$ be a bounded function on $\R^2$ whose Fourier transform is supported in
$[-\s,\s]\times[-\s,\s]$. Then both $\dg_xf$ and $\dg_yf$ are Schur multipliers and
\bay
\label{anfunne}
\|\dg_xf\|_{\fM_{\C,\C}}+\|\dg_yf\|_{\fM_{\C,\C}}\le\const\s\|f\|_{L^\be}
\ey
\end{thm}

The proof of Theorem \ref{fnner} is based on the following lemma whose proved can be found in \cite{APPS}.

\begin{lem}
\label{razl}
Let $f$ be a bounded function on $\R$ whose Fourier transform is supported in $[-\s,\s]$. Then
$$
\frac{f(x)-f(y)}{x-y}=\sum_{n\in\Z}\s\frac{f(x)-f\big(\pi n\s^{-1}\big)}{\s x-\pi n}
\cdot\frac{\sin(\s y-\pi n)}{\s y-\pi n}
$$
Moreover,
\bay
\label{2vy}
\sum_{n\in\Z}\frac{\big|f(x)-f\big(\pi n\s^{-1}\big)\big|^2}{(\s x-\pi n)^2}
\le 3\|f\|_{L^\be(\R)}^2\quad\mbox{and}
\quad\sum_{n\in\Z}\frac{\sin^2 \s y}{(\s y-\pi n)^2}=1,~ x,\,y\in\R.
\ey
\end{lem}

\medskip

{\bf Proof of Theorem \ref{fnner}.} Clearly, it suffices to consider the case $\s=1$.
By Lemma \ref{razl}, we have
\bey
\big(\dg_xf\big)(z_1,z_2)=\frac{f(x_1,y_2)-f(x_2,y_2)}{x_1-x_2}=
\sum_{n\in\Z}\frac{f(\pi n,y_2)-f(x_2,y_2)}{\pi n-x_2}\cdot\frac{\sin(x_1-\pi n)}{x_1-\pi n}
\eey
and

\bey
\big(\dg_yf\big)(z_1,z_2)=\frac{f(x_1,y_1)-f(x_1,y_2)}{y_1-y_2}=\sum_{n\in\Z}\frac{f(x_1,y_1)-f(x_1,\pi n)}{y_1-\pi n}\cdot\frac{\sin(y_2-\pi n)}{y_2-\pi n}.
\eey

By Lemma \ref{razl}, we have
\bey
\sum_{n\in\Z}\frac{|f(x_1,y_1)-f(x_1,\pi n)|^2}{(y_1-\pi n)^2}\le3\|f(x_1,\cdot)\|_{L^\be(\R)}^2
\le3\|f\|_{L^\be(\C)}^2,\\
\sum_{n\in\Z}\frac{|f(\pi n,y_2)-f(x_2,y_2)|^2}{(\pi n-x_2)^2}\le3\|f(\cdot,y_2)\|_{L^\be(\R)}^2
\le3\|f\|_{L^\be(\C)}^2,\\
\eey
and
$$
\sum_{n\in\Z}\frac{\sin^2(x_1-n\pi)}{(x_1-\pi n)^2}=\sum_{n\in\Z}\frac{\sin^2(y_2-n\pi)}{(y_2-\pi n)^2}=1.
$$
It remains to observe that \rf{2vy} give us desired estimates of the norm of $\dg_xf$ and $\dg_yf$ in the Haagerup tensor product $L^\be\!\otimes_{\rm h}\!L^\be$. The result follows nor from Theorem \ref{tomSc}.
 $\bl$
 
Theorem \ref{fnner} implies the following result:

\begin{thm}
\label{nooL}
Let $f\in B_{\be,1}^1(\R^2)$. Then both $\dg_xf$ and $\dg_yf$ belong to $\fM_{\C,\C}$
and
$$
\|\dg_xf\|_{\fM_{\C,\C}}+\|\dg_yf\|_{\fM_{\C,\C}}\le\const\s\|f\|_{B_{\be,1}^1}.
$$
Moreover, if $N_1$ and $N_2$ are normal operators with bounded $N_1-N_2$,
then formula {\em\rf{irnoro}} holds and
$$
\|f(N_1)-f(N_2)\|\le\const\|f\|_{B_{\be,1}^1}\|N_1-N_2\|.
$$
\end{thm}

In other words, if $f\in B_{\be,1}^1(\R^2)$, then $f$ is an {\it operator Lipschitz function}.

To prove that $\dg_xf,~\dg_yf\in\fM_{\C,\C}$, it suffices to apply Theorem \ref{fnner} to each function $f*W_n$ (see \rf{Wn}). Formula \rf{irnoro} can be proved by analogy with the proof of formula \rf{fDK} for functions of self-adjoint operators. The operator Lipschitzness of $f$ follows immediately from formula \rf{irnoro}. We refer the reader to \cite{APPS} for details. 

As in the case of functions of self-adjoint operators, fundamental inequality
\rf{anfunne} allows us to establish for functions of perturbed normal operators analogs of all the results of Sections \ref{oLoD}--\ref{SchvNe} except for Theorem \ref{PoSu}.
In particular, a H\"older function of order $\a$, $0<\a<1$, on $\R^2$ must be operator H\"older of order $\a$. An analog of Theorem \ref{PoSu} was obtained in \cite{KPSS}.

We refer the reader to \cite{AP8} for more results on estimates of operator differences and quasicommutators for functions of normal operators.

\

\section{Functions of commuting self-adjoint operators}
\label{enki}

\

In the previous section we considered the behavior of functions of normal operators under
perturbation. This is equivalent to considering functions of pairs of commuting self-adjoint operators. Indeed, if $N$ is a normal operator, then $\re N$ and $\im N$ are commuting self-adjoint operators. On the other hand, if $A$ and $B$ are commuting self-adjoint operators, then $A+{\rm i}B$ is a normal operator.

In this section we are going to study functions of $d$-tuples of commuting self-adjoint operators. It is natural to try to use the approach for functions of normal operators that has been used in the previous section. However, it turns out that it does not work for $d\ge3$.

Indeed, a natural analog of formula \ref{irnoro} for functions of triple of commuting self-adjoint operators would be the following formula:

\begin{align}
\label{troika}
f(A_1,A_2,A_3)-f(B_1,B_2,B_3)=&\iint(\dg_1f)(x,y)\,dE_A(x)(A_1-B_1)\,dE_B(y)\nonumber\\[.2cm]
&+\iint (\dg_2f)(x,y)\,dE_A(x)(A_2-B_2)\,dE_B(y)\nonumber\\[.2cm]
&+\iint (\dg_3f)(x,y)\,dE_A(x)(A_3-B_3)\,dE_B(y),
\end{align}
where
$$
(\dg_1f)(x,y)=\frac{f(x_1,x_2,x_3)-f(y_1,x_2,x_3)}{x_1-y_1},\quad
(\dg_2f)(x,y)=\frac{f(y_1,x_2,x_3)-f(y_1,y_2,x_3)}{x_2-y_2},
$$
$$
(\dg_3f)(x,y)=\frac{f(y_1,y_2,x_3)-f(y_1,y_2,y_3)}{x_3-y_3},
\quad x=(x_1,x_2,x_3),\quad y=(y_1,y_2,y_3),
$$
and $E_A$ and $E_B$ are the joint spectral measures of the triples $(A_1,A_2,A_3)$ and 
$(B_1,B_2,B_3)$ on the Euclidean space $\R^3$.

It can easily be shown that \rf{troika} holds if
the functions $\dg_1f$, $\dg_2f$, and $\dg_3f$ belong to the space of Schur multipliers $\fM_{\R^3\!,\R^3}$ which would imply that $f$ is an operator Lipschitz function.

The methods of \cite{APPS} that were outlined in \S~\ref{fnorm}
allow us to prove that if $f$ is a bounded function on $\R^3$ with compactly supported
Fourier transform, then
$\dg_1f$ and $\dg_3f$ do 
belong to the space of Schur multipliers $\fM_{\R^3\!,\R^3}$. However, it turns out that the function $\dg_2f$ does not have to be in $\fM_{\R^3\!,\R^3}$, and so formula \rf{troika} cannot be used to prove that bounded functions on $\R^3$ with compactly supported Fourier transform must be operator Lipschitz. This was established in \cite{NP} where the following result was proved:

\begin{thm}
\label{kontr}
Suppose that $g$ is a bounded continuous function on $\R$ such that the Fourier transform of $g$ has compact support and is not a measure. Let $f$ be the function on $\R^3$ defined by
\bay
\label{ef}
f(x_1,x_2,x_3)=g(x_1-x_3)\sin x_2,\quad x_1,\,x_2,\,x_3\in\R.
\ey 
Then $f$ is a bounded function on $\R^3$ whose Fourier transform has compact support, but $\dg_2f\not\in\fM_{\R^3\!,\R^3}$.
\end{thm}

To construct  a function $g$ satisfying the hypothesis of Theorem \ref{kontr}, one can take, for example, the function $g$ defined by
$$
g(x)=\int_0^xt^{-1}\sin t\,dt,\quad x\in\R.
$$
Obviously, $g$ is bounded 
and
its Fourier transform $\F g$ satisfies the equality:
$$
(\F g)(t)=\frac ct
$$
for a nonzero constant $c$ and  sufficiently small positive $t$. It is easy to see that this implies that $\F g$ is not a measure.

Nevertheless, it was proved in \cite{NP} by a different method that functions in the Besov class $B_{\be,1}^1(\R^d)$ are {\it operator Lipschitz} in the sense that
$$
f(A_1,\cdots,A_d)-f(B_1,\cdots,B_d)\le\const\max_{1\le j\le d}\|A_j-B_j\|.
$$

The following result proved in \cite{NP} plays the same role as fundamental inequality
\rf{fundrr} in the case of functions of one self-adjoint operator.

\begin{lem}
\label{fnerdo}
Let $f$ be a bounded function on $\R^d$ whose Fourier transform is supported in 
$[-\s,\s]^d$. Then there are Schur multipliers $\Psi_j$, $1\le j\le d$, such that
$$
\|\Psi_j\|_{\fM_{\R^d,\R^d}}\le\const\s\|f\|_{L^\be}
$$
and
$$
f(x_1,\cdots,x_d)-f(y_1,\cdots,y_d)=\sum_{j=1}^d(x_j-y_j)\Psi_j(x_1,\cdots,x_d,y_1,\cdots,y_d).
$$
\end{lem}

Lemma \ref{fnerdo} implies the following result (see \cite{NP}) that considerably improves earlier estimates obtained in \cite{F2}.

\begin{thm}
\label{dshtuk}
Let $f$ be a function in the Besov class $B_{\be,1}^1(\R^d)$. Then there are Schur multipliers $\Phi_j$, $1\le j\le d$, such that
$$
\|\Phi_j\|_{\fM_{\R^d,\R^d}}\le\const\s\|f\|_{B_{\be,1}^1}
$$
and
$$
f(x_1,\cdots,x_d)-f(y_1,\cdots,y_d)=\sum_{j=1}^d(x_j-y_j)\Phi_j(x_1,\cdots,x_d,y_1,\cdots,y_d).
$$

If $(A_1,\cdots,A_d)$ and $(B_1,\cdots,B_d)$ are $d$-tuples of commuting self-adjoint operators, then
$$
f(A_1,\cdots,A_d)-f(B_1,\cdots,B_d)=\sum_{j=1}^d\iint\Phi_j(x,y)\,
dE_A(x)(A_j-B_j)\,dE_B(y),
$$
where $E_A$ and $E_B$ are the joint spectral measures of the $(A_1,\cdots,A_d)$ and 
$(B_1,\cdots,B_d)$.
\end{thm}

\begin{cor}
Let $f\in B_{\be,1}^1(\R^d)$. Then $f$ is operator Lipschitz.
\end{cor}

Lemma \ref{fnerdo} allows us to obtain analogs of all the results of Sections \ref{oLoD}--\ref{SchvNe} except for Theorem \ref{PoSu}. An analog of Theorem \ref{PoSu} for $d$-tuples of commuting self-adjoint operators was obtained in \cite{KPSS}.


\

\begin{center}
\bf\huge Chapter 2
\end{center}

\

\begin{center}
\bf\Large Multiple operator integrals with integrands in projective tensor products and their applications
\end{center}

\

\addtocontents{toc}{\vspace*{.3cm}\textbf{{\sc Chapter 2}. Multiple operator integrals with 
	integrands in projective \\ \hspace*{2.65cm} tensor products and their applications}
\hfill\pageref{iptpL}}
\renewcommand{\thesection}{2.\arabic{section}} 
\setcounter{section}{0}
\label{iptpL}

\

Multiple operator integrals were considered by several mathematicians, see \cite{Pa}, \cite{St}. However, those definitions required very strong restrictions on the classes of functions that can be integrated. In \cite{Pe4} multiple operator integrals were defined for functions that belong to the (integral) projective tensor product of $L^\be$ spaces. Later in \cite{JTT} multiple operator integrals were defined for Haagerup tensor products of $L^\be$ spaces.

In this chapter we consider applications of multiple operator integrals with integrands in the integral projective tensor product of $L^\be$ spaces. Such multiple operator integrals have nice Schatten--von Neumann properties. In Chapter 3 we shall see that multiple operator integrals with integrands in the Haagerup tensor product of $L^\be$ spaces do not possess such properties.

We consider in this chapter applications of multiple operator integrals to higher operator derivatives and estimates of higher operator differences. We also consider connections between multiple operator integrals and trace formulae for perturbations of class $\bS_m$, where $m$ is positive integer greater than $1$.

\

\section{A brief introduction to multiple operator integrals}
\label{vvekoi}

\

Multiple operator integrals are expressions of the form
$$
\underbrace{\int\cdots\int}_m\Psi(x_1,x_2,\cdots,x_m)
\,dE_1(x_1)T_1\,dE_2(x_2)T_2\cdots T_{m-1}\,dE_m(x_m).
$$
Here $E_1,\cdots,E_m$ are spectral measures on Hilbert space, $\Psi$ is a measurable function, and $T_1,\cdots,T_{m-1}$ are bounded linear operators on Hilbert space. The function $\Psi$ is called the {\it integrand} of the multiple operator integral.

For $m\ge3$, the Birman--Solomyak approach to double operator integrals does not work.
In \cite{Pe4} multiple operator integrals were defined for functions $\Psi$ that belong to the {\it integral projective tensor product} 
$L^\be(E_1)\hat\otimes_{\rm i}\cdots\hat\otimes_{\rm i}L^\be(E_m)$. It consists of functions $\Psi$ of the form
\bay
\label{itpm}
\Psi(x_1,\cdots,x_m)=\int_\O\f_1(x_1,\o)\f_2(x_2,\o)\cdots\f_m(x_m,\o)\,d\s(\o),
\ey
where $\f_1,\f_2,\cdots,\f_m$ are measurable functions such that
\bay
\label{iptpnm}
\int_\O\|\f_1(\cdot,\o)\|_{L^\be(E_1)}\|\f_2(\cdot,\o)\|_{L^\be(E_2)}\cdots
\|\f_m(\cdot,\o)\|_{L^\be(E_m)}\,d\s(\o)<\be.
\ey
If $\Psi$ belongs to $L^\be(E_1)\hat\otimes_{\rm i}\cdots\hat\otimes_{\rm i}L^\be(E_m)$, clearly, a representation of $\Psi$ in the form \rf{itpm} is not unique. The norm
$\|\Psi\|_{L^\be\hat\otimes_{\rm i}\cdots\hat\otimes_{\rm i}L^\be}$ is, by definition, the infimum of the expressions on the left-hand side of \rf{iptpnm} over all representations of $\Psi$ in the form of \rf{itpm}.

If $\Psi\in L^\be(E_1)\hat\otimes_{\rm i}\cdots\hat\otimes_{\rm i}L^\be(E_m)$ and $\Psi$ is represented as in \rf{itpm}, the multiple operator integral is defined by
\begin{align}
\label{oprkoi}
\underbrace{\int\cdots\int}_m&\Psi(x_1,\cdots,x_m)
\,dE_1(x_1)T_1\cdots T_{m-1}\,dE_m(x_m)\nonumber\\[.2cm]
&\df
\int_\O\left(\int\f_1(x_1,\o)\,dE_1(x_1)\right)T_1\cdots T_{m-1}
\left(\int\f_m(x_m,\o)\,dE_m(x_m)\right).
\end{align}

The following result shows that the multiple operator integral is well defined.

\begin{thm}
\label{nezav}
The expression on the right-hand side of {\em\rf{oprkoi}} does not depend on the choice of a representation of the form {\em\rf{itpm}}.
\end{thm}

The following proof is based on the approach of \cite{ACDS}.

\medskip

\Pf To simplify the notation, we assume that $n=3$. In the general case the proof is the same. Consider the right-hand side of \rf{oprkoi}. 
It is easy to see that it suffices to prove its independence on the choice of \rf{itpm}
for finite rank operators $T_1$ and $T_2$. It follows that we may assume that $\rank T_1=\rank T_2=1$. Let $T_1=(\cdot,u_1)v_1$ and $T_2=(\cdot,u_2)v_2$, where $u_1,v_1,u_2$ and $v_2$ are vectors in our Hilbert space. Suppose that $w_1$ and $w_2$ are arbitrary vectors. We are going to use the following notation:
\begin{align*}
W(&T_1,T_2)\\[.2cm]
&\df\!\int_\O\!\left(\int\f_1(x_1,\o)\,dE_1(x_1)\right)T_1\!
\left(\int\f_2(x_2,\o)\,dE_2(x_2)\right)T_2\!
\left(\int\f_3(x_3,\o)\,dE_3(x_3)\right).
\end{align*}
It is easy to verify that
\begin{align*}
\big(W&(T_1,T_2)w_1,w_2\big)\\[.2cm]
&=\int_\O\left(\int\f_1(x_1,\o)\,d\nu_1(x_1)\right)\left(\int\f_2(x_2,\o)\,d\nu_2(x_2)\right)
\left(\int\f_3(x_3,\o)\,d\nu_3(x_3)\right)\,d\s(\o)
\end{align*}
where
$$
\nu_1\df(E_1v_1,w_2),\quad\nu_2\df(E_2v_2,u_1)
\quad\mbox{and}\quad\nu_3\df(E_3w_1,u_2).
$$
Thus
\begin{align*}
\big(W(T_1,T_2)w_1,&w_2\big)\\[.2cm]
=&\iiint\left(\int_\O\f_1(x_1,\o)\f_2(x_2,\o)\f_3(x_3,\o)\,d\s(\o)\right)
\,d\nu_1(x_1)\,d\nu_2(x_2)\,d\nu_3(x_3)\\
=&\iiint\Psi(x_1,x_2,x_3)\,d\nu_1(x_1)\,d\nu_2(x_2)\,d\nu_3(x_3).
\end{align*}
It follows that $W(T_1,T_2)$ does not depend on the choice of a representation
of the form  \rf{itpm}. $\bl$

The following result is an easy consequence of the above definitions.

\begin{thm}
\label{opnor}
Let $\Psi$ be a function in  
$L^\be(E_1)\hat\otimes_{\rm i}\cdots\hat\otimes_{\rm i}L^\be(E_m)$. Suppose that $T_1,\cdots,T_{m-1}$ are bounded linear operator. Then
\begin{align*}
\left\|\underbrace{\int\cdots\int}_m\Psi(x_1,\cdots,x_m)
\,dE_1(x_1)T_1\cdots T_{m-1}\,dE_m(x_m)\right\|
\le\|\Psi\|_{L^\be\hat\otimes_{\rm i}\cdots\hat\otimes_{\rm i}L^\be}\prod_{j=1}^{m-1}\|T_j\|.
\end{align*}
\end{thm}

To simplify the notation, by $\bS_\be$ we mean the space of bounded linear operators on Hilbert space. The proof of the following result is also straightforward.

\begin{thm}
\label{Schpr}
Let $\Psi$ be a function in  
$L^\be(E_1)\hat\otimes_{\rm i}\cdots\hat\otimes_{\rm i}L^\be(E_m)$. Suppose that $p_j\ge1$, $1\le j\le m$, and $1/{p_1}+1/{p_2}+\cdots+1/{p_m}\le1$. If $T_1,T_2,\cdots,T_m$ are linear operators on Hilbert space such that $T_j\in\bS_{p_j}$, $1\le j\le m$, then
$$
\underbrace{\int\cdots\int}_m\Psi(x_1,\cdots,x_m)
\,dE_1(x_1)T_1\cdots T_{m-1}\,dE_m(x_m)\in\bS_r
$$
and
\begin{align*}
\left\|\underbrace{\int\cdots\int}_m\Psi(x_1,\cdots,x_m)
\,dE_1(x_1)T_1\cdots T_{m-1}\,dE_m(x_m)\right\|_{\bS_r}
\le\|\Psi\|_{L^\be\hat\otimes_{\rm i}\cdots\hat\otimes_{\rm i}L^\be}
\prod_{j=1}^{m-1}\|T_j\|_{\bS_{p_j}},
\end{align*}
where 
$$
1/r\df1-1/{p_1}-1/{p_2}-\cdots-1/{p_m}.
$$
\end{thm}

In particular, all the above facts hold for functions $\Psi$ in the projective tensor product $L^\be(E_1)\hat\otimes\cdots\hat\otimes L^\be(E_m)$ which consists of functions of the form
$$
\Psi(x_1,x_2,\cdots,x_m)=\sum_{k\ge1}\f^{[k]}_1(x_1)\f^{[k]}_2(x_2)\cdots\f^{[k]}_m(x_m),
$$
where $\f^{[k]}_j\in L^\be(E_j)$, $1\le j\le m$, and
$$
\sum_{k\ge1}
\big\|\f^{[k]}_1\big\|_{L^\be(E_1)}\big\|\f^{[k]}_2\big\|_{L^\be(E_2)}
\cdots\big\|\f^{[k]}_m\big\|_{L^\be(E_m)}<\be.
$$

\

\section{Higher operator derivatives}
\label{vyssh}

\

In \S~\ref{oLoD} we studied the problem of differentiability of the function
$t\mapsto f(A+tK)$, $t\in\R$, for self-adjoint operators $A$ and bounded self-adjoint operators $K$. In this section we are going to consider the problem of the existence of higher derivatives of this map. 

In the paper \cite{DK} Daletskii and Krein proved that in the case when the self-adjoint operator $A$ is bounded for nice functions $f$ the map $t\mapsto f(A+tK)$ has $m$-th derivative and it can be expressed in terms a multiple operator integral whose integrand is a higher order divided difference of $f$. 

Later in \cite{Pe4} the existence of higher operator differences was proved under much less restrictive assumptions.

\medskip

{\bf Definition.} For a $k$ times differentiable function $f$
the {\it divided differences $\dg^k f$ of order $k$} are defined inductively as follows:
$$
\dg^0f\df f;
$$
if $k\ge1$, then
\begin{align*}
(\dg^{k}f)&(s_1,\cdots,s_{k+1})\\[.2cm]
\df&
\left\{\begin{array}{ll}\displaystyle
{\frac{(\dg^{k-1}f)(s_1,\cdots,s_{k-1},s_k)-
(\dg^{k-1}f)(s_1,\cdots,s_{k-1},s_{k+1})}{s_{k}-s_{k+1}}},&s_k\ne s_{k+1},\\[.4cm]
\displaystyle{\frac{\partial}{\partial t}}
\Big(\big(\dg^{k-1}f\big)(s_1,\cdots,s_{k-1},t)\Big)\Big|_{t=s_k},&s_k=s_{k+1},
\end{array}\right.
\end{align*}
(the definition does not depend on the order of the variables). Note that
$\dg\f=\dg^1\f$.

The following result was obtained in \cite{Pe4}. To state it, we denote by
$\fB(\R)$ the space of bounded Borel functions on $\R$ endowed with the norm
$$
\|\f\|_{\fB(\R)}\df\sup_{t\in\R}|\f(t)|.
$$

\begin{thm}
\label{razravpo}
Let $m$ be a positive integer and let $f$ be a function in the Besov class 
$B_{\be,1}^m(\R)$. Then 
$\dg_m f\in\underbrace{\fB(\R)\hat\otimes_{\rm i}\cdots\hat\otimes_{\rm i}\fB(\R)}_{m+1}$
and
$$
\big\|\dg_m f\|_{\fB(\R)\hat\otimes_{\rm i}\cdots\hat\otimes_{\rm i}\fB(\R)}
\le\const\|f\|_{B_{\be,1}^m}.
$$
\end{thm}

Note that the integral projective tensor product of copies of $\fB(\R)$ can be defined in the same way as the integral projective tensor product of $L^\be$ spaces.

We sketch the proof of Theorem \ref{razravpo} in the special case $m=2$. In the general case the proof is the same. Let $f$ be a bounded functions on $\R$ whose Fouriesr transform is a compact subset of $[0,\be)$.
The following identity is an analog of formula \rf{ifrr}:
\begin{align}
\label{tri}
(\dg^2f)(s_1,s_2,s_3)=&
-\iint\limits_{\R_+\times\R_+}\big(f*\mu_{u+v}\big)(s_1)e^{-{\rm i}(u+v)s_1}
e^{{\rm i}vs_2}e^{{\rm i}us_3}
\,du\,dv\nonumber\\[.2cm]
&-\iint\limits_{\R_+\times\R_+}\big(f*\mu_{u+v}\big)(s_2)e^{-{\rm i}(u+v)s_2}
e^{{\rm i}us_1}e^{{\rm i}vs_3}
\,du\,dv\nonumber\\[.2cm]
&-\iint\limits_{\R_+\times\R_+}\big(f*\mu_{u+v}\big)(s_3)e^{-{\rm i}(u+v)s_3}
e^{{\rm i}us_1}e^{{\rm i}vs_2}
\,du\,dv.
\end{align}
This is a simplified version of formula (5.6) in \cite{Pe4}.

As in \S~\ref{oDf}, it is easy to deduce from \rf{tri} the following
estimate
$$
\|\dg^2f\|_{\fB(\R)\hat\otimes_{\rm i}\fB(\R)\hat\otimes_{\rm i}\fB(\R)}\le\const\s^2\|f\|_{L^\be},
$$
whenever $f$ is a bounded function on $\R$ whose Fourier transform is supported 
in $[0,\s]$.

The following theorem about the existence of the $m$th derivative of the function 
$t\mapsto f(A_t)$, where $A_t\df A+tK$, was obtained in\cite{Pe4}.

\begin{thm}
\label{dps}
Let $m$ be a positive integer.
Suppose that $A$ is a self-adjoint operator and $K$ is a bounded self-adjoint operator.
If $f\in B_{\be1}^m(\R)\bigcap B_{\be1}^1(\R)$, then the function
$t\mapsto f(A_t)$
has $m$th derivative that is a bounded operator and
\bay
\label{myapro}
\!\!\!\!\!\frac{d^m}{dt^m}\big(f(A_t)\big)\Big|_{t=0}=
m!\underbrace{\int\cdots\int}_{m+1}(\dg^{m}f)(s_1,\cdots,s_{m+1})
\,dE_A(s_1)K\cdots K\,dE_A(s_{m+1}).
\ey
\end{thm} 

We refer the reader to \cite{Pe4} for the proof.

\medskip

{\bf Remark.} Suppose that $f\in B_{\be1}^m(\R)$, $m\ge2$, but $f$ does not necessarily belong to 
$B_{\be1}^1(\R)$. In this case we still can define the $m$th derivative of the function 
$t\mapsto f(A_t)$ in the following way. We put
\bay
\label{vryad}
\frac{d^m}{dt^m}\big(f(A_t)\big)\Big|_{t=0}=
\sum_{n\in\Z}\frac{d^m}{dt^m}\big(f_n(A_t)\big)\Big|_{t=0},
\ey
where $f_n=f*W_n$, see \rf{fn}. Then the series on the right-hand side of \rf{vryad}
converges absolutely in the norm. With this (natural) definition it can easily happen that the function $t\mapsto f(A_t)$ can have $m$th derivative, but not necessarily the first derivative. We refer the reader to \cite{Pe4} for details.

\medskip

In a similar was one can consider the problem of taking higher operator derivatives for functions of unitary operators. Note that in \cite{Pe4} the formula for the $m$th derivative of the function $t\mapsto f(e^{{\rm i}tA})U$ has an error. A correct formula is an easy consequence of the results of Section 5 of \cite{AP1}.

Note also that similar results and similar formulae can be obtained for functions of contractions and for functions of dissipative operators, see \cite{AP1} and \cite{AP6}.

\

\section{Higher operator differences}
\label{vysshra}

\

In Chapter 1 we have seen that formula \rf{fDK} plays a significant role in estimating various norms of the operator differences $f(A)-f(B)$. In this section we are going to study higher order operator differences
$$
\big(\D_K^mf\big)(A)\df
\sum_{j=0}^m(-1)^{m-j}\left(\begin{matrix}m\\j\end{matrix}\right)f\big(A+jK\big),
$$
where $A$ and $K$ are self-adjoint operators on Hilbert space. We consider here only bounded self-adjoint operators $A$ and $K$ and refer the reader to \cite{AP4} for a detailed study of the case when $A$ is an unbounded self-adjoint operators.

As in the case of operator differences, an essential role is played by integral formulae for higher operator differences. In \cite{AP2} it was shown that higher operator differences can be represented in terms of multiple operator integrals. This allowed one to obtain analogs of the results discussed in Chapter 1 for higher operator differences.

Recall that for functions $f$ in the Besov class $B_{\be,1}^m(\R)$, the divided difference 
$\dg^m f$ of order $m$ belongs to the integral projective  product
$\underbrace{\fB(\R)\hat\otimes_{\rm i}\cdots\hat\otimes_{\rm i}\fB(\R)}_{m+1}$. 
The following formula was obtained in \cite{AP2}.

\begin{thm}
\label{ifdorvp}
Let $f\in B_{\be,1}^m(\R)$ and let $A$ and $K$ be bounded self-adjoint operators on Hilbert space. Then
\begin{align*}
\big(\D_K^mf\big)&(A)\\[.2cm]
=&
m!\underbrace{\int\cdots\int}_{m+1}(\dg^{m}f)(s_1,\cdots,s_{m+1})
\,dE_A(s_1)K\,dE_{A+K}(s_2)K\cdots K\,dE_{A+mK}(s_{m+1}).
\end{align*}
\end{thm}

Let us prove Theorem \ref{ifdorvp} in the special case $m=2$.

\medskip

\Pf Let $f\in B_{\be,1}^2(\R)$. We should prove the following formula:

$$
f(A+K)-2f(A)+f(A-K)=\!2\!\iiint\!(\dg^2f)(s,t,u)\,dE_{A+K}(s)K\,dE_A(t)K\,dE_{A-K}(u).
$$
Put $T=f(A+K)-2f(A)+f(A-K)$. By \rf{fDK},

\begin{align*}
T&=
f(A+K)-f(A)-\big(f(A)-f(A-K)\big)\\[.2cm]
&=\iint(\dg f)(s,t)\,dE_{A+K}(s)K\,dE_A(t)-
\iint(\dg f)(s,t)\,dE_{A}(s)K\,dE_{A-K}(t)\\[.2cm]
&=\iint(\dg f)(s,t)\,dE_{A+K}(s)K\,dE_A(t)-
\iint(\dg f)(s,u)\,dE_{A+K}(s)K\,dE_{A-K}(u)\\[.2cm]
&+\iint(\dg f)(s,t)\,dE_{A+K}(s)K\,dE_{A-K}(t)-
\iint(\dg f)(s,t)\,dE_{A}(s)K\,dE_{A-K}(t).
\end{align*}

We have
\begin{align*}
\iint(\dg f)(s,t)\,dE_{A+K}(s)K\,dE_A(t)&-
\iint(\dg f)(s,u)\,dE_{A+K}(s)K\,dE_{A-K}(u)\\[.2cm]
=&\iiint(\dg f)(s,t)\,dE_{A+K}(s)K\,dE_A(t)\,dE_{A-K}(u)\\[.2cm]
-&\iiint(\dg f)(s,u)\,dE_{A+K}(s)K\,dE_A(t)\,dE_{A-K}(u)\\[.2cm]
=&\!\iiint\!(t-u)(\dg^2f)(s,t,u)dE_{A+K}(s)KdE_A(t)dE_{A-K}(u)\\[.2cm]
=&\iiint(\dg^2f)(s,t,u)\,dE_{A+K}(s)K\,dE_A(t)K\,dE_{A-K}(u).
\end{align*}

Similarly,
\begin{align*}
\iint(\dg f)(s,t)\,dE_{A+K}(s)K\,dE_{A-K}(t)&-
\iint(\dg f)(s,t)\,dE_{A}(s)K\,dE_{A-K}(t)\\[.2cm]
=&\iiint(\dg^2f)(s,t,u)dE_{A+K}(s)KdE_A(t)KdE_{A-K}(u).
\end{align*}

Thus
$$
T=2\iiint(\dg^2f)(s,t,u)\,dE_{A+K}(s)K\,dE_A(t)K\,dE_{A-K}(u).\quad\bl
$$

The proof of the following result obtained in \cite{AP2} is similar to the proof of Theorem \ref{Hoa}.

\begin{thm}
\label{sam}
Let $0<\a<m$ and let $f\in\L_\a(\R)$. Then there exists a constant $c>0$ such that for every self-adjoint operators  $A$ and $K$ on Hilbert space the following inequality holds:
$$
\left\|\big(\D_K^mf\big)(A)\right\|
\le c\,\|f\|_{\L_\a(\R)}\cdot\|K\|^\a.
$$
\end{thm}

In particular, in the case $\a=1$, Theorem \ref{sam} means the following: let $f$ be a function in the Zygmund class $\L_1(\R)$, i.e., $f$ is a continuos function on $\R$ such that
$$
|f(s-t)-2f(s)+f(s+t)|\le\const|t|,
$$
then
$$
\|f(A-K)-2f(A)+f(A+K)\|\le\const\|f\|_{\L_1}\|K\|.
$$

We refer the reader to \cite{AP2} for an analog of Theorem \ref{modne} for higher order moduli of continuity.

To conclude this section, we also mention that the results of \S~\ref{SchvNe} were generalized in \cite{AP3} to the case higher order operator differences. We state here the following result whose proved can be found in \cite{AP3}.

\begin{thm}
\label{posl}
Let $\a>0$, $m-1\le\a<m$, and $m<p<\be$. There exists a positive number $c$ such that for every $f\in\L_\a(\R)$, for an arbitrary self-adjoint operator $A$, and an arbitrary self-adjoint operator $K$ of class $\bS_p$, the following inequality holds:
$$
\left\|\big(\D_K^mf\big)(A)\right\|_{\bS_{p/\a}}
\le c\,\|f\|_{\L_\a(\R)}\|K\|_{\bS_p}^\a.
$$
\end{thm}

Note also that similar results hold for functions of unitary operators, functions of contractions and functions of dissipative operators, see \cite{AP2}, \cite{AP3} and 
\cite{AP6}.

\

\section{Trace formulae for perturbations of class ${\bS_m}$, ${m\ge2}$}
\label{vozmSm}

\

In \S~\ref{fsLK} we have considered the Lifshits--Krein trace formula for $f(A)-f(B)$ in the case when $B$ is a trace class perturbation of $A$. In \cite{Ko} Koplienko considered the case of Hilbert--Schmidt perturbations and he found a trace formula for the second order Taylor approximation
\bay
\label{fsK}
\trace\left(f(A+K)-f(A)-\frac{d}{dt}f(A+tK)\Big|_{t=0}\right)
=\int_\R f''(x)\eta(x)\,dx.
\ey
Here $A$ is a self-adjoint operator, $K$ is a self-adjoint operator of class $\bS_2$ and $\eta$ is a function in $L^1$ that is determined by $A$ and $K$. It is called the {\it spectral shift function of order} 2. In \cite{Ko} formula \rf{fsK} was proved for rational functions with poles off $\R$.

Formula \rf{fsK} was generalized in \cite{Pe9} to the case when $f$ is an arbitrary function in the Besov class $B_{\be,1}^2(\R)$.

In \cite{PSS} the authors considered the more general problem of perturbation of class $\bS_m$, where $m$ is an arbitrary positive integer and they obtained the following trace formula for the Taylor approximation ${\mathscr T}^{(m)}_{A,K}f$ of order $m$:
\begin{align*}
{\mathscr T}^{(m)}_{A,K}f&\df f(A+K)-f(A)\\[.2cm]
&-\frac{d}{dt}f(A+tK)\Big|_{t=0}-
\cdots-\frac1{(m-1)!}\frac{d^{m-1}}{dt^{m-1}}f(A+tK)\Big|_{t=0}.
\end{align*}
They proved in \cite{PSS} that there is a unique function $\eta_m$ in $L^1$ that depends only on $A$, $K$ and $m$ such that
\bay
\label{fsdaT}
\trace\left({\mathscr T}^{(m)}_{A,K}f\right)=\int_\R f^{(m)}(s)\eta_m(s)\,ds
\ey
for functions $f$ on $\R$ such that the derivatives $f^{(j)}$ are Fourier transforms of $L^1$ functions for $1\le j\le m$. The function $\eta_m$ is called the {\it spectral shift
function of order} $m$.

The results of \cite{PSS} were improved in \cite{AP7}. First, formula \rf{fsdaT} was extended for arbitrary functions $f$ in the Besov class $B_{\be,1}^m(\R)$. Secondly, much more general trace formulae for perturbations of class $\bS_m$ we obtained in \cite{AP7}.

It was shown in \cite{AP7} that the Taylor approximation admits the following representation in terms of the multiple operator integral:
$$
{\mathscr T}^{(m)}_{A,K}f=\underbrace{\int\!\!\cdots\!\!\int}_{m+1}(\dg^{m}f)(s_1,\cdots,s_{m+1})
\,dE_{A+K}(s_1)K\,dE_A(s_2)K\cdots K\,dE_A(s_{m+1}).
$$
Here $A$ is a self-adjoint operator, $K$ is a self-adjoint operator of class $\bS_1$ and $f\in B_{\be,1}^m(\R)$. In this formula by ${\mathscr T}^{(m)}_{A,K}f$ we mean
$$
{\mathscr T}^{(m)}_{A,K}f\df\sum_{n\in\Z}{\mathscr T}^{(m)}_{A,K}f_n,
$$
where as usual $f_n=f*W_n$, see \rf{fn}.

To establish formula \rf{fsdaT}, the authors of \cite{PSS} proved the following inequality:
\bay
\label{nelnpr}
\left|\tr\left(\frac{d^{m}}{dt^{m}}f(A_t)\Big|_{t=0}\right)\right|\le
\const\big\|f^{(m)}\big\|_{L^\be}\|K\|^m_{\bS_m},
\ey
for functions $f$ whose derivatives $f^{(j)}$, $1\le j\le m$, are Fourier transforms of $L^1$ functions. Here $A_t\df A+tK$.

To prove that formula \rf{fsdaT} holds for arbitrary functions $f$ in $B_{\be,1}^m(\R)$,
we have to extend inequality \rf{nelnpr} to the class $B_{\be,1}^m(\R)$. Recall that for 
$f\in B_{\be,1}^m(\R)$, by the $m$th derivative of the functions $t\mapsto f(A_t)$, we mean
$$
\sum_{n\in\Z}\frac{d^{m}}{dt^{m}}f_n(A_t).
$$

\begin{thm}
\label{leompro}
Let $f\in B_{\be1}^m(\R)$ and $K\in\bS_m$. Then
\bay
\label{PsS}
\left\|\tr\left(\frac{d^{m}}{dt^{m}}f(A_t)\right)\right\|_{L^\infty}\le
\const\big\|f^{(m)}\big\|_{L^\be}\|K\|^m_{\bS_m}.
\ey
\end{thm}

Note that the proof of Theorem \ref{leompro} given in \cite{AP7} contains an inaccuracy. We give here a corrected proof.

\Pf As before, it suffices to consider the case when $f$ is a bounded function on $\R$ whose Fourier transform has compact support in $(0,\be)$.
By Theorem \ref{dps}, we have
$$
\frac{d^m}{dt^m}\big(f(A_t)\big)\Big|_{t=0}=
m!\underbrace{\int\cdots\int}_{m+1}(\dg^{m}f)(s_1,\cdots,s_{m+1})
\,dE_A(s_1)K\cdots K\,dE_A(s_{m+1}).
$$

For simplicity we assume that $m=2$. In the general case the proof is the same.
Recall formula \rf{tri}:
\begin{align*}
(\dg^2f)(s_1,s_2,s_3)=&
-\iint\limits_{\R_+\times\R_+}\big(f*\mu_{u+v}\big)(s_1)e^{-{\rm i}(u+v)s_1}
e^{{\rm i}vs_2}e^{{\rm i}us_3}
\,du\,dv\\[.2cm]
&-\iint\limits_{\R_+\times\R_+}\big(f*\mu_{u+v}\big)(s_2)e^{-{\rm i}(u+v)s_2}
e^{{\rm i}us_1}e^{{\rm i}vs_3}
\,du\,dv\\[.2cm]
&-\iint\limits_{\R_+\times\R_+}\big(f*\mu_{u+v}\big)(s_3)e^{-{\rm i}(u+v)s_3}
e^{{\rm i}us_1}e^{{\rm i}vs_2}
\,du\,dv.
\end{align*}

Thus  
\begin{align}
\label{triint}
\frac{d^2}{dt^2}\big(f(A_t)\big)\Big|_{t=0}
=&-2\iint\limits_{\R_+\times\R_+}\big(f*\mu_{u+v}\big)(A)e^{-{\rm i}(u+v)A}K
e^{{\rm i}vA}Ke^{{\rm i}uA}\,du\,dv\nonumber\\[.2cm]
&-2\iint\limits_{\R_+\times\R_+}e^{{\rm i}uA}K
\big(f*\mu_{u+v}\big)(A)e^{-{\rm i}(u+v)A}Ke^{{\rm i}vA}\,du\,dv\nonumber\\[.2cm]
&-2\iint\limits_{\R_+\times\R_+}e^{{\rm i}uA}Ke^{{\rm i}vA}K
\big(f*\mu_{u+v}\big)(A)e^{-{\rm i}(u+v)A}\,du\,dv.
\end{align}

Let $\o$ be a function in $C^\be(\R)$ such that $\o(0)=1$ and $\F\o$ is a nonnegative infinitely differentiable function
with compact support. For $\e>0$, we put $f_\e(x)\df\o(\e x)f(x)$.
Then $\supp\F f_\e$ is a compact and
$$
\F f_\e\in L^1(\R)\cap C^\be(\R).
$$
Hence, by
Theorem 2.1 of \cite{PSS},
$$
\left|\tr\left(\frac{d^{2}}{dt^{2}}f_\e(A_t)\Big|_{t=0}\right)\right|\le
\const\big\|f_\e^{''}\big\|_{L^\be}\|K\|^2_{\bS_2}.
$$
It is easy to see that
$$
\big\|f_\e^{''}\big\|_{L^\be}\le C\big\|f^{''}\big\|_{L^\be},
$$
where $C$ depends only on $\Phi$.
Moreover, $\supp\F f_\e$ is a compact subset of $(0,\be)$ for sufficiently small $\e$.
It follows that 
$$
\|f_\e*\mu_{u+v}\|_{L^\be}\le\const\|f*\mu_{u+v}\|_{L^\be}\quad\mbox{and}\quad 
\lim_{\e\to0}(f_\e*\mu_{u+v})(s)=(f*\mu_{u+v})(s),\quad s\in\R.
$$
By the spectral theorem,
$$
\lim_{\e\to0}(f_\e*\mu_{u+v})(A)=(f*\mu_{u+v})(A)
$$
in the strong operator topology. Thus
$$
\lim_{\e\to0}\big(f_\e*\mu_{u+v}\big)(A)e^{-{\rm i}(u+v)A}Ke^{{\rm i}vA}Ke^{{\rm i}uA}
=\big(f*\mu_{u+v}\big)(A)e^{-{\rm i}(u+v)A}Ke^{{\rm i}vA}Ke^{{\rm i}uA}
$$
in the norm of $\bS_1$. It follows that
\begin{align*}
\lim_{\e\to0}
\iint\limits_{\R_+\times\R_+}\big(f_\e*\mu_{u+v}\big)(A)&e^{-{\rm i}(u+v)A}K
e^{{\rm i}vA}Ke^{{\rm i}uA}\,du\,dv\\[.2cm]
&=\iint\limits_{\R_+\times\R_+}\big(f*\mu_{u+v}\big)(A)e^{-{\rm i}(u+v)A}K
e^{{\rm i}vA}Ke^{{\rm i}uA}\,du\,dv
\end{align*}
in the norm of $\bS_1$. The same is true about the second and the third integral on the right-hand side of \rf{triint}. This proves that 
$$
\lim_{\e\to0}\frac{d^2}{dt^2}\big(f_\e(A_t)\big)\Big|_{t=0}
=\frac{d^2}{dt^2}\big(f(A_t)\big)\Big|_{t=0}
$$
in the norm of $\bS_1$, and so \rf{PsS} holds with $m=2$. $\bl$

Note that Theorem \ref{leompro} was used in \cite{AP7} to obtain considerably more general trace formulae. In particular trace formulae were found for 
$$
\trace \frac{d^{m}}{dt^{m}}\big(f(A_t)\big)\Big|_{t=0}\quad\mbox{and}\quad
\trace \big(\D_K^mf\big)(A).
$$


\

\begin{center}
\bf\huge Chapter 3
\end{center}

\

\begin{center}
\bf\Large Triple operator integrals, Haagerup(-like) tensor products and functions of noncommuting operators
\end{center}

\

\addtocontents{toc}{\vspace*{.3cm}\textbf{{\sc Chapter 3}. Triple operator integrals, 
Haagerup(-like) tensor products \\ \hspace*{2.65cm}
and functions of noncommuting operators}\hfill\pageref{Stoi}}
\renewcommand{\thesection}{3.\arabic{section}} 
\setcounter{section}{0}
\label{Stoi}

\

In this chapter we deal with triple operator integrals and we apply triple operator integrals to estimates of functions of perturbed noncommuting pairs of self-adjoint operators. It turns out that for this purpose it is not enough to consider triple operator integral whose integrands belong to the (integral) projective tensor product of $L^\be$ spaces. In \cite{JTT} multiple operator integrals were defined for functions that belong
to the Haagerup tensor product of $L^\be$ spaces. We define triple operator integrals for functions in the Haagerup tensor product in \S~\ref{Haagtpatroi}. However, for our purpose we have to modify the notion of the Haagerup tensor product. We define in \S~\ref{Ttoi} Haagerup-like tensor products of the first kind and of the second kind. 
We are going to use the following representation of $f(A_1,B_1)-f(A_2,B_2)$ in terms of triple operator integrals:
\begin{align}
\label{predst}
\!\!\!\!f(A_1,\!B_1)\!-\!f(A_2,\!B_2)&\!=\!
\iiint\big(\dg^{[1]}f\big)(x_1,x_2,y)\,dE_{A_1}(x_1)(A_1-A_2)\,dE_{A_2}(x_2)\,dE_{B_1}(y)
\nonumber\\[.2cm]
+\!\!&\iiint\!\!\big(\dg^{[2]}f\big)(x,y_1,y_2)dE_{A_2}(x)dE_{B_1}(y_1)(B_1-B_2)dE_{B_2}(y_2),
\end{align}
where the divided differences $\dg^{[1]}f$ and $\dg^{[2]}f$ are defined by
$$
\big(\dg^{[1]}f\big)(x_1,x_2,y)\df\frac{f(x_1,y)-f(x_2,y)}{x_1-x_2}
\quad\mbox{and}\quad
\big(\dg^{[2]}f(x,y_1,y_2)\big)\df\frac{f(x,y_1)-f(x,y_2)}{y_1-y_2}.
$$
Here $f$ is a function in the Besov class $B_{\be,1}^1(\R^2)$ and $(A_1,B_1)$ and $(A_2,B_2)$ are pairs of (not necessarily commuting) self-adjoint operators. 

It turns out that the divided differences do not have to belong to the integral projective tensor product of the $L^\be$ spaces. That is why we have to consider triple operator integrals defined for other classes of functions. In \S~\ref{Haagtpatroi} we define 
the Haagerup tensor product of $L^\be$ spaces and triple operator integrals for such functions. It turned out, however, that the divided differences do not have to belong to the Haagerup tensor product. To overcome the problems, we introduce in \S~\ref{Ttoi} Haagerup-like tensor products of the first kind and of the second kind. We will see in 
\S~\ref{ddiff} that for functions $f$ in $B_{\be,1}^1(\R^2)$, the divided difference 
$\dg^{[1]}f$ belongs to the Haagerup-like tensor product of the first kind, while the divided difference $\dg^{[2]}f$ belongs to the Haagerup-like tensor product of the second kind.

We obtain in \S~\ref{ple2} Lipschitz type estimates for functions of noncommuting self-adjoint operators in the norm of $\bS_p$ with $p\in[1,2]$. It turns out that such Lipschitz type estimates in the norm of $\bS_p$ for $p>2$ and in the operator norm do not hold.

Finally, we use in \S~\ref{2c} triple operator integrals with integrands in Haagerup-like tensor products to estimates trace norms of commutators of functions of almost commuting operators. 

In the first section of this chapter we define functions of noncommuting self-adjoint operators.

\medskip

\section{Functions of noncommuting self-adjoint operators}
\label{fisdlynekomop}

\medskip

Let $A$ and $B$ be self-adjoint operators on Hilbert space and let $E_A$ and $E_B$ be their spectral measures. Suppose that $f$ is a function of two variables that is defined at least on $\s(A)\times\s(B)$, where $\s(A)$ and $\s(B)$ are the spectra of $A$ and $B$. If $f$ is a Schur multiplier with respect to the pair $(E_A,E_B)$, we define the function $f(A,B)$ of $A$ and $B$ by
\bay
\label{fAB}
f(A,B)\df\iint f(x,y)\,dE_A(x)\,dE_B(y).
\ey
Note that the map $f\mapsto f(A,B)$ is linear, but not multiplicative.

If we consider functions of bounded operators, without loss of generality we may deal with periodic functions with a sufficiently large period. Clearly, we can rescale the problem and assume that our functions are $2\pi$-periodic in each variable.

If $f$ is a trigonometric polynomial of degree $N$, we can represent $f$ in the form
$$
f(x,y)=\sum_{j=-N}^Ne^{{\rm i}jx}\left(\sum_{k=-N}^N\hat f(j,k)e^{{\rm i}ky}\right).
$$
Thus 
$f$ belongs to the projective tensor product $L^\be\hat\otimes L^\be$ and
$$
\|f\|_{L^\be\hat\otimes L^\be}\le\sum_{j=-N}^N\sup_y
\left|\sum_{k=-N}^N\hat f(j,k)e^{{\rm i}ky}\right|
\le(1+2N)\|f\|_{L^\be}
$$
It follows easily from \rf{Bperf} that every periodic function $f$ on $\R^2$ of Besov class $B_{\be1}^1$ of periodic functions belongs to 
$L^\be\hat\otimes L^\be$, and so the operator $f(A,B)$ is well defined by \rf{fAB}.

Note that the above definitions of functions of noncommuting operators is related to the Maslov theory, see \cite{Ma}. If $A$ and $B$ are self-adjoint operators, we can consider the transformer $\cL_A$ of left multiplication by $A$ and the transformer $\cR_B$ of right multiplication by $B$:
$$
\cL_AT\df AT,\quad \cR_BT\df TB.
$$
Clearly, the transformers $\cL_A$ and $\cR_B$ commute.

We can consider the transformers $\cL_A$ and $\cR_B$ defined on the Hilbert Schmidt class $\bS_2$. In this case they are commuting self-adjoint operators on $\bS_2$ and the spectral theorem allows us to define functions $f(\cL_A,\cR_B)$ for all bounded Borel functions $f$ on $\R^2$.

If our Hilbert space is finite-dimensional, the definition of $f(A,B)$ given by \rf{fAB} is equivalent to the following one:
$$
f(A,B)=f(\cL_A,\cR_B)I,
$$
where $I$ is the identity operator, and so the definition of functions of noncommuting operators can be reduced to the functional calculus for the commuting self-adjoint operators on the Hilbert Schmidt class.

If our Hilbert space $\h$ is infinite-dimensional, we cannot apply $f(\cL_A,\cR_B)$ to the identity operator, which does not belong to the Hilbert Schmidt class. In this case we can consider the transformers $\cL_A$ and $\cR_B$ as commuting bounded linear operators on the space $\mB(\h)$ of bounded linear operators on $\h$. However, since $\mB(\h)$ is not a Hilbert space and we cannot use the spectral theorem to define functions of $\cL_A$ and $\cR_B$. Nevertheless, if $f$ is a sufficiently nice function, we can define 
$f(\cL_A,\cR_B)$, in which case the functions $f(A,B)$ defined by \rf{fAB} coincide with
$f(\cL_A,\cR_B)I$.

\

\section{Haagerup tensor products and triple operator integrals}
\label{Haagtpatroi}

\

We proceed now to the approach to multiple operator integrals based on the Haagerup tensor product of $L^\be$ spaces. We refer the reader to the book \cite{Pi} for detailed \lb information about Haagerup tensor products of operator spaces.
The {\it Haagerup tensor product} 
$L^\be(E_1)\!\otimes_{\rm h}\!L^\be(E_2)\!\otimes_{\rm h}\!L^\be(E_3)$ of $L^\be$ spaces
is defined as the space of functions $\Psi$ of the form
\bay
\label{htr}
\Psi(x_1,x_2,x_3)=\sum_{j,k\ge0}\a_j(x_1)\b_{jk}(x_2)\g_k(x_3),
\ey
where $\a_j$, $\b_{jk}$, and $\g_k$ are measurable functions such that
\bay
\label{ogr}
\{\a_j\}_{j\ge0}\in L_{E_1}^\be(\ell^2), \quad 
\{\b_{jk}\}_{j,k\ge0}\in L_{E_2}^\be({\mathcal B}),\quad\mbox{and}\quad
\{\g_k\}_{k\ge0}\in L_{E_3}^\be(\ell^2),
\ey
where ${\mathcal B}$ is the space of matrices that induce bounded linear operators on $\ell^2$ and this space is equipped with the operator norm. In other words,
$$
\|\{\a_j\}_{j\ge0}\|_{L^\be(\ell^2)}\df
E_1\mbox{-}\ess\sup\left(\sum_{j\ge0}|\a_j(x_1)|^2\right)^{1/2}<\be,
$$
$$
\|\{\b_{jk}\}_{j,k\ge0}\|_{L^\be({\mathcal B})}\df
E_2\mbox{-}\ess\sup\|\{\b_{jk}(x_2)\}_{j,k\ge0}\|_{{\mathcal B}}<\be,
$$
and
$$
\|\{\g_k\}_{k\ge0}\|_{L^\be(\ell^2)}\df
E_3\mbox{-}\ess\sup\left(\sum_{k\ge0}|\g_k(x_3)|^2\right)^{1/2}<\be.
$$
The norm of $\Psi$ in 
$L^\be\!\otimes_{\rm h}\!L^\be\!\otimes_{\rm h}\!L^\be$ is, by definition, 
the infimum of
$$
\|\{\a_j\}_{j\ge0}\|_{L^\be(\ell^2)}\|\{\b_{jk}\}_{j,k\ge0}\|_{L^\be({\mathcal B})}
\|\{\g_k\}_{k\ge0}\|_{L^\be(\ell^2)}
$$
over all representations of $\Psi$ of the form \rf{htr}.

It is well known that $L^\be\hat\otimes L^\be\hat\otimes L^\be\subset
L^\be\!\otimes_{\rm h}\!L^\be\!\otimes_{\rm h}\!L^\be$. Indeed,
suppose that $\Psi$ admits a representation
$$
\Psi(x_1,x_2,x_3)=\sum_n\f_n(x_1)\psi_n(x_2)\chi_n(x_3)
$$
with
$$
\sum_n\|\f_n\|_{L^\be(E_1)}\|\psi_n\|_{L^\be(E_2)}\|\chi_n\|_{L^\be(E_3)}<\be.
$$ 
Without loss of generality we may assume 
$$
c_n\df\|\f_n\|_{L^\be}\|\psi_n\|_{L^\be}\|\chi_n\|_{L^\be}\ne0
\quad\mbox{for every}\quad n.
$$
We define $\a_j$, $\b_{j,k}$ and $\g_k$ by
$$
\a_j(x_1)=\frac{\sqrt{c_j}}{\|\f_j\|_{L^\be}}\f_j(x_1),\quad
\g_k(x_3)=\frac{\sqrt{c_k}}{\|\chi_k\|_{L^\be}}\chi_j(x_3)
$$
and
$$
\b_{jk}(x_2)=\left\{\begin{array}{ll}\psi_j(x_2)\|\psi_j\|^{-1}_{L^\be},&j=k\\
[.2cm]0,&j\ne k.
\end{array}\right.
$$
Clearly, \rf{htr} holds,
$$
\|\{\a_j\}\|_{L^\be(\ell^2)}
\le\left(\sum_jc_j\right)^{1/2}<\be,\quad
\|\{\g_k\}\|_{L^\be(\ell^2)}
\le\left(\sum_kc_k\right)^{1/2}<\be
$$
and
$$
\big\|\{\b_{jk}(x_2)\}_{j,k\ge0}\big\|_{{\mB}}\le1.
$$

In \cite{JTT} multiple operator integrals were defined for functions in the Haagerup tensor product of $L^\be$ spaces. Suppose that $\Psi$ has a representation of the form \rf{htr} and \rf{ogr} holds and suppose that $T$ and $R$ are bounded linear operators on Hilbert space. Then the triple operator integral 
\bay
\label{troinoi}
\iiint\Psi(x_1,x_2,x_3)\,dE_1(x_1)T\,dE_2(x_2)R\,dE_3(x_3)
\ey
can be defined in the following way.

Consider the spectral measure $E_2$. It is defined on a $\s$-algebra $\Sigma$ of subsets of
$\X_2$.
We can represent our Hilbert space $\h$ as the direct integral
\bay
\label{din}
\h=\int\limits_{\X_2}\bigoplus\sG(x)\,d\mu(x),
\ey
associated with $E_2$.
Here $\mu$ is a finite measure on $\X_2$, $x\mapsto\sG(x)$, is a measurable Hilbert family. The Hilbert space $\h$ consists of measurable functions $f$ such that $f(x)\in\sG(x)$, $x\in\X_2$, and
$$
\|f\|_\h\df\left(\;\int\limits_{\X_2}\|\f(x)\|_{\sG(x)}^2\,d\mu(x)\right)^{1/2}<\be.
$$
Finally, for $\D\in\Sigma$, $E(\D)$ is multiplication by the characteristic function of $\D$. We refer the reader to \cite{BS0}, Ch. 7 for an introduction to direct integrals of Hilbert spaces.

Suppose that $\Psi$ belongs to the Haagerup tensor product 
$L^\be\!\otimes_{\rm h}\!L^\be\!\otimes_{\rm h}\!L^\be$ and  
\rf{htr} holds. The triple operator integral \rf{troinoi} is defined by
\begin{align}
\label{htraz}
\iiint\Psi&(x_1,x_2,x_3)\,dE_1(x_1)T\,dE_2(x_2)R\,dE_3(x_3)\nonumber\\[.2cm]
=&
\sum_{j,k\ge0}\left(\int\a_j\,dE_1\right)T\left(\int\b_{jk}\,dE_2\right)
R\left(\int\g_k\,dE_3\right)\nonumber\\[.2cm]
=&\lim_{M,N\to\be}~\sum_{j=0}^N\sum_{k=0}^M
\left(\int\a_j\,dE_1\right)T\left(\int\b_{jk}\,dE_2\right)
R\left(\int\g_k\,dE_3\right).
\end{align}

Let us show that the series on the right converges in the weak operator topology.
Let $f$ and $g$ be vectors in $\h$. Put
\bay
\label{ukvj}
u_k\df R\left(\int\g_k\,dE_3\right)f\quad\mbox{and}\quad
v_j\df T^*\left(\int\ov{\a_j}\,dE_1\right)g.
\ey
We consider the vectors $v_j$ and $u_k$ as elements of the direct integral \rf{din}, i.e., vector functions on $\X_2$.

We have
\begin{align*}
&\left|\left(\sum_{j,k\ge0}\left(\int\a_j\,dE_1\right)T\left(\int\b_{jk}\,dE_2\right)
R\left(\int\g_k\,dE_3\right)f,g
\right)\right|\\[.2cm]
&=\left|\sum_{j,k\ge0}\left(\left(\int\b_{jk}\,dE_2\right)u_k,v_j\right)\right|
=\left|\sum_{j,k\ge0}\int\limits_{\X_2}\big(\b_{jk}(x)u_k(x),v_j(x)\big)_{\sG(x)}
\,d\mu(x)\right|
\\[.2cm]
&\le\int\limits_{\X_2}
\|\{\b_{jk}(x)\}_{j,k\ge0}\|_{\mB}\cdot\|\{u_k(x)\}_{k\ge0}\|_{\ell^2}
\cdot\|\{v_j(x)\}_{j\ge0}\|_{\ell^2}d\mu(x)
\\[.2cm]
&\le\|\{\b_{jk}\}_{j,k\ge0}\|_{L^\be(\mB)}
\left(~\int\limits_{\X_2}\Big(\sum_{k\ge0}|u_k(x)|^2\Big)d\mu(x)\right)^{1/2}\!
\!\left(~\int\limits_{\X_2}\Big(\sum_{j\ge0}|v_j(x)|^2\Big)d\mu(x)\right)^{1/2}
\\[.2cm]
&=\|\{\b_{jk}\}_{j,k\ge0}\|_{L^\be(\mB)}
\left(\sum_{k\ge0}\|u_k\|^2_\h\right)^{1/2}
\left(\sum_{j\ge0}\|v_j\|^2_\h\right)^{1/2}.
\end{align*}
Keeping \rf{ukvj} in mind, we see that the last expression is equal to
\begin{align*}
&\|\{\b_{jk}\}_{j,k\ge0}\|_{L^\be(\mB)}
\left(\sum_{k\ge0}\left\|R\left(\int\g_k\,dE_3\right)f\right\|^2_\h\right)^{1/2}
\left(\sum_{j\ge0}\left\|T^*\left(\int\ov{\a_j}\,dE_1\right)g\right\|^2_\h\right)^{1/2}
\\[.2cm]
&\le\|\{\b_{jk}\}_{j,k\ge0}\|_{L^\be\!(\mB)}\|R\|\!\cdot\!\|T\|\!
\left(\sum_{k\ge0}\left\|\left(\int\g_k\,dE_3\right)\!f\right\|^2\right)^{1/2}\!\!\!\!
\left(\sum_{j\ge0}\left\|\left(\int\ov{\a_j}\,dE_1\right)\!g\right\|^2\right)^{1/2}\!\!.
\end{align*}
By properties of integrals with respect to spectral measures,
$$
\sum_{k\ge0}\left\|\left(\int\g_k\,dE_3\right)f\right\|^2=
\left(\int\left(\sum_{k\ge0}|\g_k|^2\right)(dE_3f,f)\right)
\le\|\{\g_k\}_{k\ge0}\|^2_{L^\be(\ell^2)}\|f\|^2.
$$
Similarly,
$$
\sum_{j\ge0}\left\|\left(\int\ov{\a_j}\,dE_1\right)g\right\|^2=
\left(\int\left(\sum_{j\ge0}|\a_j|^2\right)(dE_1g,g)\right)
\le\|\{\a_j\}_{j\ge0}\|^2_{L^\be(\ell^2)}\|g\|^2.
$$
This implies that 
\begin{align*}
&\left|\left(\sum_{j,k\ge0}\left(\int\a_j\,dE_1\right)T\left(\int\b_{jk}\,dE_2\right)
R\left(\int\g_k\,dE_3\right)f,g
\right)\right|\\[.2cm]
&\le
\|\{\b_{jk}\}_{j,k\ge0}\|_{L^\be(\mB)}\cdot
\|\{\a_j\}_{k\ge0}\|_{L^\be(\ell^2)}\cdot
\|\{\g_k\}_{k\ge0}\|_{L^\be(\ell^2)}\|f\|\cdot\|g\|.
\end{align*}
It follows that the series in \rf{htraz} converges absolutely in the weak operator topology.


The above inequalities show that
\bay
\label{opno}
\!\!\!\!\left\|\iiint\Psi(x_1,x_2,x_3)\,dE_1(x_1)T\,dE_2(x_2)R\,dE_3(x_3)\right\|
\le\|\Psi\|_{L^\be\!\otimes_{\rm h}\!L^\be\!\otimes_{\rm h}\!L^\be}
\|T\|\cdot\|R\|.
\ey

Note that the triple operator integral is well defined by \rf{htraz}, i.e.,
the sum of the series in \rf{htraz} does not depend on the choice of a representation \rf{htr}, see \cite{JTT} and \cite{ANP3}.

It is easy to verify that if $\Psi$ is a function that belongs to the projective tensor product 
$L^\be(E_1)\hat\otimes L^\be(E_2)\hat\otimes L^\be(E_3)$, then the above definition 
coincides with the definition of the triple operator integral given in Chapter
2.

It turns out, however, that unlike in the case when the integrand belongs to the projective tensor product $L^\be\hat\otimes L^\be\hat\otimes L^\be$ (see Theorem \ref{Schpr}), triple operator integrals with integrands in the Haagerup tensor product 
$L^\be\!\otimes_{\rm h}\!L^\be\!\otimes_{\rm h}\!L^\be$ {\it do not possess the property} 
$$
T\in\bS_p,\quad R\quad\mbox{is bounded}\quad\Longrightarrow 
\iiint\Psi\,dE_1T\,dE_2R\,dE_3\in \bS_p
$$
with $p<2$; this will be established in \S~\ref{2c}. We will see in \S~\ref{SvNtoi} that for integrands $\Psi$ in $L^\be\!\otimes_{\rm h}\!L^\be\!\otimes_{\rm h}\!L^\be$, 
$$
T\in\bS_p,\quad R\in \bS_q,\quad \frac1p+\frac1q\le\frac12\quad\Longrightarrow 
\iiint\Psi\,dE_1T\,dE_2R\,dE_3\in \bS_r,\quad\frac1r=\frac1p+\frac1q.
$$
We do not know whether this can be true if $1/p+1/q>1/2$.

\

\section{Schatten--von Neumann properties}
\label{SvNtoi}

\

In this section we study Schatten--von Nemann properties of triple operator integrals with integrands in the Haagerup tensor product 
$L^\be\!\otimes_{\rm h}\!L^\be\!\otimes_{\rm h}\!L^\be$.
First, we consider the case when one of the operators is bounded and the other one belongs to the Hilbert--Schmidt class. Then we use an interpolation theorem for bilinear operators to a considerably more general situation.

The following result was established in \cite{ANP1} and its detailed proof was published in \cite{ANP3}.

\begin{thm}
\label{bS2}
Let $E_1$, $E_2$, and $E_3$ be spectral measures on Hilbert space and let $\Phi$
be a function in the Haagerup tensor product 
$L^\be(E_1)\!\otimes_{\rm h}\!L^\be(E_2)\!\otimes_{\rm h}\!L^\be(E_3)$.
Suppose that $T$ is a bounded linear operator and $R$ is an operator that belongs 
to the Hilbert--Schmidt class $\bS_2$. Then
\bay
\label{W}
W\df\int\limits_{\X_1}\int\limits_{\X_2}\int\limits_{\X_3}
\Psi(x_1,x_2,x_3)\,dE_1(x_1)T\,dE_2(x_2)R\,dE_3(x_3)\in\bS_2
\ey
and 
\bay
\label{WS_2}
\|W\|_{\bS_2}\le
\|\Psi\|_{L^\be\!\otimes_{\rm h}\!L^\be\!\otimes_{\rm h}\!L^\be}
\|T\|\cdot\|R\|_{\bS_2}.
\ey
\end{thm}

It is easy to see that Theorem \ref{bS2} implies the following fact:

\begin{cor}
\label{S2b}
Let $E_1$, $E_2$, $E_3$, and $\Psi$ satisfy the hypotheses of Theorem 
{\em\ref{bS2}}. If $T$ is a Hilbert Schmidt operator and $R$ is a bounded linear operator, then the operator $W$ defined by {\em\rf{W}} belongs to $\bS_2$ and
$$
\|W\|_{\bS_2}\le
\|\Psi\|_{L^\be\!\otimes_{\rm h}\!L^\be\!\otimes_{\rm h}\!L^\be}
\|T\|_{\bS_2}\|R\|.
$$
\end{cor}

Clearly, to deduce Corollary \ref{S2b} from Theorem \ref{bS2}, it suffices to consider the adjoint operator $W^*$.

\medskip

{\bf Proof of Theorem \ref{bS2}.} For simplicity we consider the case when $E_3$ is a discrete spectral measure and we refer the reader to \cite{ANP3} for the general case. Under this assumption, there exists an orthonormal basis 
$\{e_m\}_{m\ge0}$, the spectral measure $E_3$ is defined on the $\s$-algebra of all subsets
of $\Z_+$, and 
$E_3(\{m\})$ is the orthogonal projection onto the one-dimensional space spanned by $e_m$.
In this case the function $\Psi$ has the form
$$
\Psi(x_1,x_2,m)=\sum_{j,k\ge0}\a_j(x_1)\b_{jk}(x_2)\g_k(m),\quad
x_1\in\X_1,~x_2\in\X_2,~m\in\Z_+,
$$
where 
$$
\{\a_j\}_{j\ge0}\in L_{E_1}^\be(\ell^2), \quad 
\{\b_{jk}\}_{j,k\ge0}\in L_{E_2}^\be({\mathcal B}),
$$
and
$$
\sup_{m\ge1}\sum_{k\ge0}|\g_k(m)|^2<\be.
$$
Then
$$
W=\iint\sum_{m\ge0}\Psi(x_1,x_2,n)\,dE_1(x_1)T\,dE_2(x_2)R\,(\cdot,e_m)e_m.
$$
We have
\bay
\label{Zn}
\sum_{m\ge0}\|We_m\|^2=\sum_{m\ge0}\|Z_mRe_m\|^2,
\ey
where
\begin{align*}
Z_m&\df\iint\Psi(x_1,x_2,m)\,dE_1(x_1)T\,dE_2(x_2)\\[.2cm]
&=\iiint\Psi_m(x_1,x_2,m)\,dE_1(x_1)T\,dE_2(x_2)I\,d\E_m.
\end{align*}
Here $\E_m$ is the spectral measure defined on the one point set $\{m\}$
and the function $\Psi_m$ is defined on $\X_1\times\X_2\times\{m\}$ by
$$
\Psi_m(x_1,x_2,m)=\Psi(x_1,x_2,m),\quad x_1\in\X_1,x_2\in\X_2.
$$

It is easy to see that 
$$
\|\Psi_m\|_{L^\be(E_1)\otimes_{\rm h}L^\be(E_2)\otimes_{\rm h}L^\be(\E_m)}
\le\|\Psi\|_{L^\be(E_1)\otimes_{\rm h}L^\be(E_2)\otimes_{\rm h}L^\be(E_3)},
\quad m\ge0.
$$
It follows now from \rf{opno} that 
$$
\|Z_m\|\le\|\Psi\|_{L^\be\otimes_{\rm h}L^\be\otimes_{\rm h}L^\be}\|T\|,
$$
and by \rf{Zn}, we obtain
\begin{align*}
\sum_{m\ge0}\|We_m\|^2&\le\sum_{n\ge0}\|Z_m\|^2\|Re_m\|^2\\[.2cm]
&\le\|\Psi\|^2_{L^\be\otimes_{\rm h}L^\be\otimes_{\rm h}L^\be}
\|T\|^2\sum_{m\ge0}\|Re_m\|^2
\\[.2cm]
&=\|\Psi\|^2_{L^\be\otimes_{\rm h}L^\be\otimes_{\rm h}L^\be}\|T\|^2\|R\|_{\bS_2}^2.
\end{align*}
It follows that $W\in\bS_2$ and inequality \rf{WS_2} holds. $\bl$

We are going to use Theorem 4.4.1 from \cite{BL} on complex interpolation of bilinear operators. Recall that the Schatten--von Neumann classes $\bS_p$, $p\ge1$, and the space of bounded linear operators ${\mathcal B}(\h)$ form a complex interpolation scale: 
\bay
\label{Spint}
(\bS_1,{\mathcal B}(\h))_{[\theta]}=\bS_{\frac1{1-\theta}},\quad1<\theta<1.
\ey 
This fact is well known. For example, it follows from Theorem 13.1 of Chapter III of \cite{GK}.

The following result was established in \cite{ANP1} and its proof was published in \cite{ANP3}.

\begin{thm}
\label{SNSp}
Let $\Psi\in L^\be(E_1)\!\otimes_{\rm h}\!L^\be(E_2)\!\otimes_{\rm h}\!L^\be(E_3)$.
Then the following holds:
\newline
{\em(i)} if $p\ge2$, $T\in\B(\h)$, and $R\in\bS_p$, then the triple operator integral in {\em\rf{W}} belongs to $\bS_p$ and
\bay
\label{boSp}
\|W\|_{\bS_p}\le\|\Psi\|_{L^\be\!\otimes_{\rm h}\!L^\be\!\otimes_{\rm h}\!L^\be}
\|T\|\cdot\|R\|_{\bS_p};
\ey
\newline
{\em(ii)} if $p\ge2$, $T\in\bS_p$, and $R\in\B(\h)$, then the triple operator integral in {\em\rf{W}} belongs to $\bS_p$ and
$$
\|W\|_{\bS_p}\le\|\Psi\|_{L^\be\!\otimes_{\rm h}\!L^\be\!\otimes_{\rm h}\!L^\be}
\|T\|_{\bS_p}\|R\|;
$$
\newline
{\em(iii)} if $1/p+1/q\le1/2$, $T\in\bS_p$, and $R\in\bS_q$, 
then the triple operator integral in {\em\rf{W}} belongs to $\bS_r$
with $1/r=1/p+1/q$ and
$$
\|W\|_{\bS_r}\le\|\Psi\|_{L^\be\!\otimes_{\rm h}\!L^\be\!\otimes_{\rm h}\!L^\be}
\|T\|_{\bS_p}\|R\|_{\bS_q}.
$$
\end{thm}

\medskip

We will see in \S~\ref{2c} that neither (i) nor (ii) holds for $p>2$.

\medskip

{\bf Proof of Theorem \ref{SNSp}.} Let us first prove (i). Clearly, to deduce (ii) from (i), it suffices to consider $W^*$. 

Consider the bilinear operator $\W$ defined by
$$
\W(T,R)=\iiint\Psi(x_1,x_2,x_3)\,dE_1(x_1)T\,dE_2(x_2)R\,dE_3(x_3).
$$
By \rf{opno}, $\W$ maps $\mB(\h)\times\mB(\h)$ into $\mB(\h)$ and
$$
\|\W(T,R)\|\le\|T\|\cdot\|R\|.
$$
On the other hand, by Theorem \ref{bS2}, $\W$ maps $\mB(\h)\times\bS_2$
into $\bS_2$ and
$$
\|\W(T,R)\|_{\bS_2}\le\|T\|\cdot\|R\|_{\bS_2}.
$$
It follows from the complex interpolation theorem for linear operators 
(see \cite{BL}, Theorem 4.1.2 that) $\W$ maps $\mB(\h)\times\bS_p$, $p\ge2$,
into $\bS_p$ and
$$
\|\W(T,R)\|_{\bS_p}\le\|T\|\cdot\|R\|_{\bS_p}.
$$

Suppose now that $1/p+1/q\le1/2$ and $1/r=1/p+1/q$. It follows from statements (i) and (ii) (which we have already proved) that $\W$ maps $\mB(\h)\times\bS_r$ into
$\bS_r$ and $\bS_r\times\mB(\h)$ into $\bS_r$, and
$$
\|\W(T,R)\|_{\bS_r}\le\|T\|\cdot\|R\|_{\bS_r}
\quad\mbox{and}\quad
\|\W(T,R)\|_{\bS_r}\le\|T\|_{\bS_r}\cdot\|R\|.
$$
It follows from  Theorem 4.4.1 of \cite{BL} on interpolation of bilinear operators,
$\W$ maps $(\mB(\h),\bS_r)_{[\theta]}\times(\bS_r,\mB(\h))_{[\theta]}$ into
$\bS_r$ and 
$$
\|\W(T,R)\|_{\bS_r}\le
\|T\|_{(\mB(\h),\bS_r)_{[\theta]}}\|R\|_{(\bS_r,\mB(\h))_{[\theta]}}.
$$
It remains to observe that for $\theta=r/p$,
$$
(\mB(\h),\bS_r)_{[\theta]}=\bS_p\quad\mbox{and}\quad
(\bS_r,\mB(\h))_{[\theta]}=\bS_q,
$$
which is a consequence of \rf{Spint}. $\bl$

\

\section{Haagerup-like tensor products and triple operator integrals}
\label{Ttoi}

\

As we have mentioned in the introduction to this chapter, we are going to use a representation of $f(A_1,B_1)-f(A_2,B_2)$
in terms of triple operator integrals that involve the divided differences
$\dg^{[1]}f$ and $\dg^{[2]}f$. However, we will see in \S~\ref{2c} that the divided differences $\dg^{[1]}f$ and $\dg^{[2]}f$ do not have to belong to the
Haagerup tensor product \lb$L^\be\!\otimes_{\rm h}\!L^\be\!\otimes_{\rm h}\!L^\be$ for an arbitrary function $f$ in the Besov class $B_{\be,1}^1(\R^2)$. In addition to this, representation \rf{predst} involves operators of class $\bS_p$ with $p\le2$. However, we will see in \S~\ref{2c} that statements (i) and (ii) of Theorem \ref{SNSp} do not hold for $p<2$.

This means that we need a new approach to triple operator integrals. In this section we introduce Haagerup-like tensor products and define triple operator integrals whose integrands belong to such Haagerup-like tensor products. Note that the Haagerup-like tensor products were defined in \cite{ANP1} and \cite{AP10}, see also \cite{ANP3}.

\medskip

{\bf Definition 1.} 
A function $\Psi$ is said to belong to the Haagerup-like tensor product 
$L^\be(E_1)\!\otimes_{\rm h}\!L^\be(E_2)\!\otimes^{\rm h}\!L^\be(E_3)$ of the first kind if it admits a representation
\bay
\label{yaH}
\Psi(x_1,x_2,x_3)=\sum_{j,k\ge0}\a_j(x_1)\b_{k}(x_2)\g_{jk}(x_3),\quad x_j\in\X_j,
\ey
with $\{\a_j\}_{j\ge0},~\{\b_k\}_{k\ge0}\in L^\be(\ell^2)$ and 
$\{\g_{jk}\}_{j,k\ge0}\in L^\be(\mB)$. As usual, 
$$
\|\Psi\|_{L^\be\otimes_{\rm h}\!L^\be\otimes^{\rm h}\!L^\be}
\df\inf\big\|\{\a_j\}_{j\ge0}\big\|_{L^\be(\ell^2)}
\big\|\{\b_k\}_{k\ge0}\big\|_{L^\be(\ell^2)}
\big\|\{\g_{jk}\}_{j,k\ge0}\big\|_{L^\be(\mB)},
$$
the infimum being taken over all representations of the form \rf{yaH}.

\medskip

{\bf Definition 2.}
Let $1\le p\le2$. For 
$\Psi\in L^\be(E_1)\!\otimes_{\rm h}\!L^\be(E_2)\!\otimes^{\rm h}\!L^\be(E_3)$, for a bounded linear operator $R$, and for an operator $T$ of class $\bS_p$, we define the triple operator integral
\bay
\label{WHft}
W=\iint\!\!\upint\Psi(x_1,x_2,x_3)\,dE_1(x_1)T\,dE_2(x_2)R\,dE_3(x_3)
\ey
as the following continuous linear functional on $\bS_{p'}$,
$1/p+1/p'=1$ (on the class of compact operators in the case $p=1$):
\bay
\label{fko}
Q\mapsto
\trace\left(\left(
\iiint
\Psi(x_1,x_2,x_3)\,dE_2(x_2)R\,dE_3(x_3)Q\,dE_1(x_1)
\right)T\right).
\ey

\medskip

Clearly, the triple operator integral in \rf{fko} is well defined because the function
$$
(x_2,x_3,x_1)\mapsto\Psi(x_1,x_2,x_3)
$$ 
belongs to the Haagerup tensor product 
$L^\be(E_2)\!\otimes_{\rm h}\!L^\be(E_3)\!\otimes_{\rm h}\!L^\be(E_1)$. It follows easily from statement (i) of Theorem \ref{SNSp} that
$$
\|W\|_{\bS_p}\le\|\Psi\|_{L^\be\otimes_{\rm h}\!L^\be\otimes^{\rm h}\!L^\be}
\|T\|_{\bS_p}\|R\|,\quad1\le p\le2,
$$
(see Theorem \ref{ftHtp}).

It is easy to see that in the case when $\Psi$ belongs to the projective tensor product $L^\be(E_1)\hat\otimes L^\be(E_2)\hat\otimes L^\be(E_3)$, the definition of the triple operator integral given above is consistent with the definition of the triple operator integral given in Chapter 2. Indeed, it suffices to verify this for functions $\Psi$ of the form
$$
\Psi(x_1,x_2,x_3)=\f(x_1)\psi(x_2)\chi(x_3),\quad\f\in L^\be(E_1),\quad
\psi\in L^\be(E_2),\quad\chi\in L^\be(E_3),
$$
in which case the verification is obvious.

We also need triple operator integrals in the case when $T$ is a bounded linear operator and $R\in\bS_p$, $1\le p\le2$.

\medskip

{\bf Definition 3.} A function $\Psi$ is said to belong to the Haagerup-like tensor product $L^\be(E_1)\!\otimes^{\rm h}\!L^\be(E_2)\!\otimes_{\rm h}\!L^\be(E_3)$
of the second kind if
$\Psi$ admits a representation
\bay
\label{preds}
\Psi(x_1,x_2,x_3)=\sum_{j,k\ge0}\a_{jk}(x_1)\b_{j}(x_2)\g_k(x_3)
\ey
where $\{\b_j\}_{j\ge0},~\{\g_k\}_{k\ge0}\in L^\be(\ell^2)$, 
$\{\a_{jk}\}_{j,k\ge0}\in L^\be(\mB)$. The norm of $\Psi$ in 
the space $L^\be\otimes^{\rm h}\!L^\be\otimes_{\rm h}\!L^\be$ is defined by
$$
\|\Psi\|_{L^\be\otimes^{\rm h}\!L^\be\otimes_{\rm h}\!L^\be}
\df\inf\big\|\{\a_j\}_{j\ge0}\big\|_{L^\be(\ell^2)}
\big\|\{\b_k\}_{k\ge0}\big\|_{L^\be(\ell^2)}
\big\|\{\g_{jk}\}_{j,k\ge0}\big\|_{L^\be(\mB)},
$$
the infimum being taken over all representations of the form \rf{preds}.

\medskip

{\bf Definition 4.}
Suppose now that 
$\Psi\in L^\be(E_1)\!\otimes^{\rm h}\!L^\be(E_2)\!\otimes_{\rm h}\!L^\be(E_3)$,
$T$ is a bounded linear operator, and $R\in\bS_p$, $1\le p\le2$. The continuous linear functional 
$$
Q\mapsto
\trace\left(\left(
\iiint\Psi(x_1,x_2,x_3)\,dE_3(x_3)Q\,dE_1(x_1)T\,dE_2(x_2)
\right)R\right)
$$
on the class $\bS_{p'}$ (on the class of compact operators in the case $p=1$) 
determines an operator $W$ of class $\bS_p$, which
we call the triple operator integral
\bay
\label{WHst}
W=\upint\!\!\!\iint\Psi(x_1,x_2,x_3)\,dE_1(x_1)T\,dE_2(x_2)R\,dE_3(x_3).
\ey

\medskip

Moreover,
$$
\|W\|_{\bS_p}\le
\|\Psi\|_{L^\be\otimes^{\rm h}\!L^\be\otimes_{\rm h}\!L^\be}
\|T\|\cdot\|R\|_{\bS_p}.
$$

As above, in the case when 
$\Psi\in L^\be(E_1)\hat\otimes L^\be(E_2)\hat\otimes L^\be(E_3)$, the definition of the triple operator integral given above is consistent with the definition of the triple operator integral given in Chapter 2.

The following result can easily be deduced from 
Theorem \ref{SNSp}, see \cite{ANP3}.

\begin{thm}
\label{ftHtp}
Let $\Psi\in L^\be\!\otimes_{\rm h}\!L^\be\!\otimes^{\rm h}\!L^\be$.
Suppose that $T\in\bS_p$ and $R\in\bS_q$, where
$1\le p\le2$ and $1/p+1/q\le1$. Then the operator $W$ in {\em\rf{WHft}} belongs to $\bS_r$, $1/r=1/p+1/q$, and
\bay
\label{rpq}
\|W\|_{\bS_r}\le\|\Psi\|_{L^\be\otimes_{\rm h}\!L^\be\otimes^{\rm h}\!L^\be}
\|T\|_{\bS_p}\|R\|_{\bS_q}.
\ey
If $T\in\bS_p$, $1\le p\le2$, and $R$ is a bounded linear operator, then $W\in\bS_p$
and
\bay
\label{pB}
\|W\|_{\bS_p}\le\|\Psi\|_{L^\be\otimes_{\rm h}\!L^\be\otimes^{\rm h}\!L^\be}
\|T\|_{\bS_p}\|R\|.
\ey
\end{thm}

In the same way we can prove the following theorem:

\begin{thm}
\label{stHtp}
Let $\Psi\in L^\be\!\otimes^{\rm h}\!L^\be\!\otimes_{\rm h}\!L^\be$.
Suppose that $p\ge1$, $1\le q\le2$,  and $1/p+1/q\le1$. If $T\in\bS_p$, $R\in\bS_q$, then the operator $W$ in {\em\rf{WHst}} belongs to $\bS_r$, $1/r=1/p+1/q$, and
$$
\|W\|_{\bS_r}\le\|\Psi\|_{L^\be\otimes_{\rm h}\!L^\be\otimes^{\rm h}\!L^\be}
\|T\|_{\bS_p}\|R\|_{\bS_q}.
$$
If $T$ is a bounded linear operator and $R\in\bS_p$, $1\le p\le2$, then $W\in\bS_p$
and
$$
\|W\|_{\bS_p}\le\|\Psi\|_{L^\be\otimes_{\rm h}\!L^\be\otimes^{\rm h}\!L^\be}
\|T\|_{\bS_p}\|R\|.
$$
\end{thm}

\

\section{Conditions for ${\dg^{[1]}f}$ and 
${\dg^{[2]}}f$ to be in Haagerup-like tensor products}
\label{ddiff}

\

As we have already mentioned before, for functions $f$ in $B_{\be,1}^1(\R^2)$,
the divided differences $\dg^{[1]}f$ and $\dg^{[2]}f$,
$$
\big(\dg^{[1]}f\big)(x_1,x_2,y)\df\frac{f(x_1,y)-f(x_2,y)}{x_1-x_2}
\quad\mbox{and}\quad
\big(\dg^{[2]}f\big)(x,y_1,y_2)\df\frac{f(x,y_1)-f(x,y_2)}{y_1-y_2},
$$
do not have to belong to the Haagerup tensor product $L^\be\!\otimes_{\rm h}\!L^\be\!\otimes_{\rm h}\!L^\be$. This will be seen in \S~\ref{2c}.

In this section we will see that for $f\in B_{\be,1}^1(\R^2)$, the divided difference
$\dg^{[1]}f$
belongs to the tensor product 
$L^\be(E_1)\!\otimes_{\rm h}\!L^\be(E_2)\!\otimes^{\rm h}\!L^\be(E_3)$,
while the divided difference 
$\dg^{[2]}f$
belongs to the tensor product 
$L^\be(E_1)\!\otimes^{\rm h}\!L^\be(E_2)\!\otimes_{\rm h}\!L^\be(E_3)$
for arbitrary Borel spectral measures $E_1$, $E_2$, and $E_3$ on $\R$.

This will allow us to prove in the next section that if $(A_1,B_1)$ and $(A_2,B_2)$ are pairs of self-adjoint operators on Hilbert space, $(A_2,B_2)$ is an $\bS_p$ perturbation of $(A_1,B_1)$, $1\le p\le2$, and $f\in B_{\be,1}^1(\R^2)$, then the following integral formula holds:
\begin{align*}
f(A_1,B_1)-f(A_2,B_2)=&
\iint\!\!\upint\frac{f(x_1,y)-f(x_2,y)}{x_1-x_2}
\,dE_{A_1}(x_1)(A_1-A_2)\,dE_{A_2}(x_2)\,dE_{B_1}(y),\\[.2cm]
+&\upint\!\!\!\iint\frac{f(x,y_1)-f(x,y_2)}{y_1-y_2}
\,dE_{A_2}(x)\,dE_{B_1}(y_1)(B_1-B_2)\,dE_{B_2}(y_2).
\end{align*}
 
The following theorem contains a formula that is crucial for our estimates. It was established in \cite{ANP1} and \cite{AP10}, its detailed proof was published in \cite{ANP3}.

\begin{thm}
\label{rrp}
Let $f$ be a bounded function on $\R^2$ whose Fourier transform is supported 
in the ball $\{\xi\in\R^2:~\|\xi\|\le1\}$. Then
\bay
\label{crfor}
\frac{f(x_1,y)-f(x_2,y)}{x_1-x_2}=
\sum_{j,k\in\Z}\frac{\sin(x_1-j\pi)}{x_1-j\pi}\cdot\frac{\sin(x_2-k\pi)}{x_2-k\pi}
\cdot\frac{f(j\pi,y)-f(k\pi,y)}{j\pi-k\pi},
\ey
where for $j=k$, we assume that 
$$
\frac{f(j\pi,y)-f(k\pi,y)}{j\pi-k\pi}
=\frac{\partial f(x,y)}{\partial x}\Big|_{(j\pi,y)}.
$$

Moreover,
\bay
\label{sinusy}
\sum_{j\in\Z}\frac{\sin^2(x_1-j\pi)}{(x_1-j\pi)^2}
=\sum_{k\in\Z}\frac{\sin^2(x_2-k\pi)}{(x_2-k\pi)^2}=1,
\quad x_1~x_2\in\R,
\ey
and
\bay
\label{comnor}
\sup_{y\in\R}\left\|\left\{\frac{f(j\pi,y)-f(k\pi,y)}{j\pi-k\pi}
\right\}_{j,k\in\Z}\right\|_{\mB}\le\const\|f\|_{L^\be(\R)}.
\ey
\end{thm}

\medskip

Formula \rf{crfor} can be deduced from Lemma \ref{razl}, see \cite{ANP3} for details, identities \rf{sinusy} are well-known, see, e.g., \cite{Ti}, 3.3.2, Example IV.

To estimate the operator norm of the matrix
$$
\left\{\frac{f(j\pi,y)-f(k\pi,y)}{j\pi-k\pi}
\right\}_{j,k\in\Z},
$$
we represent this matrix as the sum of the matrices $C_y=\{c_{jk}(y)\}_{j,k\in\Z}$
and $D_y=\{d_{jk}(y)\}_{j,k\in\Z}$, where
$$
c_{jk}(y)=\left\{\begin{array}{ll}\frac{f(j\pi,y)-f(k\pi,y)}{j\pi-k\pi},&j\ne k\\[.2cm]
0,&j=k
\end{array}\right.
\quad\mbox{and}\quad
d_{jk}(y)=\left\{\begin{array}{ll}0,&j\ne k\\[.2cm]
\frac{\partial f(x,y)}{\partial x}\Big|_{(j\pi,y)},&j=k.
\end{array}\right.
$$

It is easy to see that $C_y$ is the commutator of the discrete Hilbert transform ${\mathcal H}_{\rm d}$ and an operator of multiplication by a bounded sequence on 
$\ell^2$ and $\|C_y\|\le\const\|f\|_{L^\be(\R^2)}$.

On the other hand,
$$
\|D_y\|=\sup_{j\in\Z}\left|\frac{\partial f(x,y)}{\partial x}\Big|_{(j\pi,y)}\right|
\le\|f\|_{L^\be(\R^2)}
$$
by Bernstein's inequality. This completes the proof of \rf{comnor}.

We refer the reader to \cite{ANP3} for details.

The following result can be deduced easily from Theorem \ref{rrp}, see \cite{ANP3}.

\begin{cor}
\label{sle}
Let $f$ be a bounded function on $\R^2$ such that its Fourier transform is supported in $\{\xi\in\R^2:~\|\xi\|\le\s\}$, $\s>0$. Then the divided differences
$\dg^{[1]}f$ and $\dg^{[2]}f$ have the following properties: 
$$
\dg^{[1]}f\in L^\be(E_1)\!\otimes_{\rm h}\!L^\be(E_2)\!\otimes^{\rm h}\!L^\be(E_3)
\quad\mbox{and}\quad
\dg^{[2]}f\in L^\be(E_1)\!\otimes^{\rm h}\!L^\be(E_2)\!\otimes_{\rm h}\!L^\be(E_3)
$$
for arbitrary Borel spectral measures $E_1$, $E_2$ and $E_3$. Moreover,
\bay
\label{dg1}
\big\|\dg^{[1]}f\big\|_{L^\be\!\otimes_{\rm h}\!L^\be\!\otimes^{\rm h}\!L^\be}
\le\const\s\|f\|_{L^\be(\R^2)}
\ey
and
\bay
\label{dg2}
\big\|\dg^{[2]}f\big\|_{L^\be\!\otimes^{\rm h}\!L^\be\!\otimes_{\rm h}\!L^\be}
\le\const\s\|f\|_{L^\be(\R^2)}.
\ey
\end{cor}

Corollary \ref{sle} implies in turn the following theorem that was established in \cite{AP3}.

\begin{thm}
\label{Bes}
Let $f\in B_{\be,1}^1(\R^2)$. Then
$$
\dg^{[1]}f\in L^\be(E_1)\!\otimes_{\rm h}\!L^\be(E_2)\!\otimes^{\rm h}\!L^\be(E_3)
\quad\mbox{and}\quad
\dg^{[2]}f\in L^\be(E_1)\!\otimes^{\rm h}\!L^\be(E_2)\!\otimes_{\rm h}\!L^\be(E_3)
$$
for arbitrary Borel spectral measures $E_1$, $E_2$ and $E_3$. Moreover,
$$
\big\|\dg^{[1]}f\big\|_{L^\be\!\otimes_{\rm h}\!L^\be\!\otimes^{\rm h}\!L^\be}
\le\const\|f\|_{B_{\be,1}^1}
$$
and
$$
\big\|\dg^{[2]}f\big\|_{L^\be\!\otimes^{\rm h}\!L^\be\!\otimes_{\rm h}\!L^\be}
\le\const\s\|f\|_{B_{\be,1}^1}.
$$
\end{thm}

\

\section{Lipschitz type estimates in the case ${1\le p\le2}$}
\label{ple2}

\

In this section we discuss the results announced in \cite{ANP1} and \cite{ANP2}
whose detailed proofs were given in \cite{ANP3}.

We will see in this section that for functions $f$ in the Besov class 
$B_{\be,1}^1(\R^2)$, we have a Lipschitz type estimate for functions of noncommuting self-adjoint operators in the norm of $\bS_p$ with $p\in[1,2]$. 

The following integral formula plays an important role.

\begin{thm}
\label{intfor}
Let $f\in B_{\be,1}^1(\R^2)$ and $1\le p\le2$. Suppose that $(A_1,B_1)$ and $(A_2,B_2)$ are pairs of self-adjoint operators such that $A_2-A_1\in\bS_p$ and 
$B_2-B_1\in\bS_p$. Then the following identity holds:
\begin{align}
\label{osnf}
f(A_1,B_1)&-f(A_2,B_2)\nonumber\\[.2cm]
&=
\iint\!\!\upint\frac{f(x_1,y)-f(x_2,y)}{x_1-x_2}
\,dE_{A_1}(x_1)(A_1-A_2)\,dE_{A_2}(x_2)\,dE_{B_1}(y),\nonumber\\[.2cm]
&+\upint\!\!\!\iint\frac{f(x,y_1)-f(x,y_2)}{y_1-y_2}
\,dE_{A_2}(x)\,dE_{B_1}(y_1)(B_1-B_2)\,dE_{B_2}(y_2).
\end{align}
\end{thm}

Note that by Theorem \ref{Bes}, the divided differences $\dg^{[1]}f$ and $\dg^{[2]}f$ belong to the corresponding Haagerup like tensor products, and so the triple operator integrals on the right make sense.

\medskip

\Pf It suffices to prove that
\begin{align}
\label{perfor}
f(A_1,B_1)&-f(A_2,B_1)\nonumber\\[.2cm]
&=\iint\!\!\upint\big(\dg^{[1]}f\big)(x_1,x_2,y)
\,dE_{A_1}(x_1)(A_1-A_2)\,dE_{A_2}(x_2)\,dE_{B_1}(y)
\end{align}
and 
\begin{align}
\label{vtfor}
f(A_2,B_1)&-f(A_2,B_2)\nonumber\\[.2cm]
&=\upint\!\!\!\iint\big(\dg^{[2]}f\big)(x,y_1,y_2)
\,dE_{A_2}(x)\,dE_{B_1}(y_1)(B_1-B_2)\,dE_{B_2}(y_2).
\end{align}

Let us establish \rf{perfor}. Formula \rf{vtfor} can be proved in exactly the same way. 

Suppose first that the function $\dg^{[1]}f$ belongs to the projective tensor product \lb$L^\be(E_{A_1})\hat\otimes L^\be(E_{A_2})\hat\otimes L^\be(E_{B_1})$. 
In this case we can write
\begin{align*}
\iint\!\!\upint\big(\dg^{[1]}f\big)(x_1,x_2,y)&
\,dE_{A_1}(x_1)(A_1-A_2)\,dE_{A_2}(x_2)\,dE_{B_1}(y)\\[.2cm]
&=\iint\!\!\upint\big(\dg^{[1]}f\big)(x_1,x_2,y)
\,dE_{A_1}(x_1)A_1\,dE_{A_2}(x_2)\,dE_{B_1}(y)\\[.2cm]
&-\iint\!\!\upint\big(\dg^{[1]}f\big)(x_1,x_2,y)
\,dE_{A_1}(x_1)A_2\,dE_{A_2}(x_2)\,dE_{B_1}(y).
\end{align*}
Note that the above equality does not make sense if 
$\dg^{[1]}f$ does not belong to $L^\be\hat\otimes L^\be\hat\otimes L^\be$ because
the operators $A_1$ and $A_2$ do not have to be compact, while the definition of
triple operator integrals with integrands in the Haagerup-like tensor product
$L^\be\!\otimes_{\rm h}\!L^\be\!\otimes^{\rm h}\!L^\be$
assumes that the operators $A_1$ and $A_2$ belong to $\bS_2$. 

It follows immediately from the definition of triple operator integrals with integrands in $L^\be\hat\otimes L^\be\hat\otimes L^\be$ that
\begin{align*}
\iint\!\!\upint\big(\dg^{[1]}f\big)(x_1,x_2,y)&
\,dE_{A_1}(x_1)A_1\,dE_{A_2}(x_2)\,dE_{B_1}(y)\\[.2cm]
=&\iint\!\!\upint x_1\big(\dg^{[1]}f\big)(x_1,x_2,y)
\,dE_{A_1}(x_1)\,dE_{A_2}(x_2)\,dE_{B_1}(y)
\end{align*}
and
\begin{align*}
\iint\!\!\upint\big(\dg^{[1]}f\big)(x_1,x_2,y)&
\,dE_{A_1}(x_1)A_2\,dE_{A_2}(x_2)\,dE_{B_1}(y)\\[.2cm]
=&\iint\!\!\upint x_2\big(\dg^{[1]}f\big)(x_1,x_2,y)
\,dE_{A_1}(x_1)\,dE_{A_2}(x_2)\,dE_{B_1}(y).
\end{align*}
Thus
\begin{align*}
\iint\!\!\upint\big(\dg^{[1]}f\big)&(x_1,x_2,y)
\,dE_{A_1}(x_1)A_1\,dE_{A_2}(x_2)\,dE_{B_1}(y)\\[.2cm]
-&\iint\!\!\upint\big(\dg^{[1]}f\big)(x_1,x_2,y)
\,dE_{A_1}(x_1)A_2\,dE_{A_2}(x_2)\,dE_{B_1}(y)\\[.2cm]
=&\iint\!\!\upint(x_1-x_2)\frac{f(x_1,y)-f(x_2,y)}{x_1-x_2}
\,dE_{A_1}(x_1)\,dE_{A_2}(x_2)\,dE_{B_1}(y)\\[.2cm]
=&\iint\!\!\upint f(x_1,y)\,dE_{A_1}(x_1)\,dE_{A_2}(x_2)\,dE_{B_1}(y)\\[.2cm]
-&\iint\!\!\upint f(x_2,y)\,dE_{A_1}(x_1)\,dE_{A_2}(x_2)\,dE_{B_1}(y)
=f(A_1,B_1)-f(A_2,B_1).
\end{align*}

Consider the functions $f_n$ defined by $f_n=f*W_n$,
$n\in\Z$, see \rf{fn}. It is easy to see from the definition of the Besov class $B_{\be,1}^1(\R^2)$ that to prove \rf{perfor}, it suffices to show that
\begin{align*}
f_n(A_1,B_1)&-f_n(A_2,B_1)\\[.2cm]
&=\iint\!\!\upint\big(\dg^{[1]}f_n\big)(x_1,x_2,y)
\,dE_{A_1}(x_1)(A_1-A_2)\,dE_{A_2}(x_2)\,dE_{B_1}(y).
\end{align*}

As we have mentioned in \S~\ref{prel}, the function $f_n$ is a restriction of an entire function of two variables to 
$\R\times\R$. Thus it suffices to establish formula \rf{perfor} in the case when $f$ is an entire function. To complete the proof, we show that for entire functions $f$ the divided differences $\dg^{[1]}f$ must belong to the projective tensor product
$L^\be(E_{A_1})\hat\otimes L^\be(E_{A_2})\hat\otimes L^\be(E_{B_1})$. 

Let $f(x,y)=\sum\limits_{j=0}^\be\Big(\sum\limits_{k=0}^\be a_{jk}x^jy^k\Big)$ be an entire function and let $R$ be a positive number such that the spectra $\s(A_1)$, 
$\s(A_2)$, and $\s(B)$ are contained in $[-R/2,R/2]$.
Clearly,
$$
\|f\|_{L^\be\hat\otimes L^\be}\le\sum_{j=0}^\be\Big(\sum_{k=0}^\be |a_{jk}|R^{j+k}\Big)
<\be
$$
and
\begin{align*}
\left\|\dg^{[1]}f\right\|_{L^\be\hat\otimes L^\be\hat\otimes L^\be}&=
\left\|\sum_{j=0}^\be\left(\sum_{k=1}^\be
\Big(\sum_{l=0}^{j-1} a_{jk}x_1^lx_2^{j-1-l}y^k\Big)
\right)\right\|_{L^\be\hat\otimes L^\be\hat\otimes L^\be}\\[.2cm]
&\le\sum\limits_{j=0}^\be\left(\sum\limits_{k=1}^\be j|a_{jk}|R^{j+k-1}\right)<+\be,
\end{align*}
where in the above expressions $L^\be$ means $L^\be[-R,R]$. This completes the proof. $\bl$
 
\begin{thm}
\label{Lipp>2}
Let $p\in[1,2]$. Then there is a positive number $C$ such that
\bay
\label{Lipnerp2}
\|f(A_1,B_1)-f(A_2,B_2)\|\le C\|f\|_{B_{\be,1}^1}
\max\big\{\|A_1-A_2\|_{\bS_p},\|B_1-B_2\|_{\bS_p}\big\},
\ey
whenever $f\in B_{\be,1}^1(\R^2)$, and $A_1$, $A_2$, $B_1$, and $B_2$ are self-adjoint operators such that
$A_2-A_1\in\bS_p$ and 
$B_2-B_1\in\bS_p$.
\end{thm}

\Pf This is an immediate consequence of Theorem \ref{intfor} and Theorems
\ref{ftHtp} and \ref{stHtp}. $\bl$

\medskip

{\bf Remark.} We have defined functions $f(A,B)$ for $f$ in 
$B_{\be,1}^1(\R^2)$ only for bounded self-adjoint operators $A$ and $B$. 
However, formula \rf{osnf} allows us to define the difference $f(A_1,B_1)-f(A_2,B_2)$ in the case when $f\in B_{\be,1}^1(\R^2)$ and the 
self-adjoint operators $A_1,\,A_2,\,B_1,\,B_2$ are possibly unbounded once we know that the pair $(A_2,B_2)$ is an $\bS_p$ perturbation of the pair $(A_1,B_1)$, 
$1\le p\le2$.
Moreover, inequality \rf{Lipnerp2} also holds for such operators.

\medskip

Note that similar results for functions of unitary operators were obtained in
\cite{ANP3} as well.

\

\section{Lipschitz type estimates cannot be extended beyond ${p\le2}$}
\label{Bp>2}

\

It was showed in \cite{ANP2} and \cite{ANP3} that there 
is no Lipschitz type inequality of the form \rf{Lipp>2} in the norm of $\bS_p$ with $p>2$ and in the operator norm for an arbitrary function $f$ in $B_{\be,1}^1(\R^2)$.
In this section we give the construction of \cite{ANP3}. 

\begin{thm} 
\label{nLte}
{\em(i)}
There is no positive number $M$ such that 
$$
\|f(A_1,B)-f(A_2,B)\|\le M\|f\|_{L^\be(\R^2)}\|A_1-A_2\|
$$ 
for all bounded functions $f$ on $\R^2$ with Fourier transform supported 
in $[-2\pi,2\pi]^2$ and for all finite rank self-adjoint operators $A_1,\,A_2,\,B$.

{\em(ii)}
Let $p>2$.
Then there is no positive number $M$ such that 
$$
\|f(A_1,B)-f(A_2,B)\|_{\bS_p}\le M\|f\|_{L^\be(\R^2)}\|A_1-A_2\|_{\bS_p}
$$ 
for all bounded functions $f$ on $\R^2$ with Fourier transform supported in $[-2\pi,2\pi]^2$ and for all finite rank self-adjoint operators $A_1,\,A_2,\,B$.
\end{thm}

\Pf Let us first prove (ii). Let $\{g_j\}_{1\le j\le N}$ and $\{h_j\}_{1\le j\le N}$
be orthonormal systems in Hilbert space. 
Consider the rank one projections $P_j$ and $Q_j$ defined by
$$
P_j v=(v,g_j)g_j\quad\mbox{and}\quad Q_jv=(v,h_j)h_j,\quad 1\le j\le N.
$$
We define the self-adjoint operators $A_1$, $A_2$, and $B$ by
$$
A_1=\sum_{j=1}^N2jP_j,\quad A_2=\sum_{j=1}^N(2j+1) P_j,\quad\mbox{and}\quad
B=\sum_{k=1}^Nk\, Q_k.
$$
Then $\|A_1-A_2\|_{\bS_p}=N^{\frac1p}$.
Put
$$
\f(x)=\frac{1-\cos2\pi x}{2\pi^2 x^2}.
$$
Clearly, 
$\supp\F\f\subset[-2\pi,2\pi]$, $\f(k)=0$ for all $k\in\Z$ such that $k\ne0$, $\f(0)=1$.
Put $\f_k(x)=\f(x-k)$.
Given a matrix $\{\tau_{jk}\}_{1\le j,k\le N}$, we define 
the function $f$ by
$$
f(x,y)=\sum_{1\le j,k\le N}\tau_{jk}\f_{2j}(x)\f_k(y).
$$ 
It is easy to see that $\f_{2j}(A_1)=P_j$, $\f_{2j}(A_2)=0$, $\f_k(B)=Q_k$ provided
$1\le j,k\le N$, and
$$
\|f\|_{L^\be(\R^2)}\le\const\max_{1\le j,k\le N}|\tau_{jk}|.
$$
Clearly,
$$
f(A_1,B)=
\sum_{1\le j,k\le N}\tau_{jk}P_jQ_k\quad\mbox{and}\quad
f(A_2,B)=\0.
$$
Note that 
$$
(f(A_1,B)h_k,g_j)=\tau_{jk}(h_k,g_j),\quad1\le j,k\le N.
$$  
Clearly, for every unitary matrix
$\{u_{jk}\}_{1\le j,k\le N}$, there exist orthonormal systems $\{g_j\}_{1\le j\le N}$ and $\{h_j\}_{1\le j\le N}$
such that  $(h_k,g_j)=u_{jk}$.
Put 
$$
u_{jk}\df\frac1{\sqrt N}\exp\left(\frac{2\pi{\rm i}jk}N\right),\quad 1\le j,k\le N.
$$
Obviously, $\{u_{jk}\}_{1\le j,k\le N}$ is a unitary matrix.
Hence, 
we may find vectors $\{g_j\}_{j=1}^N$ and $\{h_j\}_{j=1}^N$ such that $(h_k,g_j)=u_{jk}$.
Put $\tau_{jk}=\sqrt N \,\,\ov u_{jk}$. Then
$$
\|f(A_1,B)\|_{\bS_p}=\|\{|u_{jk}|\}_{1\le j,k\le N}\|_{\bS_p}=\|\{|u_{jk}|\}_{1\le j,k\le N}\|_{\bS_2}=\sqrt N
$$
because $\rank\{|u_{jk}|\}_{1\le j,k\le N}=1$.
So for each positive integer $N$ we have constructed a function $f$ and operators $A_1$, $A_2$, $B$ 
such that $|f|\le\const$, $\supp\F f\subset[-2\pi,2\pi]^2$, $\|A_1-A_2\|_{\bS_p}=N^{\frac1p}$ and
$\|f(A_1,B)-f(A_2,B)\|_{\bS_p}=\sqrt N$. It remains to observe that $\lim_{N\to\be}N^{\frac12-\frac1p}=\be$
for $p>2$. 

Exactly the same construction works to prove (i). It suffices to replace in the above construction the $\bS_p$ norm with the operator norm and observe that 
$\|A_1-A_2\|=1$ and $\|f(A_1,B)-f(A_2,B)\|=\sqrt N$. $\bl$

Theorem \ref{nLte} implies that there is no Lipschitz type estimate in the operator norm and in the $S_p$ norm with $p>2$. Note that in the construction given in the proof the norms of $A_1-A_2$ cannot get small. The following result shows that we can easily overcome this problem. 

\begin{thm}
\label{BnLte}
There exist a sequence $\{f_n\}_{n\ge0}$ of functions in $B_{\be,1}^1(\R^2)$ and sequences of self-adjoint finite rank operators $\big\{A_1^{(n)}\big\}_{n\ge0}$, 
$\big\{A_2^{(n)}\big\}_{n\ge0}$, 
and $\big\{B^{(n)}\big\}_{n\ge0}$
such that the norms $\|f_n\|_{B_{\be,1}^1}$ do not depend on $n$,
$$
\lim_{n\to\be}\big\|A_1^{(n)}-A_2^{(n)}\big\|\to0,\quad\mbox{but}\quad
\|f_n(A_1,B)-f_n(A_2,B)\|\to\be.
$$
The same is true in the norm of $\bS_p$ for $p>2$.
\end{thm}

\Pf The existence of such sequences can be obtained easily from the construction in the proof of 
Theorem \ref{nLte}. It suffices to make the following observation. Let $f$, $A_1$, $A_2$ and $B$ be as in the proof of Theorem \ref{nLte} and let $\e>0$. Put 
$f_\e(x,y)\df\e f\big(\frac{x}{\e},\frac{y}{\e}\big)$. Then 
$$
\|f_\e\|_{B_{\be 1}^1}=\|f\|_{B_{\be 1}^1}, 
\quad \|f_\e(\e A_1,\e B)-f_\e(\e A_2,\e B)\|=\e N^{1/2},
\quad\mbox{and}\quad \|\e A_1-\e A_2\|=\e.
$$
If $p>2$, then
$$
\|f_\e(\e A_1,\e B)-f_\e(\e A_2,\e B)\|_{\bS_p}=\e N^{1/2}
\quad\mbox{and}\quad \|\e A_1-\e A_2\|_{\bS_p}=\e N^{1/p}.\quad\bl
$$

\medskip

{\bf Remark.} The construction given in the proof of Theorem \ref{nLte} shows that for every
positive number $M$ there exist a function $f$ on $\R^2$ whose Fourier transform is supported in 
$[-2\pi,2\pi]^2$ such that $\|f\|_{L^\be(\R)}\le\const$ and self-adjoint operators of finite rank $A_1$, $A_2$, $B$ such that $\|A_1-A_2\|=1$, but $\|f(A_1,B)-f(A_2,B)\|>M$.
It follows that unlike in the case of commuting self-adjoint operators (see \cite{APPS}), the fact that $f$ is a H\"older function of order $\a\in(0,1)$ on $\R^2$ does not imply the H\"older type estimate 
$$
\|f(A_1,B_1)-f(A_2,B_2)\|\le\const\max\big\{\|A_1-A_2\|^\a,\|B_1-B_2\|^\a\big\}.
$$

\

\section{Counterexamples}
\label{2c}

\

We use the results of the previous section to show that statements (i) and (ii) of Theorem \ref{SNSp} do not hold for 
$p\in[1,2)$. We also deduce from the results of \S~\ref{Bp>2} that the divided differences $\dg^{[1]}f$ and $\dg^{[2]}f$ do not have to belong to
the Haagerup tensor product $L^\be\!\otimes_{\rm h}\!L^\be\!\otimes_{\rm h}\!L^\be$ for an arbitrary function $f$ in $B_{\be,1}^1(\R^2)$. The results of this sections were obtained in \cite{AP1}, \cite{AP2} and \cite{AP3}.

\begin{thm}
\label{kpp<2}
Let $1\le p<2$.
There exist an operator $Q$ in $\bS_p$, 
spectral measures $E_1$, $E_2$ and $E_3$ on Borel subsets of $\R$ and a function $\Phi$ in the Haagerup tensor product 
$L^\be(E_1)\!\otimes_{\rm h}\!L^\be(E_2)\!\otimes_{\rm h}\!L^\be(E_3)$ and an operator $Q$ in $\bS_p$ such that
$$
\iiint\Phi(x_1,x_2,x_2)\,dE_1(x_1)\,dE_2(x_2)Q\,dE_3(x_3)\not\in\bS_p.
$$
\end{thm}

\Pf Assume the contrary. Then the linear operator
$$
Q\mapsto\iiint\Phi(x_1,x_2,x_2)\,dE_1(x_1)\,dE_2(x_2)Q\,dE_3(x_3)
$$
is bounded on $\bS_p$ for arbitrary Borel spectral measures $E_1$, $E_2$, and $E_3$
and for an arbitrary function $\Phi$ in 
$L^\be(E_1)\!\otimes_{\rm h}\!L^\be(E_2)\!\otimes_{\rm h}\!L^\be(E_3)$. 
Suppose now that $\Psi$ belongs to the Haagerup-like tensor product
$L^\be(E_1)\!\otimes_{\rm h}\!L^\be(E_2)\!\otimes^{\rm h}\!L^\be(E_3)$
of the first kind. For a finite rank operator $T$ consider the triple operator integral
$$
W=\iint\!\!\upint\Psi(x_1,x_2,x_3)\,dE_1(x_1)T\,dE_2(x_2)\,dE_3(x_3).
$$
We define the function $\Phi$ defined by
$$
\Phi(x_2,x_3,x_1)=\Psi(x_1,x_2,x_3).
$$
Let $Q\in\bS_p$. We have
\begin{align*}
\trace(WQ)&=
\trace\left(\left(
\iiint
\Psi(x_1,x_2,x_3)\,dE_2(x_2)\,dE_3(x_3)Q\,dE_1(x_1)
\right)T\right)\\[.2cm]
&=\trace\left(\left(
\iiint
\Phi(x_2,x_3,x_1)\,dE_2(x_2)\,dE_3(x_3)Q\,dE_1(x_1)
\right)T\right)
\end{align*}
(see the definition of triple operator integrals with integrands in the 
Haagerup-like tensor product of the first kind in \S~\ref{Ttoi}).

Thus
\begin{align*}
|\trace(WQ)|&=
\left|\trace\left(\left(
\iiint
\Phi(x_2,x_3,x_1)\,dE_2(x_2)\,dE_3(x_3)Q\,dE_1(x_1)
\right)T\right)\right|\\[.2cm]
&\le
\left\|\left(
\iiint
\Phi(x_2,x_3,x_1)\,dE_2(x_2)\,dE_3(x_3)Q\,dE_1(x_1)
\right)\right\|_{\bS_p}\|T\|_{\bS_{p'}}\\[.2cm]
&\le\|\Phi\|_{L^\be\!\otimes_{\rm h}\!L^\be\!\otimes_{\rm h}\!L^\be}
\|Q\|_{\bS_p}\|T\|_{\bS_{p'}}
\end{align*}
(throughout the proof of this theorem in the case $p=1$, the norm in $\bS_{p'}$ has to be replaced with the operator norm).

It follows that
\begin{align}
\label{WSp'}
\|W\|_{\bS_{p'}}&=
\left\|
\iint\!\!\upint\Psi(x_1,x_2,x_3)\,dE_1(x_1)T\,dE_2(x_2)\,dE_3(x_3)
\right\|_{\bS_{p'}}\nonumber\\[.2cm]
&\le\|\Psi\|_{L^\be\!\otimes_{\rm h}\!L^\be\!\otimes^{\rm h}\!L^\be}\|T\|_{\bS_{p'}}.
\end{align}

By Theorem \ref{Bes}, 
$\dg^{[1]}f\in L^\be\!\otimes_{\rm h}\!L^\be\!\otimes^{\rm h}\!L^\be$ for every
$f$ in $B_{\be,1}^1(\R^2)$ and by \rf{perfor},
\begin{align*}
f(A_1,B)&-f(A_2,B)\nonumber\\[.2cm]
&=\iint\!\!\upint\big(\dg^{[1]}f\big)(x_1,x_2,y)
\,dE_{A_1}(x_1)(A_1-A_2)\,dE_{A_2}(x_2)\,dE_B(y)
\end{align*}
for arbitrary finite rank self-adjoint operators $A_1$, $A_2$, and $B$.
it remains to observe that by inequality \rf{WSp'},
\begin{align*}
\|f(A_1,B_1)-f(A_2,B)\|_{\bS_{p'}}
&\le
\|\dg^{[1]}f\|_{L^\be\!\otimes_{\rm h}\!L^\be\!\otimes^{\rm h}\!L^\be}
\|A_1-A_2\|_{\bS_{p'}}\\[.2cm]
&\le\const\|f\|_{B_{\be,1}^1}\|A_1-A_2\|_{\bS_{p'}}
\end{align*}
which contradicts Theorem \ref{BnLte}. $\bl$

If we pass to the adjoint operator, we can see that for $p\in[1,2)$, there exist a function $\Psi$ in the Haagerup tensor product 
$L^\be\!\otimes_{\rm h}\!L^\be\!\otimes_{\rm h}\!L^\be$ and an operator $Q$ in $\bS_p$ such that
$$
\iiint\Phi(x_1,x_2,x_2)\,dE_1(x_1)Q\,dE_2(x_2)\,dE_3(x_3)\not\in\bS_p.
$$

The following application of Theorem \ref{BnLte} shows that for functions
$f$ in $B_{\be,1}^1(\R^2)$, the divided differences $\dg^{[1]}f$ and $\dg^{[2]}f$
do not have to belong to the Haagerup tensor product 
$L^\be\!\otimes_{\rm h}\!L^\be\!\otimes_{\rm h}\!L^\be$. We state the result for
$\dg^{[1]}f$. 

\begin{thm}
\label{Bnd12}
There exists a function $f$ in the Besov class $B_{\be,1}^1(\R^2)$ such that the divided difference $\dg^{[1]}f$ does not belong to 
$L^\be\!\otimes_{\rm h}\!L^\be\!\otimes_{\rm h}\!L^\be$.
\end{thm}

\Pf Assume the contrary. Then the map
$$
f\mapsto\dg^{[1]}f
$$
is a bounded linear operator from $B_{\be,1}^1(\R^2)$ to 
$L^\be\!\otimes_{\rm h}\!L^\be\!\otimes_{\rm h}\!L^\be$.

By \rf{perfor},
\begin{align*}
f(A_1,B)&-f(A_2,B)\nonumber\\[.2cm]
&=\iint\!\!\upint\big(\dg^{[1]}f\big)(x_1,x_2,y)
\,dE_{A_1}(x_1)(A_1-A_2)\,dE_{A_2}(x_2)\,dE_B(y)
\end{align*}
for arbitrary finite rank self-adjoint operators $A_1$, $A_2$, and $B$.
It follows now from inequality \rf{opno} that
$$
\|f(A_1,B)-f(A_2,B)\|
\le\big\|\dg^{[1]}f\big\|_{L^\be\!\otimes_{\rm h}\!L^\be\!\otimes_{\rm h}\!L^\be}
\|A_1-A_2\|
\le\const\|f\|_{B_{\be,1}^1}\|A_1-A_2\|
$$
which contradicts Theorem \ref{BnLte}. $\bl$

\

\section{Functions of almost commuting operators, an extension of the Helton--Howe trace formula}
\label{HHf}

\

Operators $A$ and $B$ are called {\it almost commuting} if the commutator
$$
[A,B]\df AB-BA
$$
belongs to $\bS_1$. In \cite{HH} the following trace formula was discovered:
\bay
\label{HeHo}
\trace\big({\rm i}\big[\big(\f(A,B),\psi(A,B)\big]\big)=\frac1{2\pi}
\iint_{\R^2}\left(\frac{\partial\f}{\partial x}\frac{\partial\psi}{\partial y}-
\frac{\partial\f}{\partial y}\frac{\partial\psi}{\partial x}\right)g(x,y)\,dx\,dy,
\ey
where $A$ and $B$ are almost commuting self-adjoint operators, $\f$ and $\psi$ are polynomials and $g$ is the Pincus principal function which is uniquely determined by $A$ and $B$ and which was introduced in \cite{Pin}.

The problem considered in \cite{Pe} was to extend the Helton--Howe trace formula for a reasonably big class of functions. It was shown in \cite{Pe} that under natural assumptions it is impossible to extend formula \rf{HeHo} to the class of all continuously differentiable functions. On the other hand, a sufficiently big class of functions 
$\mC$ was such that formula \rf{HeHo} holds for all functions $\f$ and $\psi$ in
$\mC$.

In the paper \cite{AP10} it was proved that formula \rf{HeHo} admits an extension to arbitrary functions $\f$ and $\psi$ in the Besov class $B_{\be,1}^1(\R^2)$ which considerably improved the sufficient condition $\f,\,\psi\in\mC$ found in \cite{Pe}.
The main tools used in \cite{AP10} are Haagerup-like tensor products and triple operator integrals.

The following results were obtained in \cite{AP10}.

\begin{thm}
\label{komut}
Let $A$ and $B$ be self-adjoint operators and let $Q$ be a bounded linear operator such that $[A,Q]\in\bS_1$ and $[B,Q]\in\bS_1$. 
Suppose that $\f\in B_{\be,1}^1(\R^2)$.
Then $[\f(A,B),Q\big]\in\bS_1$,
\begin{align*}
\big[\f(A,B),Q\big]&=
\iint\!\!\upint\frac{\f(x_1,y)-\f(x_2,y)}{x_1-x_2}\,dE_A(x_1)[A,Q]\,dE_A(x_2)\,dE_B(y)
\nonumber
\\[.2cm]
&+\upint\!\!\!\iint\frac{\f(x,y_1)-\f(x,y_2)}{y_1-y_2}\,dE_A(x)\,dE_B(y_1)[B,Q]\,dE_B(y_2)
\end{align*}
and
$$
\big\|[\f(A,B),Q\big]\big\|_{\bS_1}
\le\const\|\f\|_{B_{\be,1}^1(\R^2)}\big(\big\|[A,Q]\big\|_{\bS_1}+
\big\|[B,Q]\big\|_{\bS_1}\big).
$$
\end{thm}

If we apply Theorem \ref{komut} to the operator $Q=\psi(A,B)$, we obtain the following result:

\begin{thm}
\label{glav}
Let $A$ and $B$ be almost commuting self-adjoint operators and let $\f$ and $\psi$ be functions in the Besov class $B_{\be,1}^1(\R^2)$. Then
\begin{align*}
\big[\f(A,B),\psi(A,B)\big]&=
\iint\!\!\upint
\frac{\f(x_1,y)-\f(x_2,y)}{x_1-x_2}\,dE_A(x_1)[A,\psi(A,B)]\,dE_A(x_2)\,dE_B(y)
\nonumber
\\[.2cm]
&+
\upint\!\!\!\iint
\frac{\f(x,y_1)-\f(x,y_2)}{y_1-y_2}\,dE_A(x)\,dE_B(y_1)[B,\psi(A,B)]\,dE_B(y_2)
\end{align*}
and
$$
\big\|[\f(A,B),\psi(A,B)\big]\big\|_{\bS_1}
\le\const\|\f\|_{B_{\be,1}^1(\R^2)}\|\psi\|_{B_{\be,1}^1(\R^2)}
\big\|[A,B]\big\|_{\bS_1}.
$$
\end{thm}

Theorem \ref{glav} allows us to extend the Helton--Howe trace formula to functions in the Besov class $B_{\be,1}^1(\R^2)$.

\begin{thm}
Let $A$ and $B$ be almost commuting self-adjoint operators and let $\f$ and $\psi$ be functions in the Besov class $B_{\be,1}^1(\R^2)$. Then the following formula holds:
$$
\trace\big({\rm i}\big[\big(\f(A,B),\psi(A,B)\big]\big)=\frac{1}{2\pi}
\iint_{\R^2}\left(\frac{\partial\f}{\partial x}\frac{\partial\psi}{\partial y}-
\frac{\partial\f}{\partial y}\frac{\partial\psi}{\partial x}\right)g(x,y)\,dx\,dy,
$$
where $g$ is the Pincus principal function associated with the operators 
$A$ and $B$.
\end{thm}

We refer the reader to \cite{AP10} for more detail.

\

\

\noindent
\begin{tabular}{p{5cm}p{4.5cm}p{4.6cm}}
V.V. Peller \\
Department of Mathematics  \\
Michigan State University\\
East Lansing, Michigan 48824 \\
USA\\
\end{tabular}

\end{document}